\documentclass[12pt,english]{article}
\usepackage[T1]{fontenc}
\usepackage[latin9]{inputenc}
\usepackage{url}
\usepackage{amsthm}
\usepackage{amsmath}
\usepackage{amssymb}
\usepackage{graphicx}
\usepackage{setspace}
\usepackage{esint}
\usepackage[bookmarks=false]{hyperref}
\makeatletter

\let\SF@@footnote\footnote
\def\footnote{\ifx\protect\@typeset@protect
    \expandafter\SF@@footnote
  \else
    \expandafter\SF@gobble@opt
  \fi
}
\expandafter\def\csname SF@gobble@opt \endcsname{\@ifnextchar[
  \SF@gobble@twobracket
  \@gobble
}
\edef\SF@gobble@opt{\noexpand\protect
  \expandafter\noexpand\csname SF@gobble@opt \endcsname}
\def\SF@gobble@twobracket[#1]#2{}

  \theoremstyle{plain}
  \newtheorem{assumption}{\protect\assumptionname}
\theoremstyle{plain}
\newtheorem{thm}{\protect\theoremname}
  \theoremstyle{plain}
  \newtheorem{lem}{\protect\lemmaname}
  \theoremstyle{plain}
  \newtheorem{prop}{\protect\propositionname}
  \theoremstyle{plain}
  \newtheorem{cor}{\protect\corollaryname}
\newenvironment{lyxlist}[1]
{\begin{list}{}
{\settowidth{\labelwidth}{#1}
 \setlength{\leftmargin}{\labelwidth}
 \addtolength{\leftmargin}{\labelsep}
 }}
{\end{list}}

\usepackage{fullpage}
\usepackage{titlesec}
\usepackage{lscape}
\usepackage{graphicx}
\usepackage{algpseudocode}
\usepackage{subfig}
\date{April 15, 2017}

\titleformat{\section}{\large\sc\center}{\thesection}{1em}{}

\allowdisplaybreaks[4]

\@ifundefined{showcaptionsetup}{}{%
 \PassOptionsToPackage{caption=false}{subfig}}
\usepackage{subfig}
\makeatother

\usepackage{babel}
  \providecommand{\assumptionname}{Assumption}
  \providecommand{\lemmaname}{Lemma}
  \providecommand{\propositionname}{Proposition}
\providecommand{\corollaryname}{Corollary}
\providecommand{\theoremname}{Theorem}

\begin{document}

\title{Uncertainty in Economic Growth and Inequality%
\thanks{I would so much like to thank those whom I know or I haven't known
yet supporting me in this exploration.%
}}

\author{Zhengyuan Gao%
\thanks{Center for Operations Research and Econometrics (CORE), Universit�
catholique de Louvain, Voie du Roman Pays 34, B-1348, Louvain-la-Neuve,
Belgium. E-mail: \texttt{zhengyuan.gao@uclouvain.be}%
}}
\maketitle
\begin{abstract}
A step to consilience, starting with a deconstruction of the causality
of uncertainty that is embedded in the fundamentals of growth and
inequality, following a construction of aggregation laws that disclose
the invariance principle across heterogeneous individuals, ending
with a reconstruction of metric models that yields deeper structural
connections via U.S. GDP and income data.
\end{abstract}

\section{Motivation}

With a simple belief about facts that must be, Laozi in the 6th BCE
deems the phenomena that are wholly in harmony behaving in a completely
natural and uncontrived way. He says, ``Nature does nothing, but
everything is done.''%
\footnote{Tao Te Ching - Chapter 37.%
} A similar economic ideology, \emph{laissez faire} used by the physiocrats
in the 17th century or later \emph{invisible hand} used by Adam Smith
on income distribution and production, has been deeply rooted in the
core of modern economic analysis. In this perspective, individuals
optimize decisions, markets rationalize deals, different forces stabilize
interferences, and so equilibria are established, as an outcome of
these self-organized dynamics, which in principle should lead us to
a well-acceptable status.

However when one turns to another perspective viewing these economic
elements as an integrated, evolving and continuous entity, one may
lose oneself in a chaotic, intricate, and unordered retrospective
of economic history. Economic system is consistently vibrating. Individuals
are incapable of seizing their own fates, markets tend to generate
various bubbles, political forces weave a gravity field making things
fallen. The commonly accepted equilibria seem to be in an unreachable
state. The dynamical process, instead of driving us to a utopia, constructs
a predicament of mankind concealing the problems origins. Joseph Schumpeter
in \cite{Schumpeter1950} writes the following words, ``The essential
point to grasp is that in dealing with capitalism we are dealing an
evolutionary process. It may seem strange that anyone can fail to
see so obvious a fact which moreover was long ago emphasized by Karl
Marx. {[}...{]} Capitalism, then, is by nature a form or method of
economic change and not only never is but never can be stationary.'' 

One ostensible reason of this inconsistency is that models and their
assumptions are merely abstractions or oversimplifications of the
reality, thus it would be natural that they fail to consider some
of these phenomena. It would be even more natural to attribute these
unexplained imperfections to a category caused by uncertainty and
then let it explain them, since by definition uncertainty is uncertain
and is hard to be merged into a system that is certain. But what if
there is no inconsistency? What if everything is in order and one
unpleasant outcome in an economy is just a consequence of a sequence
of pleasant decisions well-accepted amongst individuals? What if in
this sequence uncertainty will certainly bring in the denouement? 

A crux is the role of uncertainty. Current treatment on uncertainty
is to separate it from the fundamentals of economic theory while expecting
it to conciliate the contradictory theories about how to interpret
the existing phenomena. The gap between two perspectives, as long
as it does not affect the self-consistencies in either domain, is
considered to be harmless. Nevertheless, there may exist another perspective:
all things return to uncertainty. It could be the uncertainty that
manipulates the economy, that shapes the social order, that generates
and wipes out our expectations. In such a case, uncertainty is no
longer the bridge over the gap, it eliminates the gap. This faith
is simple: diverse individual decisions reflected by their own beliefs
about the laws of this uncertain world, after rationalization, integration
and evolution, converge to a field that is characterized by a certain
law, an invariant of the individuals'. Uncertainty, therefore, is
certain.

Return to the controversial fundamentals. Two basic opposites co-exist
along the economic history: growth and inequality. About their causal
relationship, \cite{Kuznets1955} asks the following question ``Does
inequality in the distribution of income increase or decrease in the
course of a country's economic growth?'' He did not answer this dualism
question in \cite{Kuznets1955}, nor elsewhere to my best knowledge.
The dilemma is that the obvious answer seems to contradict the righteousness
in theory while in the rigid dichotomy the alternative answer is too
ambiguous to be believed as a law. In the history of economy, when
new products are introduced there is an intense amount of research
and development which leads to dramatic improvements in quality and
reductions in cost. This leads to a period of rapid economic growth.
However, the owners of these new products, at the time of their occurrences,
stand in sharp contrast to the great mass of their contemporaries.
Every advance first comes into being innovation developed by one or
few persons, only to become, after a time, the indispensable necessity
taken for granted by the majority. Equality does not seem to exist
at the very beginning of the growth. The opportunities of quality
improvement and cost reduction are gradually exhausted, then the products
are in widespread use and the trend turns into alleviating such an
inequality but also growth. Stimulus in production caused by the inequality
is considered as an abstract source, full of uncertainties. Thus entanglement
of growth and inequality, like other controversial fundamentals, masks
the causality behind complex stochastic patterns. 

Economic growth and inequality, from the most causal matters of divergences
in incomes and wages to the profoundest concentration of capital,
is from my perspective a man-made configuration. Growth is a human
right of surviving and propagating while equality is a human desire
rooting deeply in the conscious and belief. How would the invisible
hands dispose the relationship of these two basic elements? Enlightened
by theories and methodology in \cite{Simon1955}, this paper considers
uncertainty as a composition of more elementary stochastic components
available amongst individuals then formulate concrete dynamics using
laws extracted from individuals. Deconstruction of uncertainty makes
it feasible to epistemologically discuss previous questions. Abstraction
of inequality as a limit aggregation law of economic growth concretizes
its representation. All these are done with the preservation of two
basic elements, growth and equality, over individuals. Uncertainty
of growth entails the laws by which inequality is invariantly revealed.
The latter is so determined by its intrinsic equality conditions that
there is few choice left to an attempt of wiping it out, even at the
fundamental level.

\section{Guide}

Section \ref{sec:A-non-abstract-model} is a preliminary where the
connection with mainstream growth model is established but the substance
emphasizes the use of new assumptions. Section \ref{sec:Abstract-Models}
extends the set of new assumptions to a more general setting in order
to incorporate with uncertainty. Equilibrium is shown to be in existence
and its law is unique for each individual. Individuals are endowed
with stationary probabilistic laws for the equilibrium growth. Section
\ref{sec:Aggregation} shifts the attention to aggregates. Infinite
divisible law is introduced. Individual's law becomes aggregatable.
Non-stationary aggregates are distinguished from stationary aggregates.
The latter ones follow the same law as individuals' while the former
ones need some additional characterization. Section \ref{sec:Structural-Aggregation}
exploits structural relationship between individual and aggregate
laws. Aggregation of laws is represented by summing and scaling of
structural parameters. These parameters reveal possible latent movements.
Some of them conflict with superficial facts. Section \ref{sec:USdata}
uses data to illustrate some theoretical arguments. Some forward-looking
thoughts are refined in Section \ref{sec:Bird's-eye-View}.  Proofs,
if necessary, are given in the appendix. Some proofs are chosen to
be heuristic so that they are accessible to readers with general mathematical
background. Each main section is followed by general remarks summarizing
its main theme.

\section{A non-abstract model\label{sec:A-non-abstract-model}}

It is sometimes said that societies have to choose between greater
equality and economic growth. This section shows an ideal growth model
that can achieve the equal status regardless of initial heterogeneities
amongst different individuals. It reveals that growth does not necessarily
have negative correlation with equality. Once we separate the causality
between growth and equality, we examine different objectives of growth
and equality and see how these differences consequently consolidate
and become coherent in the equilibrium. All these arguments rely on
three assumptions.

\subsection{Scarcity, growth, and equality}

Let $X_{t}(\omega^{i})$ denote an economic variable of individual
$i$ at time $t$ with an initial state $X_{0}^{i}=\omega^{i}$, such
as a capital intensity. Suppose $\omega^{i}\in\Omega$ and $X_{t}(\omega^{i})\in\Omega$
for any $t$ and $i$. For the simplest case, we assume that the set
$\Omega$ is on a closed and bounded interval. When there is no necessity
for emphasizing a specific individual, sometime $X_{t}(\omega)$ may
refer to any individual.
\begin{assumption}
\label{assu:scarcity_1}The set $\Omega$ is on a closed and bounded
interval in $\mathbb{R}$. Without loss of generality, we assume $\Omega=[0,1]$.
\end{assumption}
This assumption reflects that each individual faces the scarcity of
economic resources. The assumption applies to both initial endowments
$\omega^{i}$ and future states $\{X_{t}(\omega^{i})\}_{t>0}$ of
this individual. It eliminates the possibility of unbounded growth
for any individual $i$.
\begin{assumption}
\label{assu:growth_1}Assume the growth function $f$ of this economic
variable as a real-valued monotone non-decreasing function $f(\cdot):\Omega\rightarrow\Omega$
where $X_{t}(\omega^{i})=f^{(t)}(\omega^{i})=f(f^{(t-1)}(\omega^{i}))$
for $t=1,2,\dots$. 
\end{assumption}
This assumption gives a non-decreasing growth function. With Assumption
\ref{assu:scarcity_1}, these two assumptions induce that any individual
in this economy will ultimately reach his or her equilibrium status.
Because $f(\cdot)$ is a monotone function on $[0,1]$, any sequence
of $f^{(t)}(\cdot)$ converges as $t\rightarrow\infty$.%
\footnote{Any monotone sequence on a bounded and closed set in $\mathbb{R}$
always converges.%
} Otherwise, either there exists a divergent sequence that will violate
the condition of monotonicity or the sequence $f^{(t)}$ will go beyond
the bounded interval $[0,1]$. 

Consider two households with initial endowment $0$ and $1$ respectively.
At time $t$, their economic variables become $X_{t}(0)=f^{(t)}(0)$
and $X_{t}(1)=f^{(t)}(1)$. The convergent property of $f$ entails
$X_{t}(0)\uparrow X^{*}(0)$ and $X_{t}(1)\downarrow X^{*}(1)$ where
$X^{*}(0)$ and $X^{*}(1)$ are the convergent limits and $X^{*}(0)\leq X^{*}(1)$. 
\begin{assumption}
\label{assu:equality_1}Two limits $X^{*}(0)$ and $X^{*}(1)$ of
the convergent sequences $X_{t}(0)$ and $X_{t}(1)$ are equivalent.
\end{assumption}
This assumption implies an ideology of equality, as the individual
with the smallest initial endowment is assumed to be able to achieve
the same equilibrium state as those with the biggest endowment $X^{*}(0)=X^{*}(1)$.
Because everyone in this economy has a convergent $X^{*}(\omega)$
and because any $X^{*}(\omega)$ should be not smaller than $X^{*}(0)$
and not larger than $X^{*}(1)$, we have $X^{*}(\omega)=X^{*}(0)=X^{*}(1)$.
The fact is that for any individual with any initial state $\omega\in[0,1]$,
when $t\rightarrow\infty$ such an individual will attain the same
limit $X^{*}(\omega)$ as the others. Thus $X^{*}(0)$ is the unique
fixed point of the limit of $\{X_{t}(\omega)\}_{t}$ for any $\omega\in\Omega$.
This result is summarized in the following theorem.
\begin{thm}
\label{thm:preliminary}Given Assumption \ref{assu:scarcity_1} to
\ref{assu:equality_1}, the economy reach an equal status, namely
any $X_{t}(\omega^{i})$ converges to the same $X^{*}$ . 
\end{thm}
These three assumptions play different roles. Assumption \ref{assu:scarcity_1}
allows for heterogenous initial states for all individuals and allows
us to identify their wealth levels especially those at the extremes.
Assumption \ref{assu:growth_1} imposes the same growth function $f(\cdot)$
to all individuals with different endowments. The monotone non-decreasing
$f(\cdot)$ implies that growth happens for all individuals. Assumption
\ref{assu:equality_1} is analogous to the statement ``all men are
created equal'', because it assumes that the poor and the rich will
come to the same state. The inequality exists at the very beginning
due to the heterogenous initial endowments, but with the same growth
function, eventually the inequality vanishes in this economy.

\subsection{A Growth Model for $N$ Individuals}

Although it looks simple, the previous specification can illustrate
some essences in framework of the classic Solow growth model. In Solow's
framework, the economic variable $X_{t}$ refers to capital per labor
such that $X_{t}=K_{t}/L_{t}$. Provided that the production function
$f_{P}(K_{t},L_{t})$ has constant returns to scale, one can set the
production per labor to $f_{P}(X_{t},1)$. The growth function refers
to the fundamental equation of the Solow model%
\footnote{The corresponding differential equation is $dX/dt=\beta_{1}f_{P}(X_{t},1)-\beta_{2}X_{t}$
which is derived from Cobb-Douglas production function. This equation
models the capital stock for an economy in which technology and the
supply of labor do not change.%
}
\[
X_{t+1}=f(X_{t})=\beta_{1}f_{P}(X_{t},1)-(\beta_{2}-1)X_{t}
\]
where $\beta_{1}$ is the saving rate and $\beta_{2}$ is the effective
depreciation rate of capital per labor. 

Given an arbitrary number of individuals $N$ in this economy, these
individuals provide labor and rent capital in a competitive labor
and capital market, they have access to the same neoclassical technology,
and produce a homogeneous output. For a monotone $f(\cdot)$ and $X_{t}=K_{t}/L_{t}$
in $\Omega$, there exists one $X^{*}$. If the poorest and the richest
can reach the same level of production in the long run, then Theorem
\ref{thm:preliminary} assures a unique fixed point for all initial
endowments $\omega^{i}\in\Omega$. Therefore for any individual, $f(X_{t})\rightarrow X^{*}$
as $t\rightarrow\infty$. In the equilibrium, the production function
becomes
\[
f(X^{*})=\beta_{1}A_{S}^{*}+\left(\beta_{1}B_{S}^{*}+1-\beta_{2}\right)X^{*}
\]
where $B_{S}^{*}$ and $A_{S}^{*}$ are equilibrium coefficients satisfying
the fixed point result $f(X^{*})=X^{*}$ which coincides with the
result in the classic Solow model.%
\footnote{The exact forms are $B_{S}^{*}=f_{P}^{'}(X^{*})$ and $A_{S}^{*}=f_{P}(X^{*})-f_{P}^{'}(X^{*})X^{*}$.
The result comes from the maximization of a quasi-linear utility function.
Illustration is given in Appendix \ref{sub:Solow-Model}. %
}The value $X^{*}$ is independent of the initial endowment $\omega$.
One can conclude that for any initial state $\omega\in[0,1]$, the
equilibrium $X^{*}$ makes all individuals access to the same production
type as they use exactly the same capital-labor ratio. 

Moreover, a linear aggregation is simplified as a summation of $N$
identical $X^{*}$s when $t\rightarrow\infty$. The growth function
\[
f\left(N\times X^{*}\right)=N\times X^{*}.
\]
induces the stability of the structure. It refers to the fact that
the aggregate growth acts as if it was a non aggregated function of
individuals. The aggregation process does not distort the individual
growth function in the limit.%
\footnote{We should distinguish it from the capital deepening situation in the
Solow model where population growth is represented by an increase
of $L$. In this case, the population increases and the capital intensity
$K/L$ decreases, so economic expansion will not continue indefinitely.
Normally, in the Solow model, capital deepening is considered as a
necessary but not a sufficient condition for the economic development.%
}

\subsection{Remarks}

Assumption \ref{assu:scarcity_1} to \ref{assu:equality_1} restrict
our attention to one simple model with $N$-individuals. In this model,
individuals share the same growth function but at the end they reach
the same status. This model implies that individual deterministic
growth itself may not generate inequality, on the contrary, it may
alleviate inequality. In this $N$-individuals economy, when an equilibrium
capital-labor ratio exists and its existence is independent of the
initial state $\omega$ for any individual, then eventually the productions
of individuals converge to the same level. This economy is homogenized
as the inequality together with heterogeneous growth features vanish
amongst the individuals. Another implication is that the linear aggregation
may not distort the growth if the deterministic economy stays in its
equilibrium.%
\footnote{In fact, the social welfare maximized by the equilibrium $X^{*}$
is equivalent to the sum of maximum individual utility that is the
same across all individuals. A social planner can choose $X^{*}$
to maximize the social welfare meanwhile all individuals also agree
with $X^{*}$ as it maximizes their utilities.\emph{ }Thus in this
model, the decentralized market economy and the centralized economy
are isomorphic.%
}

In reality, economic growth often accompanies with greater inequality.
Inequality may be enlarged by some relevant factors driven by the
growth as well as the initial heterogenous states. Given their unrealistic
implications, Assumption \ref{assu:scarcity_1} to \ref{assu:equality_1}
will \emph{not} be used afterwards. However, the enlightenment of
these assumptions is that without uncertainty some specification can
drive a growing economy with heterogeneous individuals to an equality
status. With this enlightenment in mind, we introduce uncertainty
to the economy in order to accommodate more general and realistic
dynamics. The new model accommodates the criteria of scarcity, growth
and equality in a wider sense.

\section{An Abstract Model\label{sec:Abstract-Models}}

There are several ways of introducing uncertainty to growth models.
One may think of multiple equilibria or unstable dynamics%
\footnote{For example, Lotka-Volterra type differential equations, known as
the predator-prey equations.%
} or statistical errors. These methods usually need some rigid conditions
for the aggregation and can only compile limited types of uncertainty.
Instead of adding uncertain features to the aggregates, we search
for probabilistic laws that relate to both individuals and aggregates.
Individuals are governed by these laws, meanwhile a similar law is
automatically attached to the aggregates. 

Such probabilistic laws are formed in an equilibrium economy.%
\footnote{The way of constructing these probabilistic laws as equilibria is
inspired by Volterra's classification for deterministic laws of physics
and de Finetti's infinite divisible laws of probability.%
} This section provides the result of this equilibrium economy. Three
new assumptions are proposed. The new assumptions are comparable with
the previous deterministic ones. The result says that the probabilistic
law of any individual growth is a stationary distribution that is
independent identical across all individuals in this equilibrium economy.

\subsection{Abstract Assumptions}

We consider an abstract dynamical economy that is a space of $N$-individuals
\[
\left((\Omega,\mathcal{B},(\mathcal{B}_{t}),\mathbb{P}),\mathcal{A},f\right)^{\otimes N}.
\]
We allow $N$ to grow so that we can represent a large economy with
possibly arbitrary number of individuals. The space $(\Omega,\mathcal{B},(\mathcal{B}_{t}),\mathbb{P})^{\otimes N}$
is an $N$-folds of a filtered probability space including a state
space $\Omega$, a filtration $\mathcal{B}_{1}\subseteq\mathcal{B}_{2}\subseteq\cdots$
such that $\lim_{t\rightarrow\infty}\mathcal{B}_{t}\rightarrow\mathcal{B}$,
and an unknown probability law $\mathbb{P}$ for each individual.
The set $\mathcal{A}^{\otimes N}$ is an $N$-folds arbitrary index
set for an $N$-folds growth function $f^{\otimes N}$. Each individual
$i$ has the growth function $f:\Omega\rightarrow\Omega$ indexed
by $a_{t}\in\mathcal{A}$ such as $f^{(a_{t})}(\omega^{i})$. In this
economy, a sequence of measurements of dynamical variables (vectors)
on $\Omega$ are collected over time. For example the first measurement
of the variable $X$ of individual $i$ is denoted as $X_{1}(\omega^{i})=f^{(a_{1})}(\omega^{i})$.
If an argument is true for any $i$, there is no necessity of emphasizing
$i$. In this case, we can drop $i$ and express $X_{t}(\omega)$
for a variable of an arbitrary individual at time $t$.
\begin{assumption}
\label{assu:scarcity_2} (Abstract scarcity) In the economy $\left((\Omega,\mathcal{B},(\mathcal{B}_{t}),\mathbb{P}),\mathcal{A},f\right)^{\otimes N}$,
each individual faces an identical set $\Omega$ of potential states.
Given a total ordered set $\mathcal{P}$ such that $\Omega\subset\mathcal{P}$.
we assume that $\mathcal{P}$ is a complete lattice so that the set
$\Omega$ has an infimum $\underline{\Omega}\subset\mathcal{P}$ and
a supremum $\overline{\Omega}\subset\mathcal{P}$.
\end{assumption}
Assumption \ref{assu:scarcity_2} generalizes Assumption \ref{assu:scarcity_1}.
Since the number of individuals is allowed to grow to infinity and
we assume that each individual faces the same set of potential states,
the topological structure of $\Omega$ for all possible states in
this economy needs to be extended. This previous setting $\Omega=[0,1]$
is a special case of the total order set with a complete lattice in
$\mathbb{R}$.%
\footnote{When the total ordered set $\mathcal{P}$ is a complete lattice, the
structure of $\mathcal{P}$ induces compactness for $\Omega\subset\mathcal{P}$.
However, the countable union of $\Omega$ may not be compact.%
} The total ordered structure $(\mathcal{P},\leq)$ assures that the
individual is able to compare and rank any two elements from $\Omega\subset\mathcal{P}$.
The complete lattice structure states that the individual face scarcity
in a wider sense. In the deterministic scarcity $\omega\in[0,1]$,
$\omega$ has to be a scalar value while $\omega\in\Omega$ in Assumption
\ref{assu:scarcity_2} is endowed with a general set theoretic structure.
This assumption implies that an unlimited want of economic elements
in $\Omega$ is impossible.
\begin{assumption}
\label{assu:growth_2} (Individual stochastic growth) The index sequence
$\mathbf{a}:=\{a_{t}\}$ is a sequence of i.i.d. random variables
$a_{t}$ taking values in the set $\mathcal{A}$. For any sequence
$\mathbf{a}\in\mathcal{A}$, there exists a subsequence $\{a_{t}^{'}\}$
such that function $f^{(a_{t}^{'})}(\cdot):\Omega\rightarrow\Omega$
is a monotone non-decreasing. For any $x$ and $y$ in $\Omega$,
the set $\{\mathbf{a}:f^{(\mathbf{a})}(x)\leq y\}$ is measurable
in $(\Omega,\mathcal{B},(\mathcal{B}_{t}),\mathbb{P})$. 
\end{assumption}
Assumption \ref{assu:growth_2} is a stochastic counterpart of Assumption
\ref{assu:growth_1}. In Assumption \ref{assu:growth_1}, every $\{f^{(a)}(\cdot)\}_{a\in\mathcal{A}}$
needs to be monotone while Assumption \ref{assu:growth_2} only requires
that every $\{f^{(a)}(\cdot)\}_{a\in\mathcal{A}}$ has at least one
monotone subsequence $\{f^{(a_{t}^{'})}(\cdot)\}_{t\in\mathbb{N}}$
and the index set $\{a_{t}^{'}\}_{t\in\mathbb{N}}$ of this subsequence
is a random variable that is measurable.

The consequence of Assumption \ref{assu:scarcity_1} and \ref{assu:growth_1}
is that any $X_{t}(\omega)=f^{(t)}(\omega)$ will converge to a limit
as $t\rightarrow\infty$ regardless of its initial endowment. This
consequence does not hold in the stochastic setting. However, Assumption
\ref{assu:scarcity_2} and \ref{assu:growth_2} can induce a similar
proposition in a wider-sense. With Assumption \ref{assu:scarcity_2}
and \ref{assu:growth_2}, we have the Markov property in Lemma \ref{lem:Markov}
and convergences of the probabilities for individuals with infimum
or supremum initial endowments in Lemma \ref{lem:Convergence}. Proposition
\ref{prop:stationary} states that the convergent probabilities are
stationary for these Markov processes.
\begin{lem}
\label{lem:Markov} Given Assumption \ref{assu:scarcity_2} and \ref{assu:growth_2},
$X_{t}(\omega)$ is a Markov process.
\end{lem}

\begin{lem}
\label{lem:Convergence} For any $\{f^{(a)}(\cdot)\}_{a\in\mathcal{A}}$,
select a sequence $\mathbf{a}=\{a_{t}\}$ and denote $X_{t}(\cdot)=f^{(a_{t})}(\cdot)$.
If Assumption \ref{assu:scarcity_2} and \ref{assu:growth_2} are
true, any $\Pr(X_{t}(\underline{\Omega})\leq y)$ is non-increasing
in $t$ and $\Pr(X_{t}(\overline{\Omega})\leq y)$ is non-decreasing
in $t$. Hence, $\Pr(X_{t}(\underline{\Omega})\leq y)$ and $\Pr(X_{t}(\overline{\Omega})\leq y)$
are convergent for each $y$.\end{lem}
\begin{prop}
\label{prop:stationary}If Assumption \ref{assu:scarcity_2} and \ref{assu:growth_2}
hold and $\mathbb{P}\left(X^{*}(\underline{\Omega})\in\mathcal{X}\right)=\lim_{t\rightarrow\infty}\Pr(X_{t}(\underline{\Omega})\in\mathcal{X})$
and $\mathbb{P}\left(X^{*}(\overline{\Omega})\in\mathcal{X}\right)=\lim_{t\rightarrow\infty}\Pr(X_{t}(\overline{\Omega})\in\mathcal{X})$
for some subset $\mathcal{X}\subset\Omega$, then $\mathbb{P}\left(X^{*}(\underline{\Omega})\in\mathcal{X}\right)$
and $\mathbb{P}\left(X^{*}(\overline{\Omega})\in\mathcal{X}\right)$
are stationary. 
\end{prop}
Proposition \ref{prop:stationary} is analogous to the existence of
fixed points in the deterministic case. Stationarity is a way of characterizing
invariant property in stochastic models. The invariance property of
these probability measure implies a Markov equilibrium. A Markov equilibrium
is the stochastic analogue of a steady state in the deterministic
model. The individual growth process is Markovian and its convergent
distribution is stationary when $t\rightarrow\infty$. 

To obtain the uniqueness result of the stationary distributions of
all individuals, we need an assumption that is similar to the equality
condition $X^{*}(1)=X^{*}(0)$ in Assumption \ref{assu:equality_1}.
Because all convergences are for probability measures, the equality
concept should also be adapted to a probabilistic setting. The following
assumption supposes that the poor and the rich may be allocated to
any possible states in the future. If these re-allocations can happen
in probability, it automatically implies the others in the economy
should face similar possibilities too. Since everyone has a stochastic
monotone growth function and everyone can probably has his or her
status changed to any other states, then in the long run, everyone's
growth transition probability may converge to the same stationary
distribution as that of the initially rich or poor person. This implication
is concretized in Theorem \ref{thm:UniqueStationary} given an assumption
about the extreme initial states are not absorbing states.
\begin{assumption}
\label{assu:Equality-prob} (Equality in probability) For any subset
$\mathcal{X}\subset\Omega$, there exists some time $\tau$ such that
\[
\Pr\left(X_{t}(\underline{\Omega})\in\mathcal{X}\right)>0,\;\mbox{and},\;\Pr\left(X_{t}(\overline{\Omega})\in\Omega\backslash\mathcal{X}\right)>0
\]
when $t\geq\tau$.\end{assumption}
\begin{thm}
\label{thm:UniqueStationary} Given Assumption \ref{assu:scarcity_2}
to \ref{assu:Equality-prob}, then there is a unique stationary distribution
such that 
\[
\mathbb{P}\left(X^{*}(\underline{\Omega})\in\mathcal{X}\right)=\mathbb{P}\left(X^{*}(\overline{\Omega})\in\mathcal{X}\right)=\mathbb{P}\left(X^{*}(\omega)\in\mathcal{X}\right)
\]
where $\mathbb{P}\left(X^{*}(\omega)\in\mathcal{X}\right)=\lim_{t\rightarrow\infty}\Pr(X_{t}(\omega)\in\mathcal{X})$
for any $\omega\in\Omega$ in the economy.
\end{thm}
Theorem \ref{thm:UniqueStationary} reveals the ideology of equality
in probability or equality in opportunity. It says that in the equilibrium,
individuals with different types of initial states $\omega$ should
have the same chance to get into a certain state. The opportunities
of entering better or worse states in this equilibrium are the same
for everyone. Theorem \ref{thm:UniqueStationary} is the result of
the unique probabilistic law for an $N$ number of Markov processes. 

Some similar results existed in the economics literature. For example,
the results in \cite{HopenhaynPrescott92} are compatible to Theorem
\ref{thm:UniqueStationary}. The indispensable reason of proposing
this theorem is to show that the stationary Markov processes can be
derived from a new set of assumptions. Assumption \ref{assu:scarcity_2}
to \ref{assu:Equality-prob} are relatively more concrete than the
existing conditions and have straightforward meanings and implications
about growth and equality.%
\footnote{For example, \cite{HopenhaynPrescott92} consider stochastic dominate
condition of the distribution and supermodular condition for the transition
operators while Assumption \ref{assu:growth_2} uses the index set
of monotone production functions and Assumption \ref{assu:Equality-prob}
considers non-degenerated transitions for the growth.%
} These assumptions also induce a less technical proof of the theorem.

\subsection{Remarks }

Assumption \ref{assu:scarcity_2} gives the concept of scarcity on
measurable sets. Assumption \ref{assu:growth_2} allows heterogeneous
individuals to have independent stochastic monotone functions as their
growth functions. Assumption \ref{assu:Equality-prob} implies that
the poor and the rich have opportunities for living in any type of
states in this economy. This gives a broader sense of equality - equality
in probability. By these three assumptions, a stochastic economy is
modeled where all individuals following the same probabilistic law
in the equilibrium. The law implies that the probability of entering
any state in the growth process is the same for everyone. The initial
status plays no impact on the equilibrium probabilistic law. 

Assumption \ref{assu:Equality-prob} illustrates a way of establishing
equality in a broader sense. Equality in probability admits heterogeneous
initial status and uncertain growth, meanwhile it provides an equal
competitive criterion for individuals. Thus in this wider sense, uncertainty
in the growth process does not distort equality from the probabilistic
perspective. Then where does the inequality come from? Note that our
current focus starts from the individual growth. Equality in probability
refers to the equality of individual's probabilistic law. We will
see a completely different image in the following section when the
focus is switched to the aggregates.

\section{Aggregation\label{sec:Aggregation}}

It was pointed out by the Sonnenschein-Mantel-Debreu theorem that
there was no general restriction on the behavior of data aggregates.
An implicit implication is that the general equilibrium results of
the good behavior assumed at the micro-level can not be extended to
the aggregate level of the equilibrium. Based on this point, some
argued that the absence of aggregate empirical restrictions in general
equilibrium theory suggested that theory was incomplete as a way of
understanding economic phenomena. The results in this section show
the opposite of this argument.

In the first part of this section, individual's probabilistic law
from Theorem \ref{thm:UniqueStationary} is shown to be aggregatable
and to be in the same distribution family as the law of linear aggregates.
On the other hand, a twist of the aggregates may simultaneously affect
the laws of individuals. Thus the economy, individuals and aggregates
as a whole, establishes a self-consistent probabilistic law of growth.
This law is called infinite divisible law and it was proposed in \cite{Finetti1929b}
and \cite{Finetti1929a}. 

In the second part of this section, we consider non-linear aggregates.
The stationarity of the linear aggregation is distorted, so is the
infinite divisible law. Recent elaborations of general equilibrium
theory with aggregate features use the methods from the mean field
theory, for example in \cite{PDEmacro2014}. By the mean field theory,
the scope of aggregates in this stochastic economy can be extended
to the non-stationary case. A system of mean field equations is shown
to approximately capture the information of non-stationary aggregates.

\subsection{Infinite Divisibility}

Generally speaking, stationarity property alone does not induce any
particular type of probability distributions. Then it would be a problem
when one considers aggregating the stochastic variables, as stochastic
properties usually vary from the individual level to the aggregate
one. Fortunately the specification of stochastic economy $\left((\Omega,\mathcal{B},(\mathcal{B}_{t}),\mathbb{P}),\mathcal{A},f\right)^{\otimes N}$
and the assumptions restrict all possible stationary distributions
to one specific class in which the linear aggregation, namely the
sum of all individual variables, is a scale invariant of the individual
variable. Without any further assumption, one has the following theorem
for the probabilistic law of the linear aggregates. 
\begin{thm}
\label{thm:inifiniteDivisible} If the economy $\left((\Omega,\mathcal{B},(\mathcal{B}_{t}),\mathbb{P}),\mathcal{A},f\right)^{\otimes N}$
satisfies Assumption \ref{assu:scarcity_2} to \ref{assu:Equality-prob}
and the number of individuals $N\in\mathbb{Z}^{+}$ is arbitrary,
then the stationary distribution $\mathbb{P}(X^{*}(\omega))$ is infinite
divisible. That is, if $X^{*}(\omega^{i})$ is distributed in $\mathbb{P}$
such that $X^{*}(\omega^{i})\sim\mathbb{P}$ for any $i$, then 
\[
\left(X^{*}(\omega^{1})+X^{*}(\omega^{2})\cdots+X^{*}(\omega^{j})\right)\sim\mathbb{P},
\]
for any $0<j\leq N$.
\end{thm}
Because the economy $\left((\Omega,\mathcal{B},(\mathcal{B}_{t}),\mathbb{P}),\mathcal{A},f\right)^{\otimes N}$
allows arbitrary number of individuals, variation of the number $N$
should not have any effect on $\mathbb{P}$. Suppose at the beginning,
there was one individual Adam in the economy who faced the stationary
distribution $\mathbb{P}$ for his uncertain growth path. Then Eve
appeared. As Adam's $\mathbb{P}$ should not be affected by the newcomer
and as Adam and Eve were equal in probabilities, Eve would have the
same $\mathbb{P}$ as her probabilistic law of growth. As time passed
by, more people show up but according to the principle of induction,
their probabilistic law of growth are the same $\mathbb{P}$. Theorem
\ref{thm:inifiniteDivisible} concretizes this idea. 

Infinite divisible law implies that the economy $\left((\Omega,\mathcal{B},(\mathcal{B}_{t}),\mathbb{P}),\mathcal{A},f\right)^{\otimes N}$
can accommodate infinite individuals. The evolution of these probabilistic
laws can be summarized as an array
\[
\begin{array}{l}
\mathbb{P}(\omega^{1})\quad\mbox{Adam's law}\\
\mathbb{P}(\omega^{1}),\mathbb{P}(\omega^{2})\quad\mbox{Adam's and Eve's laws}\\
\vdots\hspace{0.5in}\vdots\hspace{0.5in}\ddots\\
\mathbb{P}(\omega^{1}),\mathbb{P}(\omega^{2}),\ldots,\ldots,\ldots,\mathbb{P}(\omega^{N})\quad\mbox{Current individuals' laws}
\end{array}
\]
namely any probabilistic law $\mathbb{P}$ of the linear aggregation
can be thought as a convolution of $N$ identical laws of $\mathbb{P}$
itself. This triangular array can be extended to infinite entities.

Although the probabilistic law $\mathbb{P}$ is invariant with the
changes of $N$, some other criteria related to $\mathbb{P}$ may
be influenced by such changes. For example, if an individual products
$x$ units, the probability of being the most productive agent for
this individual in an $N$-individuals economy is smaller than in
an $(N-1)$-individuals economy, as
\[
\Pr\left\{ \max_{i\leq N}X(\omega^{i})\leq x\right\} =\prod_{i=1}^{N}\mathbb{P}\left\{ X(\omega^{i})\leq x\right\} 
\]
decreases when $N$ increases. In addition, other possible influences
of aggregation may be on the parameters of the distributions rather
than the distribution family itself. For example, if $\mathbb{P}$
is a Poisson distribution with parameter $1$ in the Adam's stage,
then in the current stage, the parameter of the aggregate $\mathbb{P}$
is still $1$ but individuals' parameter decreases to $1/N$ which
means that a growth event is more difficult to happen. For those who
happen in the growth state, they are in the small fraction of the
whole population and hence this is an unequal economy. Systematic
discussions about these influences are given in Section \ref{sec:Structural-Aggregation}
where the structural parameters of $\mathbb{P}$ are introduced. 

Infinite divisible law includes a large class of probabilistic laws.
Examples of infinite divisible distributions include the normal distribution,
the Poisson distribution, the Cauchy distribution, the $\chi^{2}$
distribution and many others. One interesting distribution is the
stable distribution. 
\begin{cor}
\label{cor:StableID} If the growth function $f(\cdot)$ is linear
such that $X_{t}(\omega)=AX_{t-1}(\omega)+B$, then the stationary
distribution $\mathbb{P}(X^{*}(\omega))$ is a stable distribution.
\end{cor}
Stable distribution admits linear transfers. It means that if stochastic
growth function is linear, any other linear transforms of this stationary
aggregates should also be stationary. All stable distributions follow
infinite divisible law.

\subsection{Mean Field Growth and Volatility }

Non-stationarity is a common feature for growth processes. For example,
one of the most standard growth models, exponential growth model,
induces a geometric Brownian motion $X_{t}(\omega)$ which is a non-stationary
stochastic process.%
\footnote{For unit growth rate and volatility, $X_{t}(\omega)$ is the solution
of the following stochastic differential equation $dX_{t}(\omega)=X_{t}(\omega)dt+X_{t}(\omega)dB_{t}$
where $B_{t}$ is a Brownian motion. %
} Generally speaking, if the underlying system is stationary but the
observable system is non-stationary, this non-stationarity comes from
the nonlinear transformation of the underlying system. For the exponential
growth model, the growth process would be trend stationary if one
applies logarithm transform of $X_{t}$. Since we have stationary
$X_{t}(\omega)$ and $\sum_{i=1}^{j}X_{t}(\omega^{i})$ for any $0<j\leq N$
in the equilibrium, we can consider non-stationary aggregates as non-linear
transforms of some stationary $X_{t}(\omega^{i})$ or $\sum_{i=1}^{j}X_{t}(\omega^{i})$. 

Stationarity and infinite divisibility from Theorem \ref{thm:UniqueStationary}
and \ref{thm:inifiniteDivisible} cannot be preserved if one considers
a non-linear function of some or all $X^{*}(\omega^{i})$. However,
Markov property as a fundamental property about information is preserved. 
\begin{assumption}
\label{assu:AggOneToOne} Non-stationary aggregation follows $Y_{t}(\mathbf{\omega})=g\left(X_{t}(\omega)\right)$
where $X_{t}(\omega)$ represents an arbitrary individual $X_{t}(\omega^{i})$
or a sub-group $\sum_{i=1}^{j}X_{t}(\omega^{i})$ for $0<j\leq N$,
and the function $g(\cdot)$ is one-to-one.\end{assumption}
\begin{cor}
\label{cor:MarkovAgg} If $(X_{t}(\omega^{1}),\dots,X_{t}(\omega^{N}))$
are Markovian, $Y_{t}(\mathbf{\omega})$ follows Assumption \ref{assu:AggOneToOne},
then $Y_{t}(\mathbf{\omega})$ contains the exact information of $(X_{t}(\omega^{1}),\dots,X_{t}(\omega^{N}))$.
Thus $Y_{t}(\mathbf{\omega})$ is Markovian.
\end{cor}
As the state at time $t$ is known, then any information about the
process's behavior before time $t$ is irrelevant. All the relevant
information from the history is stored in $(X_{t}(\omega^{1}),\dots,X_{t}(\omega^{N}))$,
thus $Y_{t}(\mathbf{\omega})$ is a Markov process. Especially, when
$X_{t}(\omega)=X^{*}(\omega)$ is stationary at time $t$, $Y_{t}(\mathbf{\omega})$
is a homogeneous Markov process such that
\[
\Pr\left(\left.Y_{t+h}(\mathbf{\omega})=y_{2}\right|Y_{t}(\mathbf{\omega})=y_{1}\right)=\Pr\left(\left.Y_{t+2h}(\mathbf{\omega})=y_{2}\right|Y_{t+h}(\mathbf{\omega})=y_{1}\right)=\mathbb{Q}_{h}(y_{2}|y_{1})
\]
where $\mathbb{Q}_{h}(y_{2}|y_{1})$ stands for a $h$-step Markov
transition kernel between states $y_{1}$ and $y_{2}$. 

There is a general equation of Markov processes describing the dynamical
probability of varying states. The equation is given as follows:
\begin{equation}
\frac{\partial}{\partial h}\mathbb{Q}_{h}(y_{2}|y_{1})=\int\left[\mathbb{W}(y_{3}|y_{2})\mathbb{Q}_{h}(y_{2}|y_{1})-\mathbb{W}(y_{2}|y_{3})\mathbb{Q}_{h}(y_{3}|y_{1})\right]dy_{2}\label{eq:MasterEq}
\end{equation}
which is called master equation.%
\footnote{Not only is the master equation more convenient for mathematical operations
than some other analytic equations used for Markov processes, such
as Chapman-Kolmogorov, it also holds a more general role in illuminating
Markovian properties.%
} The function $\mathbb{W}(\cdot|\cdot)$, instantaneous transition
rate or simply transition rate, is the probability for a transition
during an extremely short time interval $\Delta t\rightarrow0$. Intuitively
one can think this part of the transition probability is independent
of time. So $\mathbb{W}(\cdot|\cdot)$ holds a different role from
$\mathbb{Q}_{h}(\cdot|\cdot)$ in describing the transition probability.
Appendix \ref{sub:Master-Equation} gives a heuristic derivation of
(\ref{eq:MasterEq}). 

Equation (\ref{eq:MasterEq}) holds for both inhomogeneous and homogenous
Markov processes. Since the equilibrium economy implies the stationary
$X_{t}(\omega)$ and the stationarity implies the homogenous Markov
process $g(X_{t}(\omega))$, it is better to cast the above master
equation in a more intuitive form by using homogenous property. Noting
that all transition probabilities are for a given value $y_{1}$ at
$t$, we can write a simpler expression by suppressing redundant indices
$\mathbb{Q}(y,t)=\lim_{\Delta t\rightarrow0}\mathbb{Q}_{t+\Delta t}(y|y_{1})$.
Equation (\ref{eq:MasterEq}) becomes:
\begin{equation}
\frac{\partial}{\partial t}\mathbb{Q}(y,t)=\int\left[\mathbb{W}(y|y')\mathbb{Q}(y',t)-\mathbb{W}(y'|y)\mathbb{Q}(y,t)\right]dy'.\label{eq:MasterEq-New}
\end{equation}
This is the master equation of transition distribution of $Y_{t}(\mathbf{\omega})$
to the state $Y_{t}(\mathbf{\omega})=y$. The first term is the probability
gain due to transitions from other states $y'$, and the second term
is the probability loss due to transitions into other states $y'$.
A detailed discussion about this equation can be found in \cite{Kampen2007}.

Master equation (\ref{eq:MasterEq-New}) completely determines the
probabilistic law of  $Y_{t}(\omega)$ for all $t$. Since $Y_{t}(\omega)$
is a function of some individual $X_{t}(\omega^{i})$ or groupings
$\sum_{i=1}^{j}X_{t}(\omega^{i})$, $Y_{t}(\omega)$ contains heterogenous
information. From an ordinary macro viewpoint, however, aggregation
procedure ignores fluctuations caused by heterogeneous individuals
or groups that have negligible impacts on the aggregates. This argument
has both empirical and theoretical values. Empirically, specific macroeconomic
information is recorded as a single variable such as GDP or gross
imports and exports. This variable alone indicates the dynamics of
the aggregated growth. Theoretically, although common individual's
growth can be significant in the individual level, it may contribute
very little to the total growth in the whole economy.%
\footnote{Here we do not refer to very important individual growth, such as
important innovations. It is quite likely that some of these individuals
generate non-linearities of $g(\cdot)$ function. But this concern
goes beyond the current context. %
}

Instead of describing all possible fluctuations of $Y_{t}(\omega)$
by $\mathbb{Q}(y,t)$, it is natural to consider a non-stochastic
representative trend of $Y_{t}(\mathbf{\omega})$. Because full information
of $\mathbb{Q}(y,t)$ is difficult to obtain and some $Y_{t}(\omega)\sim\mathbb{Q}(y,t)$
may be not informative. The non-stochastic representative trend is
contained in the moments of $Y_{t}(\omega)$ such as mean and variance.
For a non-stationary Markov process $Y_{t}(\mathbf{\omega})$, the
expected value $\int y\mathbb{Q}(y,t)dy$ is a mean function of $t$
that is denoted as
\[
m_{Y}(t)=\int y\mathbb{Q}(y,t)dy.
\]
This equation only considers the information contained in the first
order moment of $\mathbb{Q}(y,t)$ rather than the whole distribution.
As the master equation determines the entire probability distribution,
it is possible to derive from it the mean-field equation as an approximation
for the case that fluctuations are negligible. The evolution of $m_{Y}(t)$
w.r.t. time $t$ is described by a deterministic differential equation
called the mean-field equation: 
\begin{equation}
\frac{dm_{Y}(t)}{dt}=\int a_{1}(y)\mathbb{Q}(y,t)dy\label{eq:mean_field}
\end{equation}
where $a_{r}(y)$ stands for the $r$-th order moment $a_{r}(y)$
\[
a_{r}(y)=\int(y'-y)^{r}\mathbb{W}(y'|y)dy',\; r=0,1,\dots.
\]
The derivation of (\ref{eq:mean_field}) is given in Appendix \ref{sub:Derivation-mean-field}.
Further deliberation of $a_{1}(y)$ is to see its relation with $m_{Y}(t)$,
as $a_{1}(y)$ in (\ref{eq:mean_field}) is not a linear function
of $m_{Y}(t)$. One can investigate this non-linear relation by expanding
$a_{1}(y)$ around $m_{Y}(t)$ via Taylor series. This gives higher
order information of the evolution in the mean field. Let $\sigma_{Y}^{2}(t)=\int(y-m_{Y}(t))^{2}\mathbb{Q}(y,t)dy$,
the variance of $m_{Y}(t)$. One can think $\sigma_{Y}^{2}(t)$ as
the volatility of the mean field growth.
\begin{thm}
\label{thm:ReducedFormAgg} Given Assumption \ref{assu:scarcity_2}
to \ref{assu:AggOneToOne}, the growth equation of mean field $m_{Y}(t)$
is
\begin{equation}
\mbox{(Reduced Form) }\begin{cases}
\frac{dm_{Y}(t)}{dt} & =a_{1}(m_{Y}(t))+\frac{1}{2}\sigma_{Y}^{2}(t)\\
\frac{d\sigma_{Y}^{2}(t)}{dt} & =a_{2}(m_{Y}(t))+2a_{1}^{(1)}(m_{Y}(t))\sigma_{Y}^{2}(t)
\end{cases}\label{eq:reducedForm}
\end{equation}
where $a_{r}(m_{Y}(t))=\int(y'-m_{Y}(t))^{r}\mathbb{W}(y'|y)dy'$
for $r=1,2$ and $a_{1}^{(1)}(m_{Y}(t))$ denotes 1st-derivative of
$a_{1}(m_{Y}(t))$ w.r.t. $y$.
\end{thm}
Theorem \ref{thm:ReducedFormAgg} is a reduced form description of
the growth dynamics of the whole economy. Equations in (\ref{eq:reducedForm})
only illustrate the first and the second moment evolutions of $\mathbb{Q}$.
But they have no implication about how the equilibrium probabilistic
laws $\mathbb{P}$ influence $\mathbb{Q}$. 

The mean field technique is used to analyze a system with a large
number of components determining the collective deterministic behavior
and seeing how this behavior modifies when the system is perturbed.
In the model (\ref{eq:reducedForm}), the first two moments, $m_{Y}(t)$
and $\sigma_{Y}^{2}(t)$, extract necessary information of the non-stationary
probabilistic law $\mathbb{Q}(y,t)$ by which the collective growth
path is embedded. Although it lacks of structural interpretation,
Theorem \ref{thm:ReducedFormAgg} gives a way of describing the essential
non-stationary dynamics and (\ref{eq:reducedForm}) is easy to implement.
For a reduced form analysis, (\ref{eq:reducedForm}) can capture significant
information of the non-stationary aggregates. An empirical study using
this reduced form model is given in Section \ref{sec:USdata}.

\subsection{Remarks}

The stationary identical distributions, as an equilibrium result of
the stochastic economy, induce restrictions on the behavior of aggregates
both cross-sectionally and inter-temporally. By these restrictions,
a class of probabilistic laws, the infinite divisible law, characterizes
the uncertainty of all individuals and stationary aggregates in this
economy. The appearance of this aggregate specification suggests that
the theory is complete as a way of understanding economic growth.
Common types of aggregate behavior such as power law or Gaussian law
are explained by the infinite divisible law.

Although these laws are identical across heterogenous individuals
and they remain invariance after stationary aggregations, the stochastic
economy endowed with these laws does not alleviate inequality. When
the economy evolutes, all individuals simultaneously alter their probabilistic
laws. Inequality appears and becomes significant during this evolutionary
process. 

When the aggregation involves non-linear patterns, invariant aggregation
of $\mathbb{P}$ is violated. The relation between growth and inequality
is more ambiguous when the growth of aggregates become non-stationary.
A reduced form model is proposed to capture the first two moments'
dynamics of non-stationary aggregates. The causal effect of $\mathbb{P}$
is invisible in the reduced form analysis. To examine the impact of
$\mathbb{P}$, we introduce structural aggregation in the next section.

\section{\label{sec:Structural-Aggregation}Structural Aggregation%
\footnote{The CORE lectures of uncertainty and economic policy given by Jacques
Dr�ze well inspired me to develop several arguments in this section.%
}}

There is one important difference between individual growth and aggregate
growth: interaction effect. Interaction happening between individual
and aggregate variables should be highly asymmetric. One would expect
the aggregate growth to have significant impacts on individual's but
not in the reverse order. Because aggregate growth may benefit the
whole economy while individual growth may only benefit one's adjacency
neighbors. This difference may cause an endogenous issue of inequality
during the aggregation process. Intensive collaborations and concentration
of capitals may lead to an innovative process. This innovative process
may increase the productivity of the whole economy and thus lead to
an equal growth. On the other hand, the innovators may be the first
group of receiving benefits from this innovation. It means that even
an equal aggregate growth may generate derivatives and these derivatives
do not spread equally to the rest of the economy immediately, then
this aggregation may become an endogenous process of creating or enlarging
economic inequality. The interaction term inducing the endogeneities
can be revealed from some higher order information due to the non-linearity. 

In this section, a structural relation is established between the
probabilistic laws of individuals and aggregates. To exploit a tractable
framework of endogenous aggregation, we restrict $\mathbb{P}$ to
a smaller class where structural parameters become visible. By infinite
divisible law of $\mathbb{P}$, these structural parameters are representatives
for both individuals and aggregates. Later, this result is extended
to non-stationary aggregates.

\subsection{Deeper Parameters}

Robert Lucas in his critique suggested looking for deep parameters
that are embedded in the deep layers such as preferences of individuals.
But even if with concerns of these deep parameters, provided that
they are measurable individually, policy suggestions may still be
incomplete as the aggregation process itself has potential endogenous
effects so that the values of parameters may already vary after the
aggregation. We propose an additional requirement for structural parameters.
Not only do the structural parameters should exist in the individual
level, but these parameters also should be measurable and be in an
invariant structure under different scales of aggregates.

Infinite divisible family from Theorem \ref{thm:inifiniteDivisible}
is an ideal category for structural analysis since the distributions
of this family remain invariance under summations and scalings. There
are several ways of characterizing this family. A general characterization
is to consider Feller's semi-group. \cite[Ch 9 and 10]{Fuller1991}
gives an illustration about constructing parametrized generators of
such a semi-group. This approach sheds some light on the way of parameterizing
the structural connection in our context. 
\begin{assumption}
\label{assu:RateAndSize}For any individual $i$, from time $t$ to
$t+\Delta t$, the probability of a significant growth of $X_{t}(\omega^{i})$
is $\alpha_{i}\Delta t$. When $\alpha_{i}\Delta t\geq1$, the event
happens for sure. If the significant growth happen between $X_{t}(\omega^{i})$
and $X_{t+\Delta t}(\omega^{i})$, the total size of changes is $c_{i}\Delta t$. 
\end{assumption}
Parameter $\alpha_{i}$ and $c_{i}$ vary across individuals. When
there is no ambiguity, we drop the index $i$. The use of the term,
significant growth, is to distinguish this growth event from the trivial
individual growth event whose impacts are not strong enough to affect
the economy. The difference of these events can be better understood
when another structural parameter $\beta$ is introduced. Please find
$\beta$ in the following theorem.
\begin{thm}
\label{thm:StructuralAgg}Given Assumption \ref{assu:scarcity_2}
to \ref{assu:RateAndSize}, (I) the probabilistic law of growth becomes
an integro-partial differential equation

\begin{equation}
\frac{\partial}{\partial t}\mathbb{P}(x,t)=-c\frac{\partial\mathbb{P}(x,t)}{\partial x}+\alpha\left[\int\left(\mathbb{W}(x|x')\mathbb{P}(x',t)\right)dx'-\mathbb{P}(x,t)\right]\label{eq:MasterEqGeneral}
\end{equation}
(II) In addition, if $X_{t}(\omega)$ reaches stationary, and the
transition rate is
\[
\mathbb{W}(x|x')=\mathbb{V}(x-x')+\delta(x-x')
\]
where $\mathbb{V}(x-x')=\beta\exp(-\beta(x-x'))$ is an exponential
kernel function, (\ref{eq:MasterEqGeneral}) becomes
\begin{equation}
c(x)\frac{\partial\mathbb{P}(x)}{\partial x}=\alpha\left[\int\left(\mathbb{V}(x-x')\mathbb{P}(x')\right)dx'\right]\label{eq:StationaryGamma}
\end{equation}
where $c(x)=\int_{0}^{x}udu$. The solution of this integro-differential
equation is 
\[
\mathbb{P}(x)=\frac{\beta^{\alpha}x^{\alpha-1}e^{-\beta x}}{\Gamma(\alpha)}\sim\mbox{Gamma}(\alpha,\beta)
\]
where $\Gamma(\alpha)=\int_{0}^{\infty}x^{\alpha-1}e^{-x}dx$ is the
Gamma function.
\end{thm}
The reason of assuming exponential distribution for the kernel of
$\mathbb{W}$ is that exponential distribution is the only distribution
that has continuous memoryless property%
\footnote{Only two kinds of distributions are memoryless: exponential distributions
of non-negative real numbers and the geometric distributions of non-negative
integers. As most indicators of growth and inequality are non-negative
and continuous in $\mathbb{R}^{+}$, exponential distribution becomes
the only option.%
}:
\[
\lim_{\Delta t\rightarrow0}\Pr\left\{ a>t+\Delta t:\: f^{(a)}(x')=x,\: a>t\right\} =\lim_{\Delta t\rightarrow0}\Pr\left\{ a>\Delta t:\: f^{(a)}(x')=x\right\} 
\]
where $a\in\mathcal{A}$ in the economy $\left((\Omega,\mathcal{B},(\mathcal{B}_{t}),\mathbb{P}),\mathcal{A},f\right)^{\otimes N}$
is the random time of a growth event. A description about memoryless
and exponential distribution is given in Appendix \ref{sub:Discussion-of-Memoryless}. 

Equation (\ref{eq:MasterEqGeneral}) is one of the most important
results in this paper. It gives a specific representation of the probabilistic
laws of $X_{t}(\omega)$. Meanwhile, it keeps a general enough formulation
to cover a large amount of interesting cases. The Gamma distribution
as a stationary solution of (\ref{eq:MasterEqGeneral}) includes some
other standard distributions such as the $\chi^{2}$-square distribution,
the exponential distribution, etc, and can approximate a large class
of distributions such as the log-normal distribution and the power
law distribution. Other implications of (\ref{eq:MasterEqGeneral})
are illustrated by two corollaries below, Corollary \ref{cor:Stein}
and \ref{cor:PowerLaw}. In addition, (\ref{eq:MasterEqGeneral})
reveals the structural meaning of its parameters.
\begin{lyxlist}{00.00.0000}
\item [{$\alpha$}] is called shape parameter in Gamma distributions. In
the current model, it comes from Assumption \ref{assu:RateAndSize}.
It represents the probability of a significant growth. For the Gamma
distribution, it is known that 
\[
X_{t}(\omega^{i})+X_{t}(\omega^{j})\sim\mbox{Gamma}(\alpha_{i}+\alpha_{j},\beta)
\]
if $X_{t}(\omega^{i})\sim\mbox{Gamma}(\alpha_{i},\beta)$ and $X_{t}(\omega^{j})\sim\mbox{Gamma}(\alpha_{j},\beta)$.
Thus $\alpha$ is an aggregatable parameter. When $\beta$ is invariance
for all individuals and the stationary aggregates, $\alpha$ is the
only parameter that captures heterogeneities. For an individual, as
this parameter can influence the aggregate growth of the economy,
this individual growth is significant enough. We can attribute technological
development or innovations to this type of growth. Thus heterogeneity
accounts for the global impacts in the growth. It is natural to think
a group of unskillful new members joining the economy with very small
$\alpha$, since their participations have negligible impacts on the
aggregate production. 
\item [{$\beta$}] is called rate parameter in Gamma distributions. In
the current model, it comes from the memoryless exponential kernel
function. As the exponential kernel is assumed for any $X_{t}(\omega)$,
the value of $\beta$ remains the same for all individuals. Thus $\beta$
is a parameter characterizing those individual growth events that
has no impact on the aggregates. It is known 
\[
C\times X_{t}(\omega)\sim\mbox{Gamma}(\alpha,\beta/C)
\]
where $C$ is a factor of the total population and $X_{t}(\omega)=\sum_{i=1}^{N}X_{t}(\omega^{i})$
having a Gamma distribution $\mbox{Gamma}(\alpha,\beta)$. When a
group of new members joining the economy, say increasing one percentage
of $X_{t}(\omega)$, there may be no change of the aggregate $\alpha$,
however their participations decrease $\beta$ to $\beta/1.01$. The
value of $\beta^{-1}$ measures a scale. For exponential distribution,
the scale of an individual growth event has a mean of $1/\beta$.
Thus with the evolution of an economy, namely more producers, the
growth of individual, if it happens, has a bigger scale. 
\item [{$c(x)$}] is an accumulation from $0$ to $X_{t}(\omega^{i})=x_{i}$
in (\ref{eq:StationaryGamma}). It comes from Assumption \ref{assu:RateAndSize}.
It measures the size of a significant growth event. Although the stationary
solution does not display its contribution, from (\ref{eq:StationaryGamma})
we can see that $c(x)\frac{\partial\mathbb{P}(x)}{\partial x}$ on
the left hand side balances the probability changes in the right hand
side. It means that for individual $i$, $\alpha_{i}$ on right hand
side of (\ref{eq:StationaryGamma}) depends on the accumulation up
to the current state $X_{t}(\omega^{i})=x_{i}$. If the state of individual
$i$ is larger than individual $j$, $X_{t}(\omega^{i})>X_{t}(\omega^{j})$
and $c(x_{i})>c(x_{j})$, then it is likely that $\alpha_{i}>\alpha_{j}$.
This inequality does not happen for sure, as it also depends on $\mathbb{P}$
and its derivative. But we can see similar phenomena in economics,
intensive capital investments often lead to high productivities. 
\end{lyxlist}
Stationary aggregates $\sum_{i=1}^{N}X_{t}(\omega^{i})$ follows $\mbox{Gamma}(\sum_{i=1}^{N}\alpha_{i},\beta)$.
Its mean and variance are 
\[
\mathbb{E}\left[\sum_{i=1}^{N}X_{t}(\omega^{i})\right]=\frac{\sum_{i=1}^{N}\alpha_{i}}{\beta},\quad\mbox{Var}\left[\sum_{i=1}^{N}X_{t}(\omega^{i})\right]=\frac{\sum_{i=1}^{N}\alpha_{i}}{\beta^{2}},
\]
where $X_{t}(\omega^{i})\sim\mbox{Gamma}(\alpha_{i},\beta)$. It means
that even though the aggregates growth events can be accumulated w.r.t.
$\alpha_{i}$, this aggregation is not influenced by the homogenous
personal growth events characterized by $\beta$. On the other hand,
as population growth (increasing $N$) causes a decrease of $\beta$,
this population growth increases mean and variance of stationary aggregates
simultanesouly. More discussion about the aggregates of Gamma distribution
is given in Section \ref{sub:Remarks-StruAgg}. Early use of Gamma
density as a descriptive statistics of incomes can be found in \cite{SalemMount1974}. 

The above discussions are for stationary aggregates. For non-stationary
aggregates, we need to modify the mean field method from the previous
section. The purpose is to connect the mean field information of $\mathbb{Q}$
with the specified $\mathbb{P}$-law given by (\ref{eq:MasterEqGeneral}).
\begin{thm}
\label{thm:StructuralForm}Given the results of (\ref{eq:mean_field})
and (\ref{eq:MasterEqGeneral}), if $\mathbb{P}(x,t)$ and $\mathbb{Q}(y,t)$
are mutual continuous, such that $\mathbb{P}(x,t)=0$ implies $\mathbb{Q}(y,t)=0$
and vice versa, then there is a structural form for the non-stationary
aggregation 
\begin{equation}
\mbox{(Structural Form) }\begin{cases}
m_{Y}(t) & =\mathbb{E}\left[\left.g\left(X_{t}(\omega)\right)\mathbb{L}(t)\right|X_{t}(\omega)\right]\\
\frac{dm_{Y}(t)}{dt} & =\theta_{t}\, m_{Y}(t)+e_{t}+o(1)\\
X_{t}(\omega) & \sim\mathbb{P}(x,t)\mbox{ as a solution of eq.\eqref{eq:MasterEqGeneral}}
\end{cases}\label{eq:StructuralEq}
\end{equation}
where $\mathbb{L}(t)=[\mathbb{Q}(y,t)dy/\mathbb{P}(x,t)dx]$ is the
likelihood ratio, $\mathbb{E}[\cdot]$ is taken w.r.t. $\mathbb{P}$,
$\theta_{t}$ and $e_{t}$ are time varying constants.
\end{thm}
The first equation in (\ref{eq:StructuralEq}) establishes the relation
between the mean field of non-stationary aggregate $m_{Y}(t)$ and
stationary aggregate $X_{t}(\omega)=\sum X_{t}(\omega^{i})$. It links
the mean field of non-stationary law $\mathbb{Q}$ with the expectation
evaluated by the stationary law $\mathbb{P}$. The connection of first
order information of $\mathbb{Q}$ and $\mathbb{P}$ is contained
in this equation. Since $\mathbb{P}$ is parameterized and $m_{Y}(t)$
has a reduced form representation, it is easy to empirically estimate
this equation.

Affine structure with time varying constants often comes from point-wise
linearization of nonlinear functions. The second equation of (\ref{eq:StructuralEq})
is in this structure. Time varying parameters make the linear equation
possible to consistently capture higher order information along the
dynamics. The discrete time stochastic version of this equation is
analogous to a recursive structure called Kesten process. It has been
used in finance, see \cite{Gabaix2009}. Since we have additional
information of $\mathbb{P}$, it would be easier to structurally estimate
this affine equation rather than consider Kesten processes. Another
advantage of using affine structure in (\ref{eq:StructuralEq}) is
the tractability. Having tractable solutions for recursive estimation
is useful because it has a closed form expression for each step updating
instead of casting a black-box algorithm. Implementation of Theorem
(\ref{eq:StructuralEq}) is given in Section \ref{sec:USdata}.

By dividing the term $m_{Y}(t)$ in the second equation in (\ref{eq:StructuralEq}),
this equation can be simplified as 
\begin{equation}
\frac{d\ln m_{Y}(t)}{dt}=\theta_{t}+\frac{e_{t}}{m_{Y}(t)}\label{eq:LogMean}
\end{equation}
as $d\ln m_{Y}(t)=(m_{Y}(t))^{-1}dm_{Y}(t)$. The differentiation
of logarithm refers to the growth rate of the mean field $m_{Y}(t)$.
One can replace $m_{Y}(t)$ in (\ref{eq:LogMean}) with the first
equation of (\ref{eq:StructuralEq}). Since computation of expected
$g(\cdot)$ is feasible by using the parametrized $\mathbb{P}$-law,
one can use stationary aggregates $X_{t}(\omega)$ and the structural
parameters $(\alpha,\beta)$ from $\mathbb{P}$-law to describe the
mean field dynamical of non-stationary $\mathbb{Q}$-law. 

Changes of $\alpha$ and $\beta$ are the endogenous forces of the
growth. In the equilibrium economy, $\alpha$ and $\beta$ are fixed
on the individual level, but on the global level the aggregate impact
can increase productivities for some heterogenous individuals (cause
the change of $\alpha$) and evolution can increase the size of the
economy (cause the change of $\beta$). Because the mean and the variance
of the Gamma distribution are $\alpha/\beta$ and $\alpha/\beta^{2}$,
$\alpha$ and $\beta$ have different propagate effects over the mean
and the volatility of the growth. This endogeneity may relate to one
open question in growth theory about the correlation between growth
and growth volatility. Business cycle volatility and growth has been
extensively studied. Many models showed that growth volatility negatively
affects growth and that economies with higher volatility experience
lower growth. With the structural relation between $\mathbb{P}$ and
$\mathbb{Q}$, we can decompose the volatility effects and examine
their deep role in the growth. In Section \ref{sec:USdata}, we give
a structural estimate of the growth volatility. The result shows that
the causality between volatility and growth is mainly due to endogenous
variations of $\alpha$ and $\beta$.
\begin{cor}
\label{cor:Stein}If $X_{t}(\omega)$ reaches stationary, $\alpha=x$,
$c=-1$ and $\mathbb{W}(x|x')\equiv0$ degenerates, (\ref{eq:MasterEqGeneral})
becomes 
\[
\frac{d\mathbb{P}(x)}{dx}-x\mathbb{P}(x)=(\mathbb{S}\circ\mathbb{P})(x)=0
\]
where $\mathbb{S}$ is the Stein operator such that $(\mathbb{S}\circ f)(x)=(\partial f(x)/\partial x)-xf(x)$
for a function $f$.
\end{cor}
This operator is one of the major devices to prove central limit theorem.
Thus most, if not all, i.i.d. sums can be approximated by this operator.
Please refer to \cite{Stein1986} for its details.

\subsection{\label{sub:Remarks-StruAgg}Remarks}

This section attempts to give a structural description of the aggregation
process. The concern of the structure in the aggregation is that the
probabilistic laws of individual growth events may aggregate simultaneously
thus the parameters of these laws may vary after the aggregation.
The specification of Assumption (\ref{assu:RateAndSize}) leads to
a refined equilibrium solution, Gamma distribution. This distribution
follows the infinite divisible law. It has two different parameters,
$\alpha$ and $\beta$, that separately characterize events that have
heterogenous aggregate impacts and have homogenous individual impacts.
This specification makes the structural analysis of aggregation possible.
Parameter $\alpha$ is aggregatable so that $\alpha$ of the aggregates
contains heterogeneous individual information measured by $\alpha_{i}$.
Parameter $\beta$ evaluates the homogenous growth opportunities of
the whole economy and it also counts the evolutionary influence of
the economy structure. After parametrizing these simultaneous effects
along with the aggregate growth, one can return to epistemic state
of inequality encoded by these parameters. Because these parameters
consistently characterize the equilibrium distributions of this economy. 

\begin{figure}
\subfloat[Mean invariance]{\includegraphics[width=0.45\textwidth]{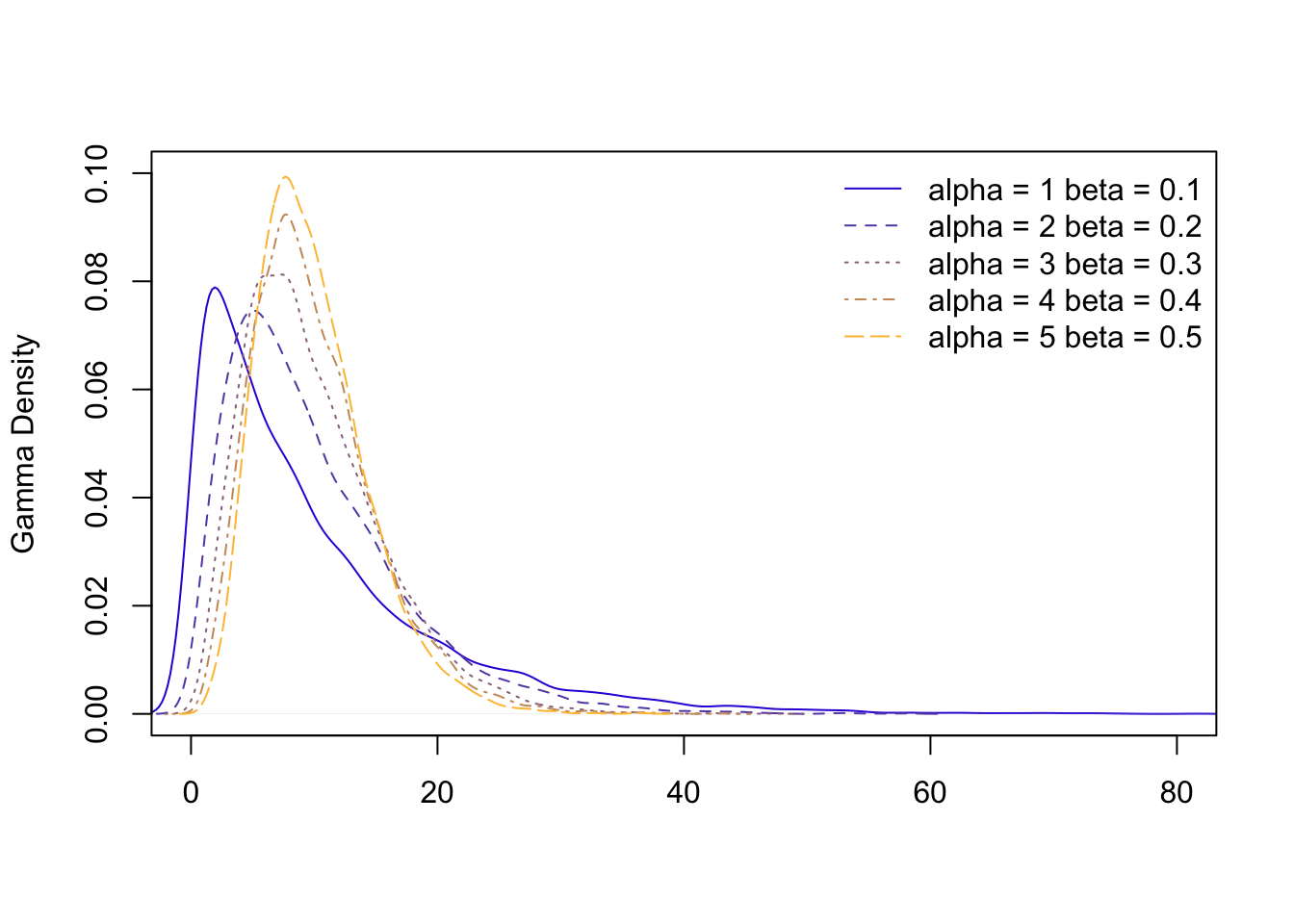}

}\hfill{}\subfloat[Variations of $\beta$]{\includegraphics[width=0.45\textwidth]{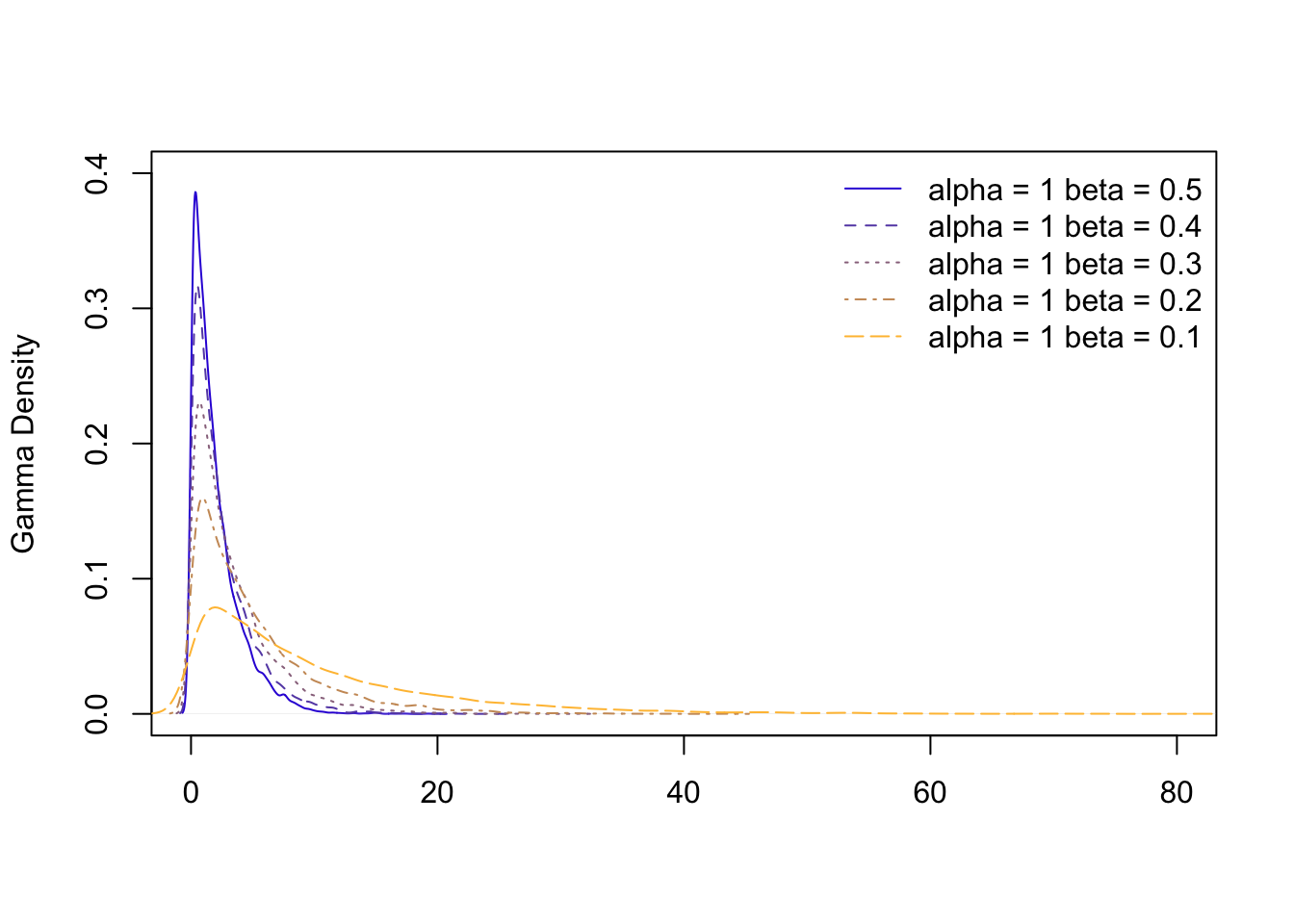}

}

\subfloat[Variations of $\alpha$]{\includegraphics[width=0.45\textwidth]{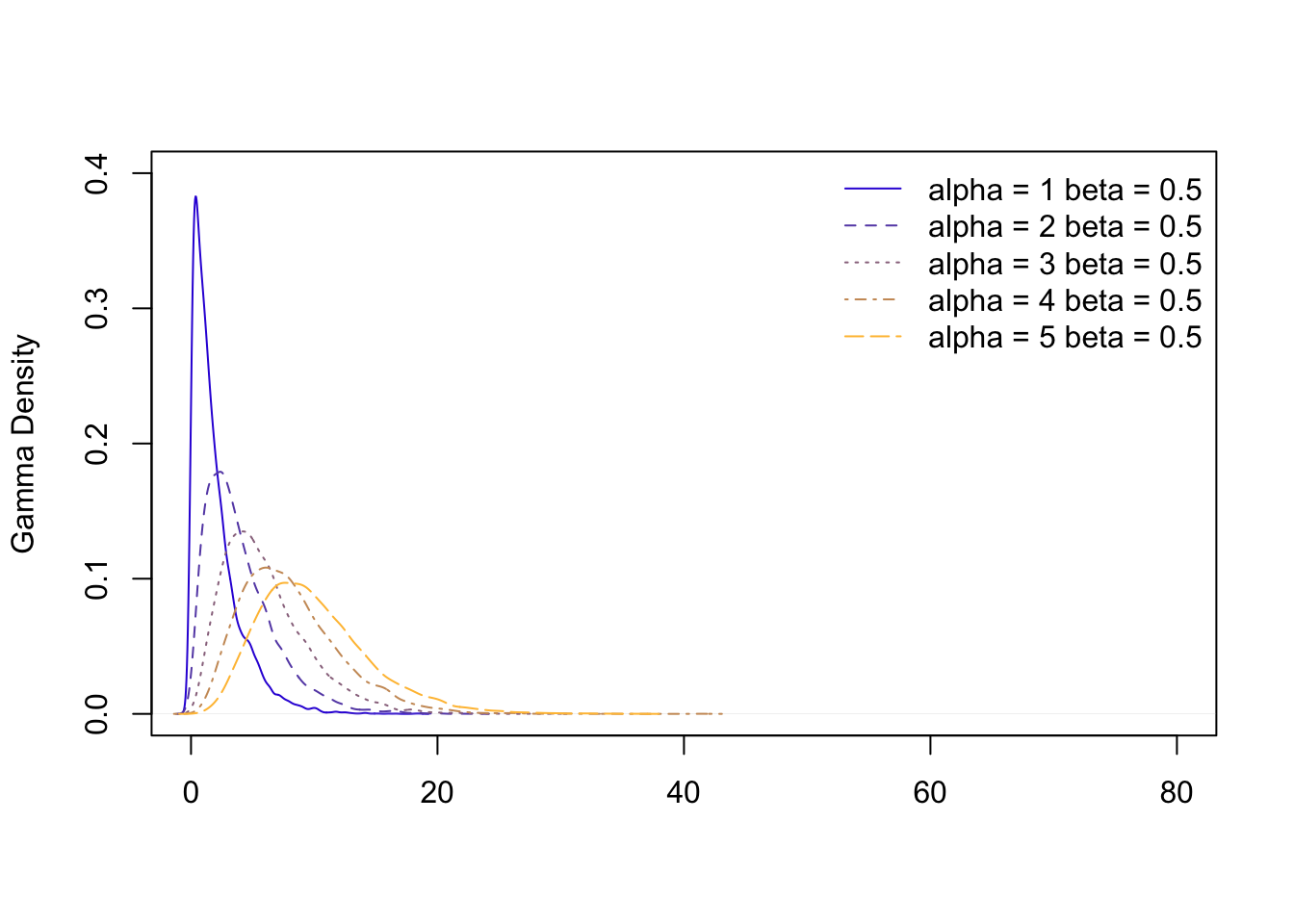}

}\hfill{}\subfloat[Variations of $(\alpha,\beta)$]{\includegraphics[width=0.45\textwidth]{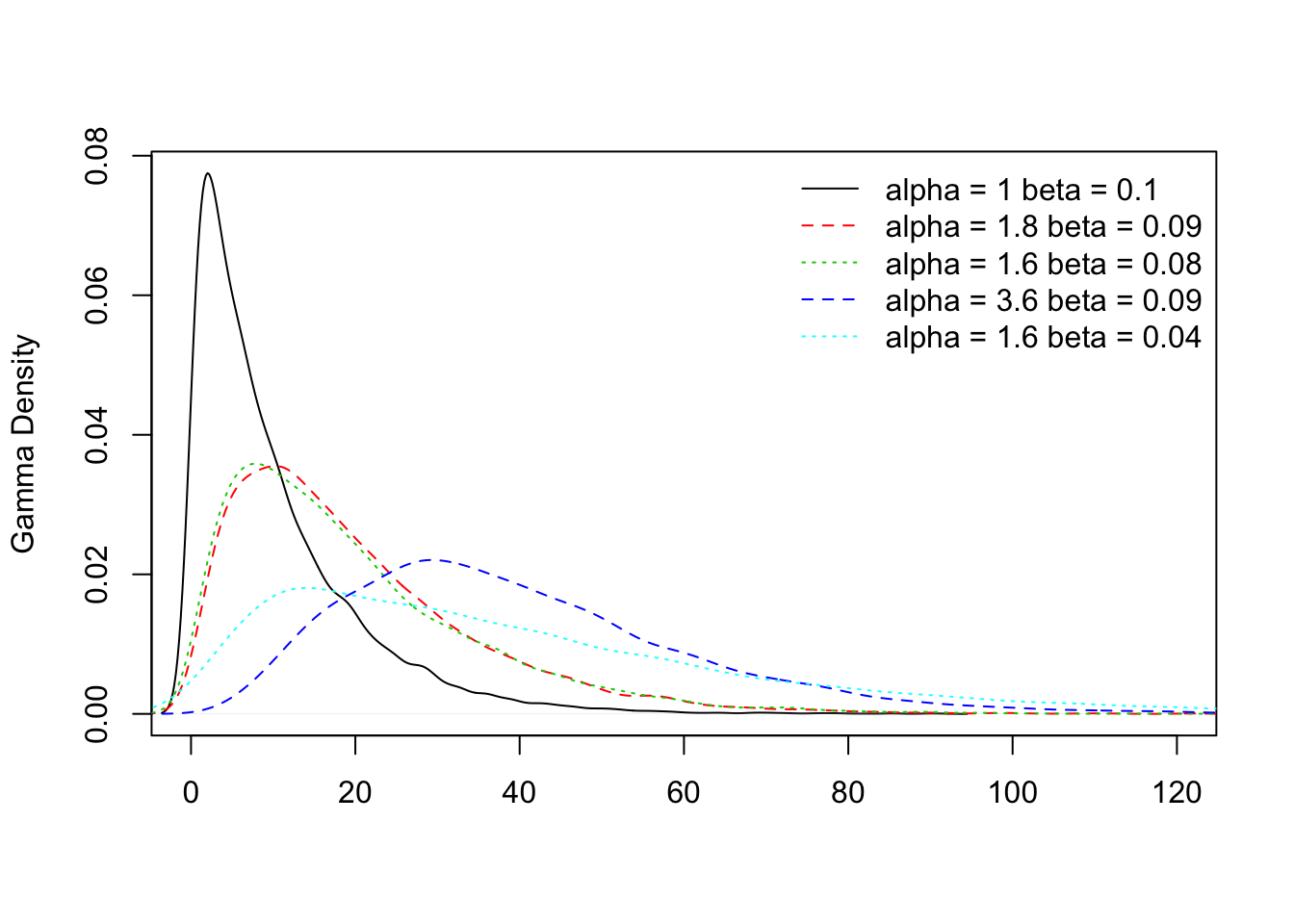}

}

\begin{centering}
\subfloat[Variations of $(\alpha,\beta)$ in the tail]{\includegraphics[width=0.6\textwidth]{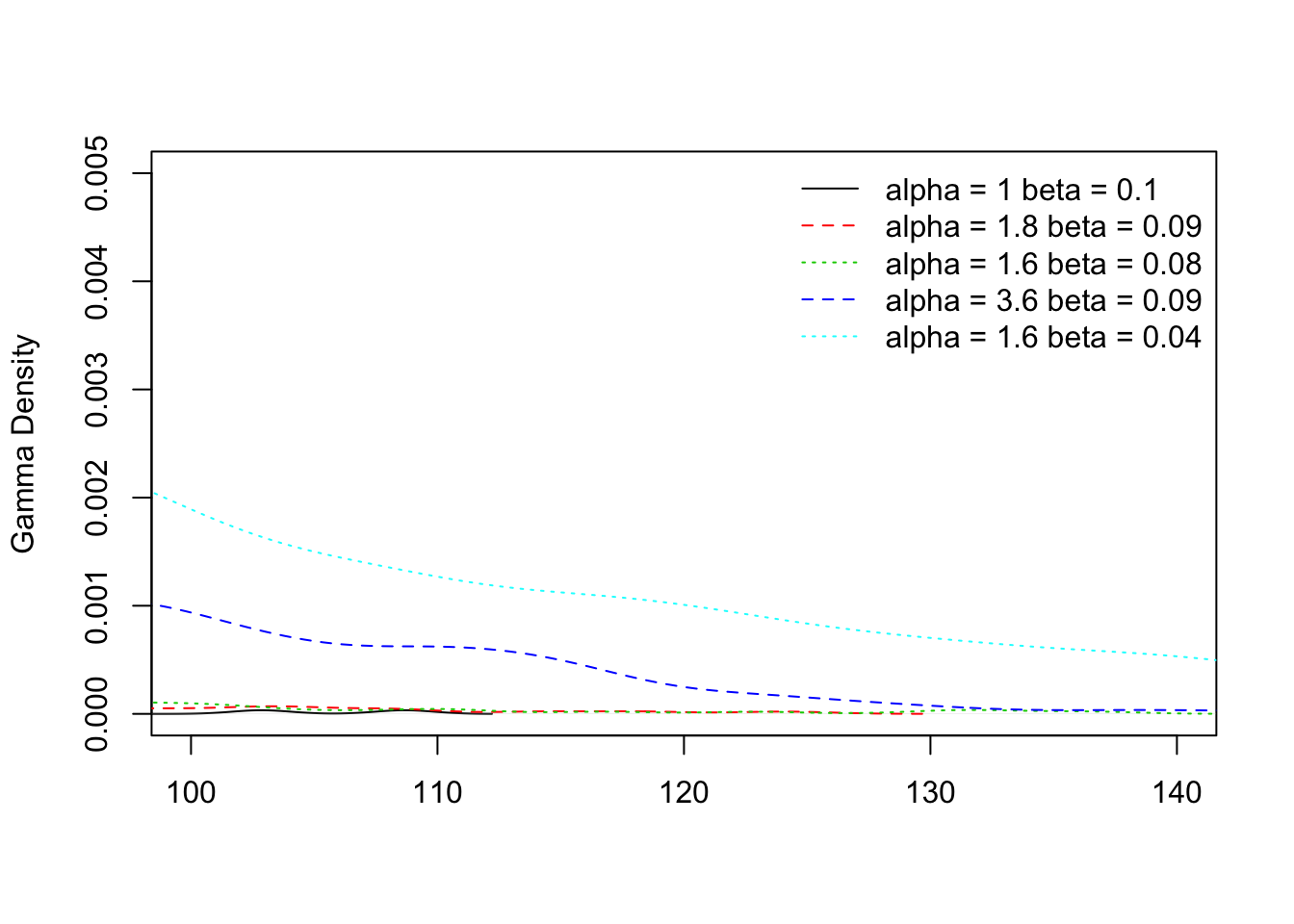}

}
\par\end{centering}

\protect\caption{\label{fig:Gamma-Distribution}Gamma Distribution}
\end{figure}

Figures \ref{fig:Gamma-Distribution} show some possibilities for
altering $\alpha=\sum\alpha_{i}$ and $\beta$ in order to give an
integrated insight of the roles of $\alpha$ and $\beta$ in growth
and inequality. From the figures, one can easily distinguish the contributions
made by $\alpha$ and $\beta$ to inequality in a dynamical setting.
Figure \ref{fig:Gamma-Distribution}(a) by holding the same mean value,
increasing both $\alpha$ and $\beta$ shifts the center to the right.
Although percentages in the high value states reduce, more is gained
as the percentages in the low value states reduce more significantly.
Figure \ref{fig:Gamma-Distribution}(b) considers changes of $\beta$.
By reducing $\beta$ from $0.5$ to $0.1$, we can see more and more
probabilities tend to be in the high value states. However, the changes
for low value states are not so significant. Figure \ref{fig:Gamma-Distribution}(c)
considers the same initial position as (b) then we increase the values
of $\alpha$. The final position of the dynamics share the same mean
as the final one in (b). However, the whole shape of the distribution
moves more significantly to the right. The movements around the low
value regions are much more significant than those happened for $\beta$,
even though the changes of means are the same. Figures \ref{fig:Gamma-Distribution}(d)
and (e) consider increasing $\alpha$ and decreasing $\beta$ simultaneously
but in different scales. By holding the same means, one trend is to
increase more $\alpha$ the other is to decrease more $\beta$. For
subtle changes, $(\alpha=1.8,\beta=.09)$ and $(\alpha=1.6,\beta=.04)$,
the effects are quite similar. But if we amplify the trends, we can
see that an increase of $\alpha$ eliminates larger percentages of
lower value states than a decrease of $\beta$, moreover, the decrease
of $\beta$ attributes more probabilities to very high value states. 

If we consider the value of states below $20$ as poor status, Figures
\ref{fig:Gamma-Distribution}(b) and (c) demonstrate that increasing
$\alpha$ can more efficiently reduce the poverty percentage. It means
that encouraging individuals to make significant contributions to
the aggregate growth can be helpful to reduce the inequality issue
of the economy. To summarize these figures, the role of $\alpha$
emphasizes individual strengths that have impacts on the aggregates,
the role of $\beta$ considers growth events that happen with an equivalent
size $1/\beta$ across all individuals in the economy. To increase
the mean of stationary aggregates, one can think of increasing $\alpha$,
decreasing $\beta$ or doing both. From the dynamical patterns in
Figure \ref{fig:Gamma-Distribution}, one can see that within the
same magnitude of mean value changes, $\alpha$ plays a more important
role in reshaping the distribution structure.

\subsection{Equality Paradox}

Previous discussions and examples of Gamma distributions suggest that
increasing $\alpha$ can more efficiently reduce the poverty percentage
than decreasing $\beta$. It suggests that encouraging individuals
to make significant contributions to the aggregate growth can be more
helpful to reduce the inequality issue even though equivalently increasing
possible growth size for all individuals seems more fair at first
glance. As change of $\beta$ refers to an equality strategy, it is
counter-intuitive to think that an economic growth by reducing $\beta$
may arise inequality. This is an equality paradox: \emph{an attempt
of achieving equality generates inequality}.

Increase of $1/\beta$ in fact only matters when growth event happens.
For example, increasing returns to all inventions only matters to
those who have made an invention, however everyone has a probability
of inventing something thus this policy is supposed to be beneficial
to everyone. What hidden in paradox is that it emphasizes the equality
of opportunity and returns (everyone has the same law and everyone
has the same return if the event happens) but it does not mention
the distortion of opportunities for individuals. Since the size of
individual growth events increase for everyone, those who are in the
motion of growing will have bigger improvements of their states. Even
if this movement has so little effect to the others that $\alpha$
remains the same, the economy could become more unequal. Because those
who seize the opportunities will get bigger returns but those who
have no opportunity in this round will face the same situation in
the future. If this trend continues, namely winners get more and losers
keep the same expectation, inequality will be enlarged. 

Things could get worse, if the trend of increasing $1/\beta$ happens
unconsciously. In this model, the evolution of economy automatically
admits new labor forces with the same stochastic production patterns
as the existing ones. By scaling property of $\mathbb{P}$-laws, it
implies an increasing trend of $1/\beta$ along this evolutionary
path. Additionally, a policy that is intended to meet social moral
criteria often claims more for equality meanwhile economic instruments
of the policy are often intended to use rewards as incentive plans.
So $1/\beta$ increases. In the process, people may not even realize
their actions of generating inequality. It is even harder to perceive
the trend if the whole economy is in an expansion. When the effect
of decreasing $\alpha$ is compensated by increasing $1/\beta$, the
mean $\alpha/\beta$ can still be the same as before or even higher.
Thus one should be aware that some parts of the growth caused by $\beta$
may be a sacrifice of the value $\alpha$. This tendency is shown
clearly in Figures \ref{fig:Gamma-Distribution}(d) and (e). Higher
$1/\beta$ increases the percentage of upper and upper middle class
but the poor almost remain the same as before. New upper and upper
middle class consist mainly of those from the previous middle class.

\subsection{Conjectural Strips }
\begin{cor}
\label{cor:PowerLaw}If $X_{t}(\omega)$ reaches stationary, $c=x$
and $\mathbb{W}(x|x')\equiv0$ degenerates, (\ref{eq:MasterEqGeneral})
becomes
\begin{equation}
x\frac{d\mathbb{P}(x)}{dx}=\alpha\mathbb{P}(x)\label{eq:PowerLaw}
\end{equation}
whose solution is $\mathbb{P}(x)=\mathbb{P}(1)x^{-\gamma}$, a power
law distribution with parameter $\gamma=\alpha^{-1}$. 
\end{cor}
Corollary \ref{cor:PowerLaw} gives a continuous version of the model
used in \cite{Simon1955} where the power law comes from a difference
equation. For more discussion about power law and its applications
in economics, please refer to the survey \cite{Gabaix2009}. From
Corollary \ref{cor:PowerLaw}, one can easily obtain the Riemann Zeta
distribution.%
\footnote{If the infinite divisible distribution comes from a counting data
set, the Zeta distribution (or empirically Zipf's law) rather than
the Gamma distribution would be a better candidate to describe the
law how the data points distribute.%
} A short description is given in Appendix \ref{sub:Zeta-Distribution}.

The Zeta distribution from (\ref{eq:PowerLaw}) and the Gamma distribution
from (\ref{eq:StationaryGamma}) can be connected by the Riemann functional
equation
\begin{equation}
\zeta(s)=2^{s}\pi^{s-1}\sin\left(\frac{s\pi}{2}\right)\Gamma(1-s)\zeta(1-s).\label{eq:RiemannFunctionalEq}
\end{equation}
It connects the Riemann Zeta function $\zeta$ from the denominator
of the Zeta distribution and the Gamma function $\Gamma$ from the
denominator of the Gamma distribution. Equation (\ref{eq:RiemannFunctionalEq})
can have multiple zeros at those $s=-2,-4,\dots$ such that $\sin(s\pi/2)=0$.
These are called trivial zero of Zeta function. As the Zeta distribution
is not defined for negative $s$, we can only consider the corresponding
Gamma distributions at $\alpha=3,5,7,\dots$. Suppose that there are
two connected systems, one is an equilibrium economy with stationary
Gamma distribution from (\ref{eq:StationaryGamma}), the other one
with an equilibrium characterized by a Zeta function $\zeta$. Assume
that these two systems are connected by the same parameter $\alpha$.
When $\alpha=3,5,\dots$, although the economy has stationary solutions,
the other system faces zeros of $\zeta$. If this $\zeta$ is the
denominator of an equilibrium solution, such as $\zeta$ in the Zeta
distribution of (\ref{eq:PowerLaw}), then the other system is arriving
at a singular point. In this case, even $\alpha$ does not cause any
significant change for the economy, economy may confront some shocks
from the other system. 

Another hypothesis to consider is the case of $s=1/2$. At $s=1/2$,
equation (\ref{eq:RiemannFunctionalEq}) gives a unique result of
$\Gamma(1/2)$ which is $\sqrt{\pi}$.%
\footnote{When $s$ is a complex number, the Riemann hypothesis asserts that
any $s$ satisfying $\zeta(s)=0$ and $s\neq2k$ (non-trivial zero)
locates at $\mbox{Re}(s)=1/2$, the so-called critical line. %
} Moreover, Stirling's formula gives a useful recurrence form for $\sqrt{\pi}$
that relates to $\Gamma(\cdot)$ 
\[
\Gamma\left(k+\frac{1}{2}\right)=\Gamma\left(\frac{1}{2}\right)\cdot\frac{1}{2}\cdot\frac{3}{2}\cdot\dots\cdot\frac{2k-1}{2}=\sqrt{\pi}\cdot\frac{(2k)!}{2^{2k}k!}.
\]
Since $\pi$ often refers to a stable and harmonic state, a natural
question is whether these values have meanings for economic dynamics?
Does $\pi$ affect the economic activities? Does an economy tend to
stay in a growth path with $\alpha=k+1/2$ such as $\alpha=3/2,\dots$
rather than other values? Here is a conjecture: $\alpha$ that indicates
the vigor of an economy has some critical values, some values may
cause singular effects in the related systems, some values may establish
stable states that are difficult to pass over. These values could
be the critical strips for the growth of human beings.

\section{\label{sec:USdata}Growth and Inequality in U.S. 1994-2015}

Two data sets are used for the illustration. For growth, the data
information comes from U.S. real GDP per capita issued by the U.S.
Bureau of Economic Analysis. It is a quarterly data set available
from $1947$ to the present. For personal income, the data comes from
Current Population Survey (CPS) issued by the U.S. Bureau of Labor
Statistics and the Census Bureau. We choose the category that records
personal total money income from persons of $15$ years old and over.
Both data sets are public and available on-line.%
\footnote{\inputencoding{latin1}Information about data and implementation is
available in the online-appendix: \inputencoding{latin9}\foreignlanguage{english}{\url{https://rpubs.com/larcenciel/UncertaintyA}}\selectlanguage{english}%
} Time period from $1994$ to $2015$ is selected to make the time
range of two data sets comparable. 

Estimates are divided into two parts. The first part considers the
reduced form mean field system that is developed in (\ref{eq:reducedForm}).
Aggregate GDP is assumed to stochastically follow $\mathbb{Q}$-law
and its realization is $m_{Y}(t)$. The main interest is to capture
non-stationary dynamics of $\mathbb{Q}$. We estimate the parameters
that specify the dynamics of the first two moments of $\mathbb{Q}$.
The second part considers the structural form system that is developed
in (\ref{eq:StructuralEq}). Income distribution is treated as an
equilibrium aggregate solution following $\mathbb{P}$-law. A structural
relationship between GDP and income distribution is exploited. We
first estimate structural parameters that characterize $\mathbb{P}$-law
and then use them to make a further estimate for the time varying
parameters that characterize $\mathbb{Q}$-law. The structural relationship
between income and GDP data eliminates an endogenous effect that is
shown in the regression pre-analysis. The scheme of both estimates
is given in Appendix \ref{sub:Scheme-of-Estimates}.

First, let $m_{Y}(t)$ in (\ref{eq:reducedForm}) represent the realized
GDP per capita in year $t$. The continuous time representation of
(\ref{eq:reducedForm}) needs to be discretized for time series data
and its regression analogy is given as follows 
\begin{align}
m_{Y}(t+1)-m_{Y}(t) & =\mbox{Constant}+\mbox{Coef}_{1}m_{Y}(t)+\mbox{Coef}_{2}\sigma_{Y}^{2}(t)+\varepsilon_{t}\label{eq:rf-1}\\
\sigma_{Y}^{2}(t+1)-\sigma_{Y}^{2}(t) & =\mbox{Constant}+\mbox{Coef}_{3}m_{Y}^{2}(t)+\mbox{Coef}_{4}\sigma_{Y}^{2}(t)+\upsilon_{t}\label{eq:rf-2}
\end{align}
where $\varepsilon_{t}$, $\upsilon_{t}$ are assumed to be white
noise. The variance $\sigma_{Y}^{2}$ is constructed by the residual
from a pre-estimate of (\ref{eq:rf-1}) without $\sigma_{Y}^{2}$.
The model shares some similarities with the conditional heteroskedasticity
models in statistics and econometrics such as the class of (generalized)
autoregressive conditional heteroskedasticity (GARCH) models. As $Y_{t}(\omega)$
is generated by nonlinear function $g(\cdot)$, one would expect severe
heteroskedasticities caused by $g(\cdot)$. This concern is well understood
in the statistics literature. But there are two crucial differences.
First, the variance term $\sigma_{Y}^{2}(t)$ does not relate to the
noise $\varepsilon_{t}$. It is the higher order information term
from the expansion of mean field equation while in heteroskedasticity
models $\sigma_{Y}^{2}(t)$ represents the variance of the noise.
Second, as $\sigma_{Y}^{2}(t)$ is an expansion term, it also depends
on the mean field function $m_{Y}(t)$. This never happens in heteroskedasticity
models. These two differences are reflected by the coefficient of
$\sigma_{Y}^{2}(t)$ in (\ref{eq:rf-1}) and the coefficient of $m_{Y}^{2}(t)$
in (\ref{eq:rf-2}). 

\begin{figure}
\subfloat[]{\includegraphics[width=0.45\textwidth]{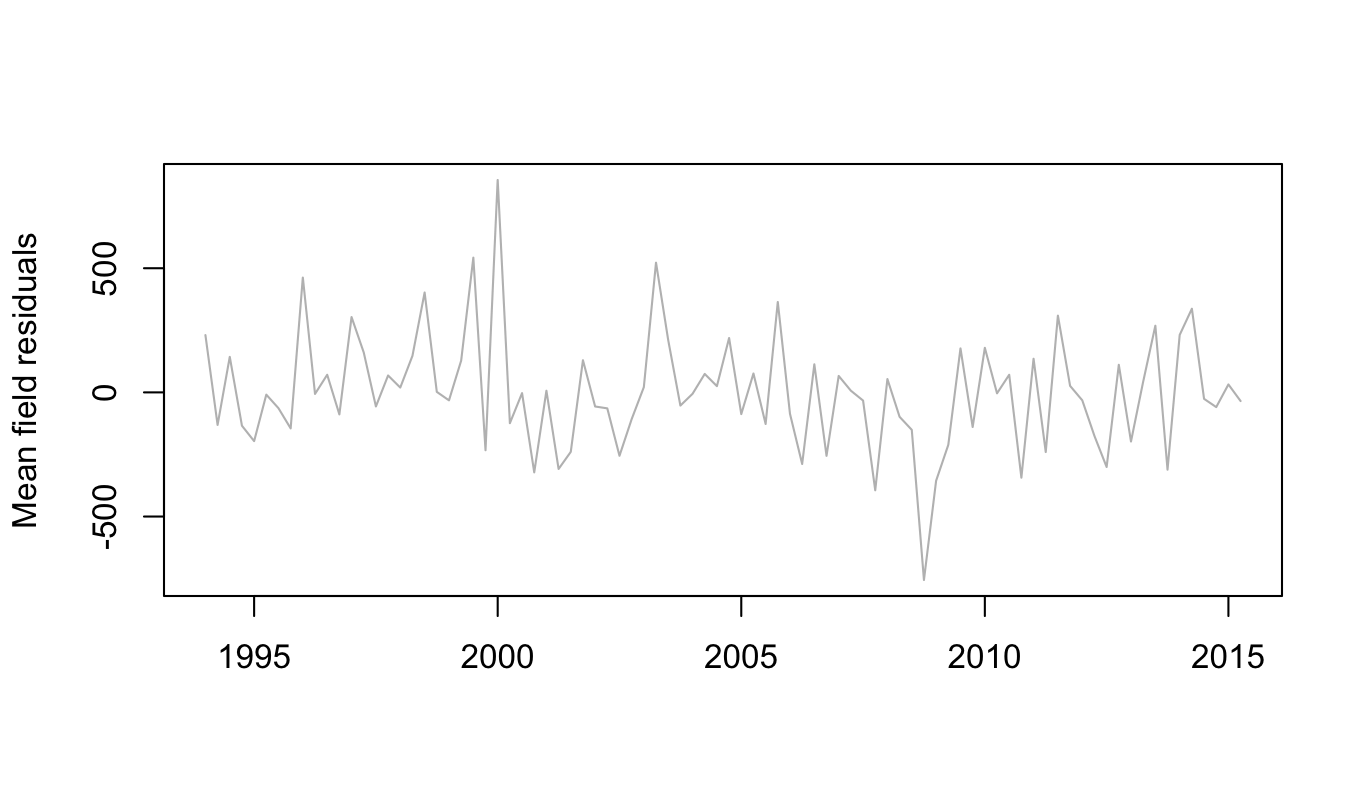}

}\hfill{}\subfloat[]{\includegraphics[width=0.45\textwidth]{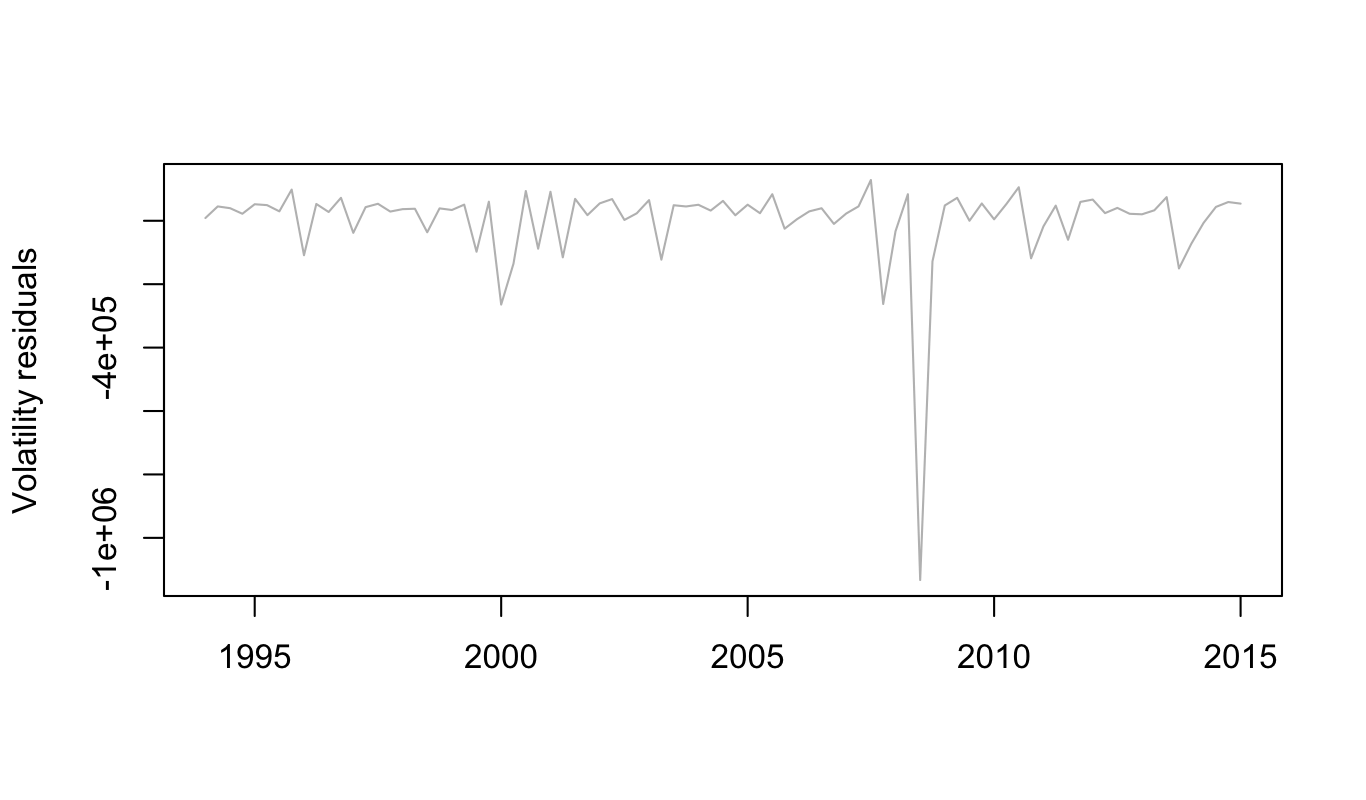}

}

\subfloat[]{\includegraphics[width=0.45\textwidth]{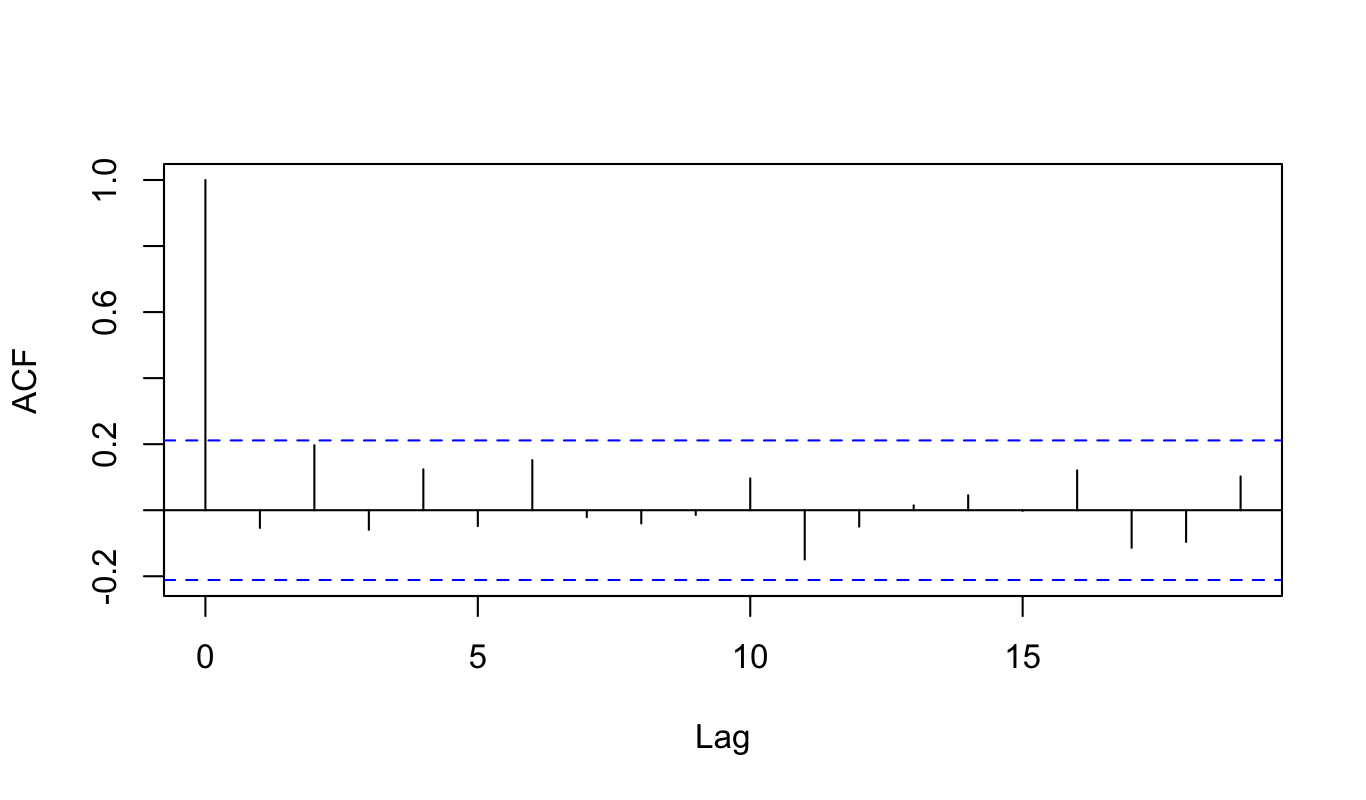}

}\hfill{}\subfloat[]{\includegraphics[width=0.45\textwidth]{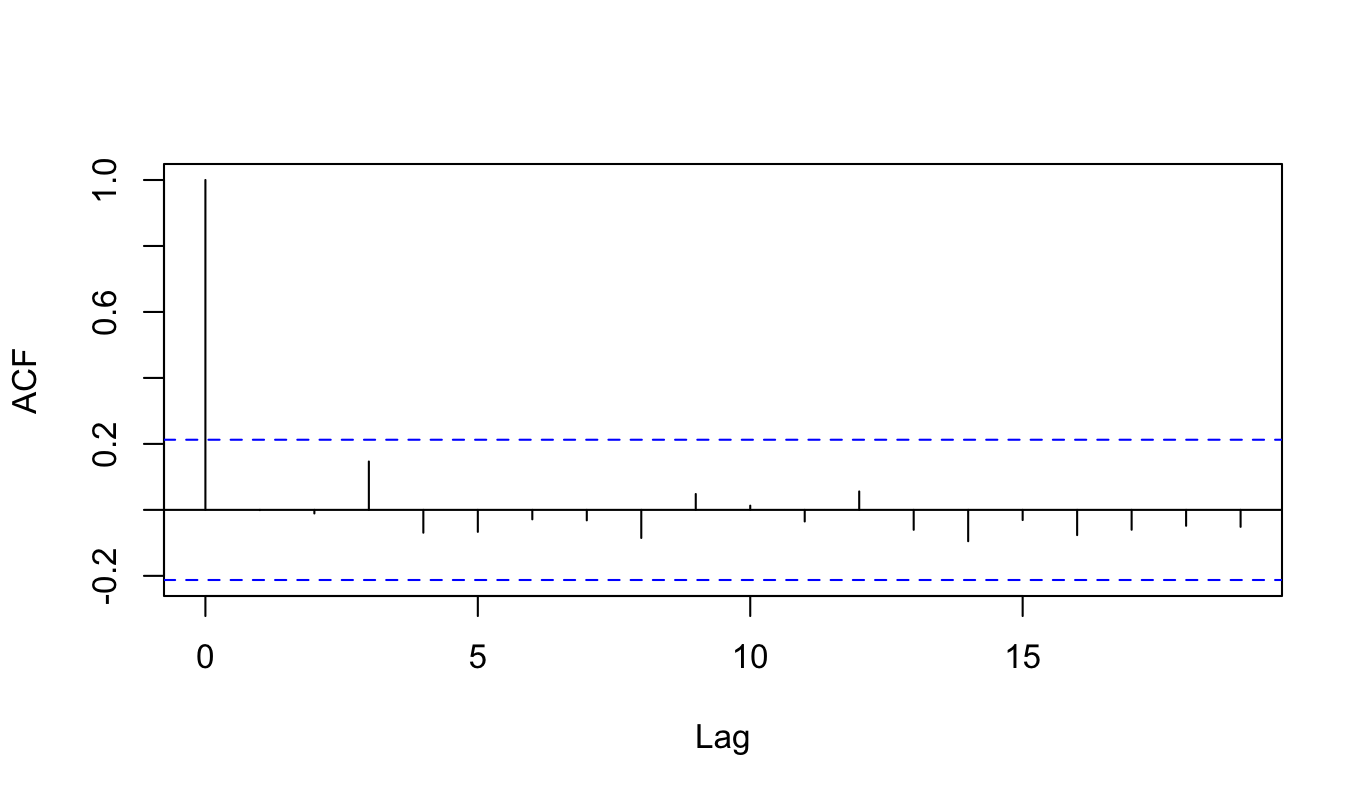}

}

\protect\caption{\label{fig:Residuals-in-Reduced-Growth}Residuals in Reduced Form
Estimate (Growth)}
\end{figure}

\begin{figure}
\subfloat[]{\includegraphics[width=0.45\textwidth]{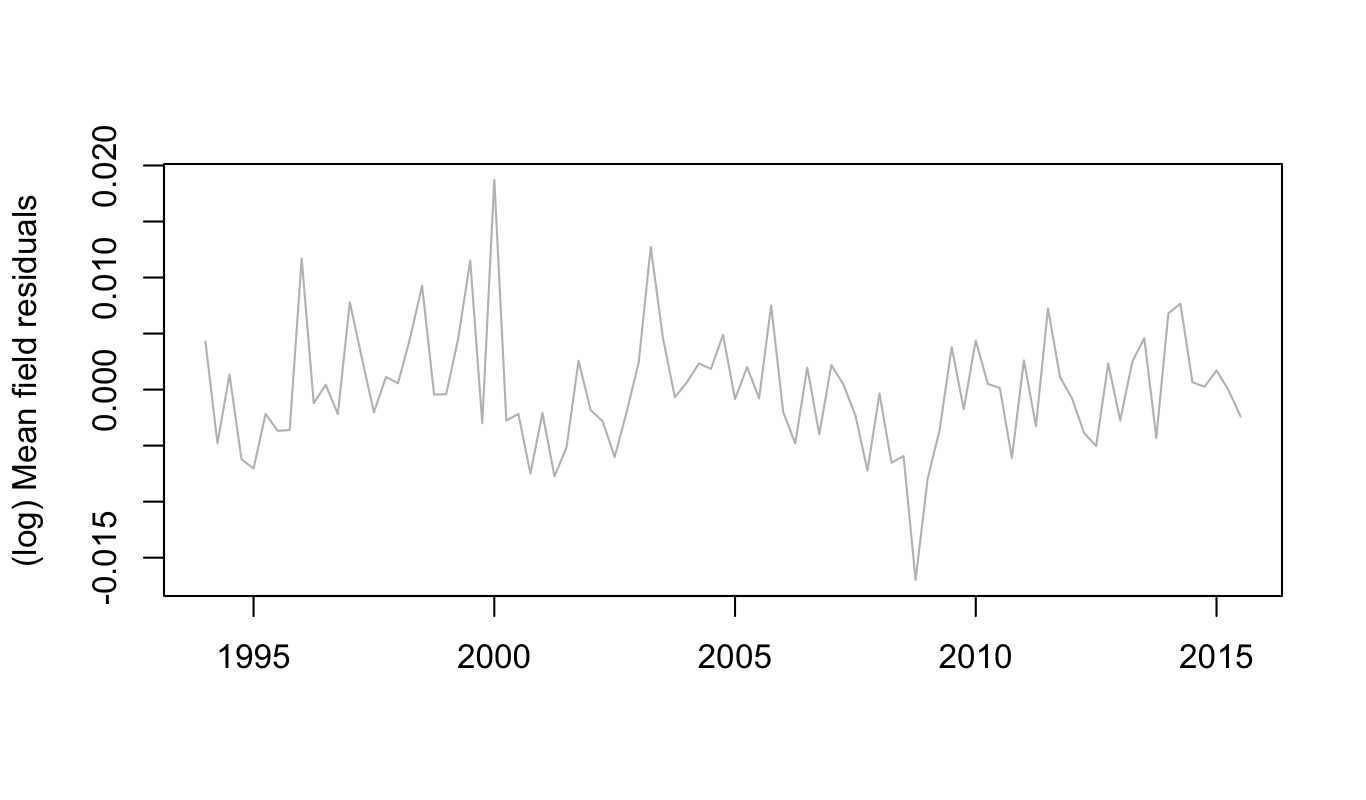}

}\hfill{}\subfloat[]{\includegraphics[width=0.45\textwidth]{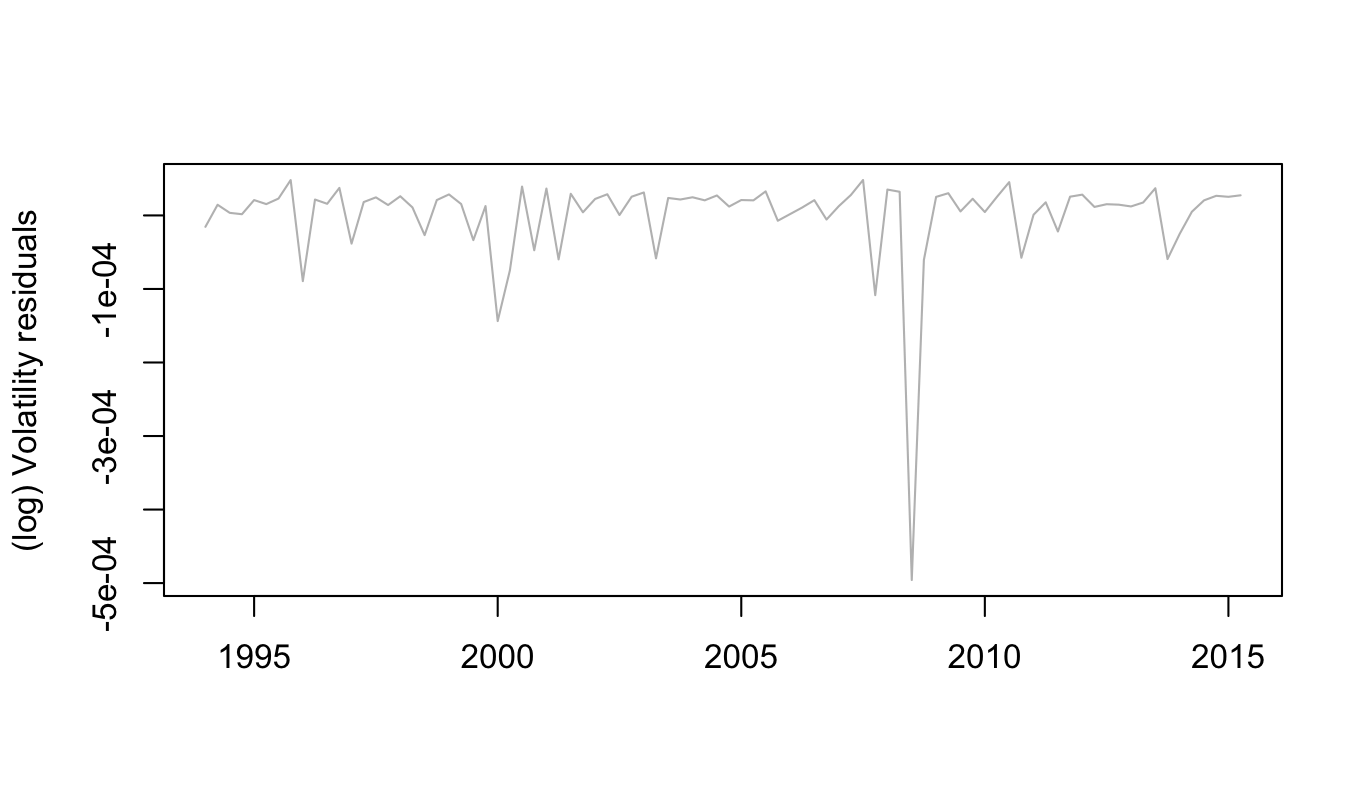}

}

\subfloat[]{\includegraphics[width=0.45\textwidth]{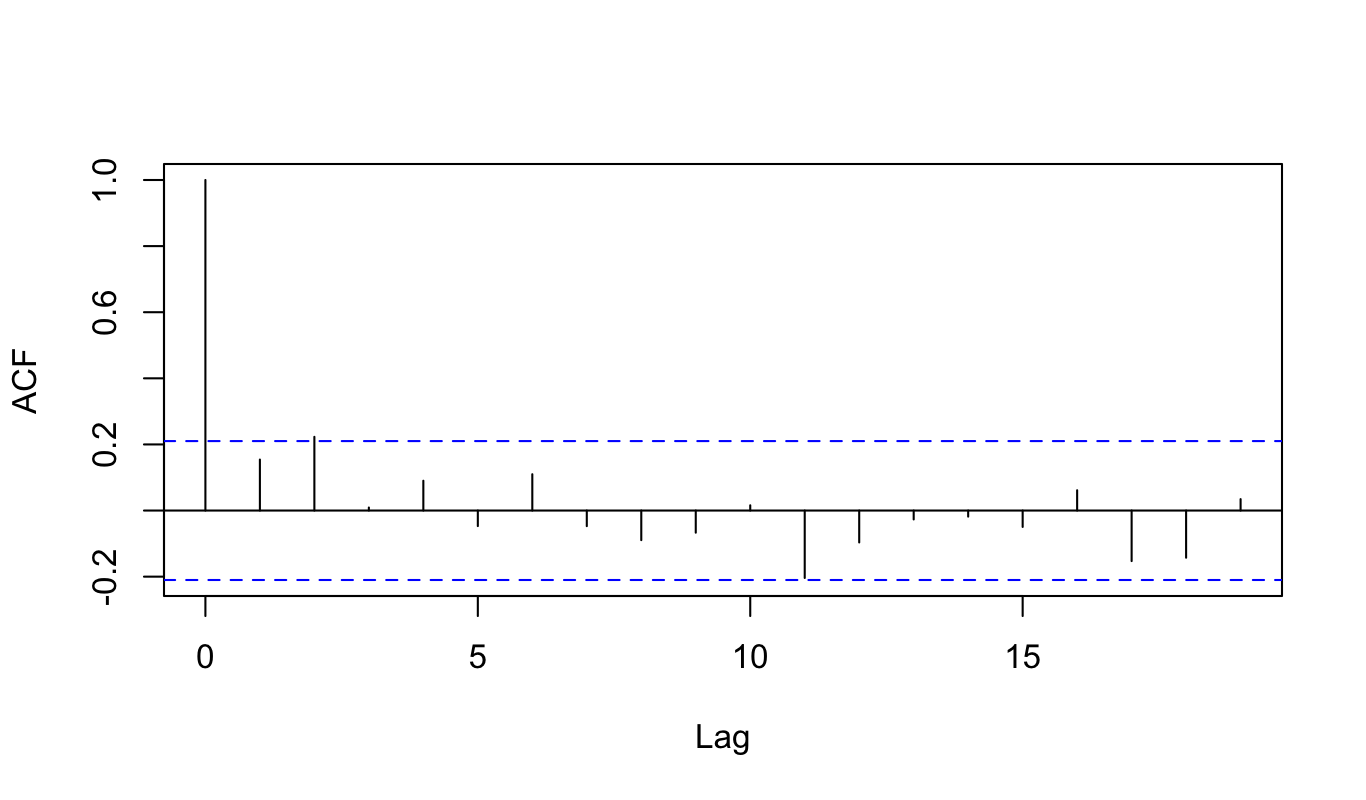}

}\hfill{}\subfloat[]{\includegraphics[width=0.45\textwidth]{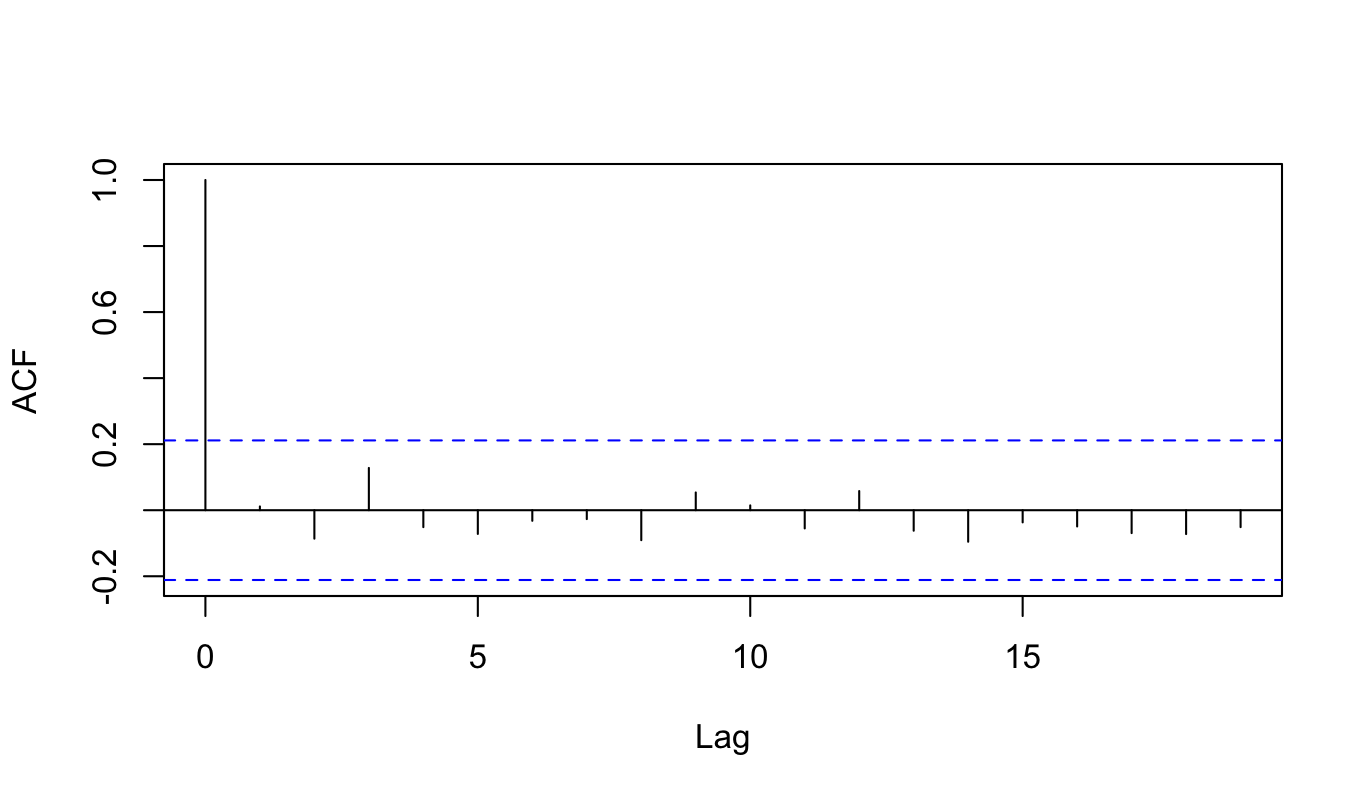}

}

\protect\caption{\label{fig:Residuals-in-Reduced-GrowthRate}Residuals in Reduced Form
Estimate (Growth rate)}
\end{figure}

Estimate results for for (\ref{eq:rf-1}) and (\ref{eq:rf-2}) are
given in Table \ref{Tab:RF}. The results contain two cases. One uses
GDP as the mean field growth variable $m_{Y}(t)$ the other one uses
logarithm of GDP as the mean field growth rate.%
\footnote{Strictly speaking, without a further justification, using logarithm
transform in a regression is quite vague, as in this transform non-linear
property disappear. If the mean field $m_{Y}(t)$ is transformed,
the expansion form of the system should be implemented using the transformed
expression. So the terms in (\ref{eq:rf-1}) and (\ref{eq:rf-2})
may be different from the original form (\ref{eq:reducedForm}). %
} From the table, we can see the coefficients of $\sigma_{Y}^{2}(t)$
in (\ref{eq:rf-1}) and (\ref{eq:rf-2}) are all significant which
means the mean field dynamics $dm_{Y}(t)$ is correlated with the
higher order information from its volatility $\sigma_{Y}^{2}(t)$.
\footnote{GARCH modules in \texttt{tseries} package of \texttt{R} provide either
non-convergent or insignificant estimates. Please refer to on-line
appendix for the output results of GARCH. %
} The residual plots for both cases are given in Figure \ref{fig:Residuals-in-Reduced-Growth}
and \ref{fig:Residuals-in-Reduced-GrowthRate} respectively. All residual
look stationary. However, the residuals in (\ref{eq:rf-2}) some dependent
structure remains. This may be due to the large perturbation around
$2008$. 

For structural estimation of (\ref{eq:StructuralEq}), some specification
of non-linear function $g(\cdot)$ is pre-required. For simplicity,
it is assumed that $g(\cdot)$ is an exponential function and (\ref{eq:MasterEqGeneral})
is stationary so that $\sum_{i=1}^{N}X_{t}(\omega^{i})\sim\mbox{Gamma}(\alpha_{t},\beta_{t})$.
The income distribution therefore is a stationary outcome but the
parameters $(\alpha_{t},\beta_{t})$ as endogenous variables adjust
to the new level each year due to the interaction and evolution effects
of the aggregation. The estimation form is given as follows
\begin{align}
\mathbb{E}\left[\sum_{i=1}^{N}X_{t}(\omega^{i})\right] & =\frac{\alpha_{t}}{\beta_{t}}=\mbox{Cons}+\mbox{Coef}_{1}\times\ln m_{Y}(t),\label{eq:sf-1}\\
\ln m_{Y}(t+1)-\ln m_{Y}(t) & =\theta(t)+\varepsilon_{1,t}+\frac{\varepsilon_{2,t}}{\sum_{i=1}^{N}X_{t}(\omega^{i})},\quad\varepsilon_{i,t}\sim GWN(0,1).\label{eq:sf-2}
\end{align}
The assumption of Gaussian white noises for $\varepsilon_{1,t}$ and
$\varepsilon_{2,t}$ and the simplification of $e_{t}/m_{Y}(t)$ in
(\ref{eq:LogMean}) are purely due to technical purposes. Because
assuming linear additive Gaussian structure can give us a feasible
filtering algorithm. In the estimation, we use the forward filtering
part of the Baum-Welch formula. For details, please see Appendix \ref{sub:Forward-Filter}.

\begin{figure}

\subfloat[]{

\includegraphics[width=0.3\textwidth]{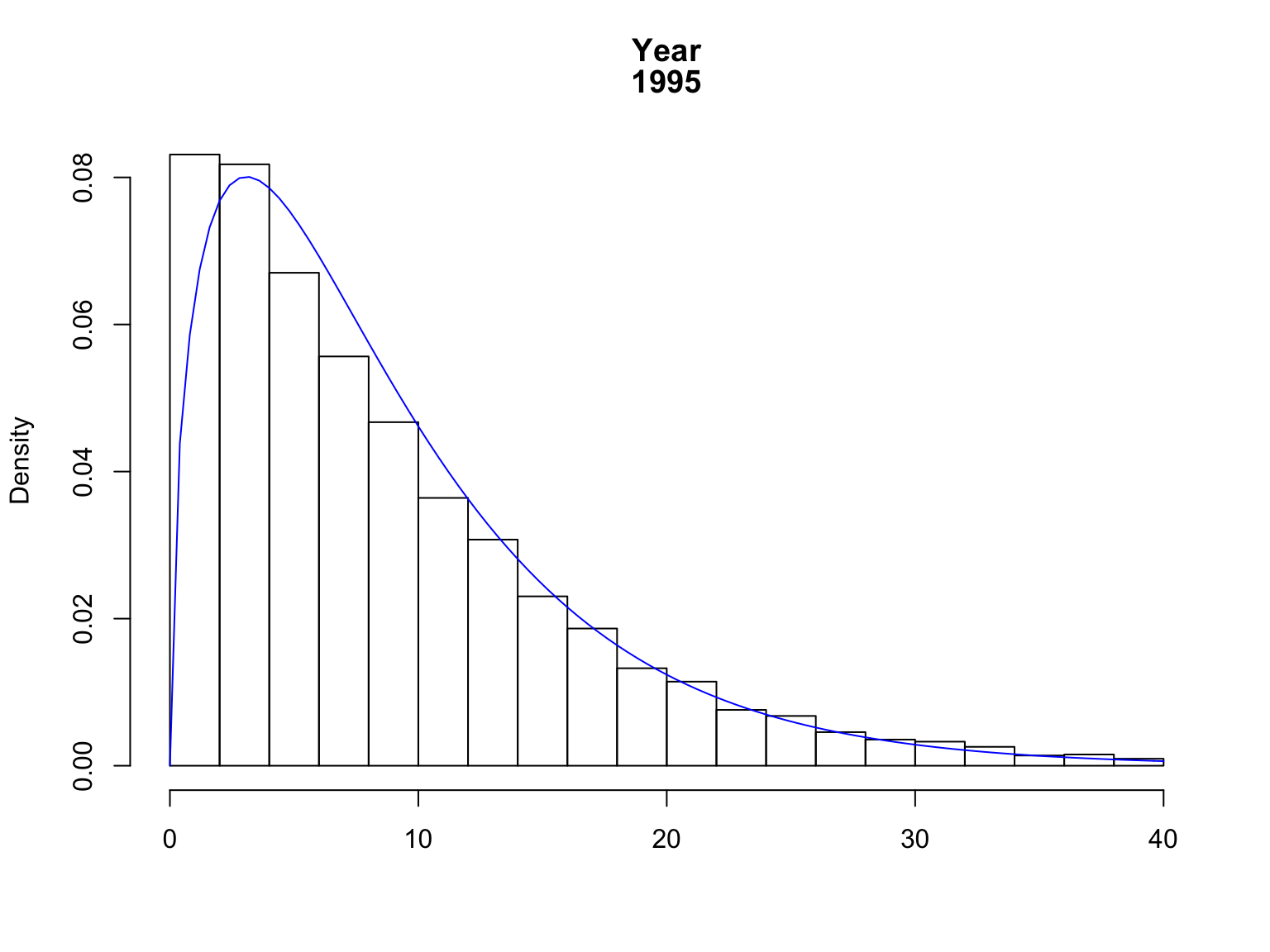}

}\hfill{}\subfloat[]{\includegraphics[width=0.3\textwidth]{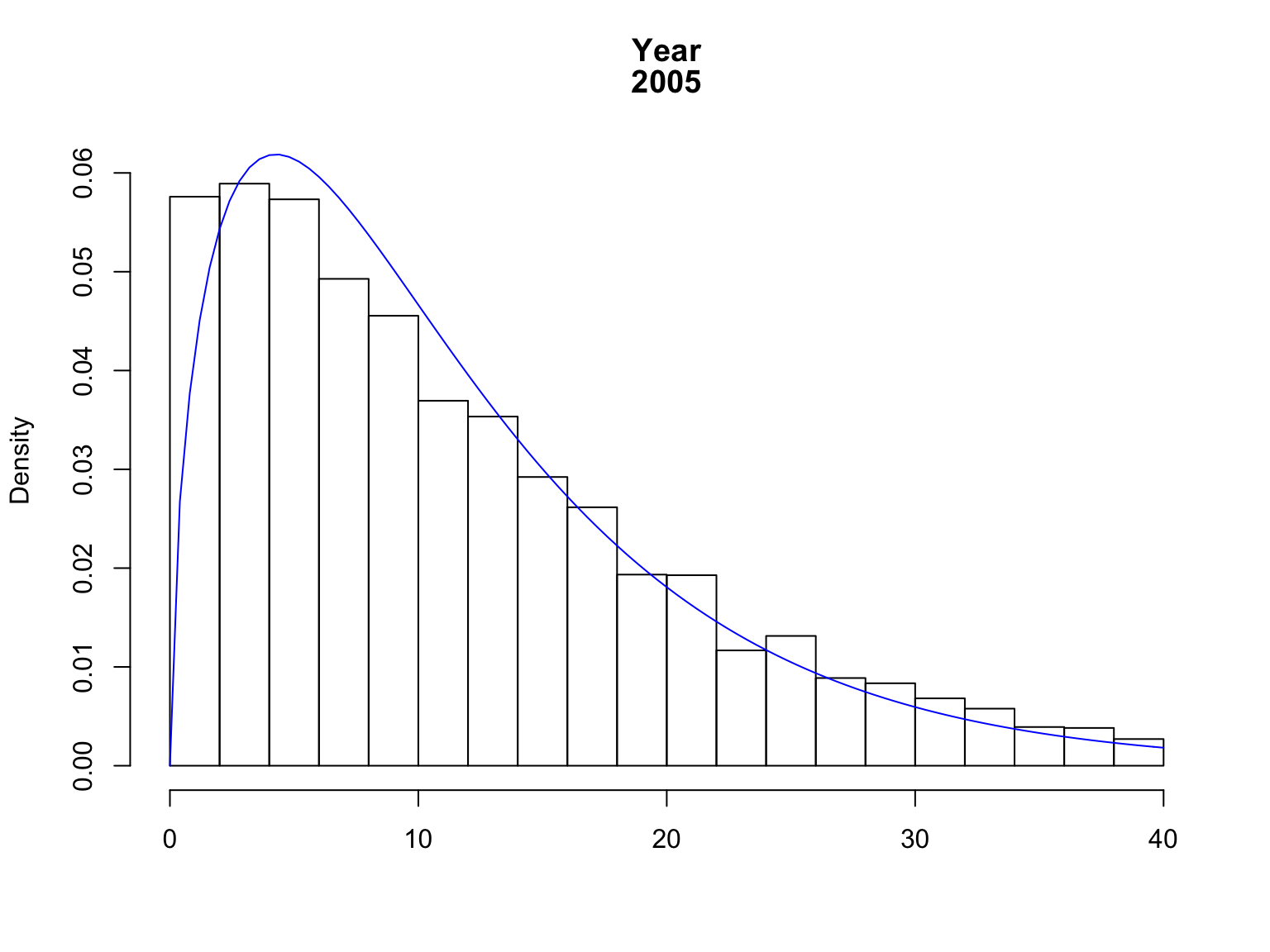}

}\hfill{}\subfloat[]{\includegraphics[width=0.3\textwidth]{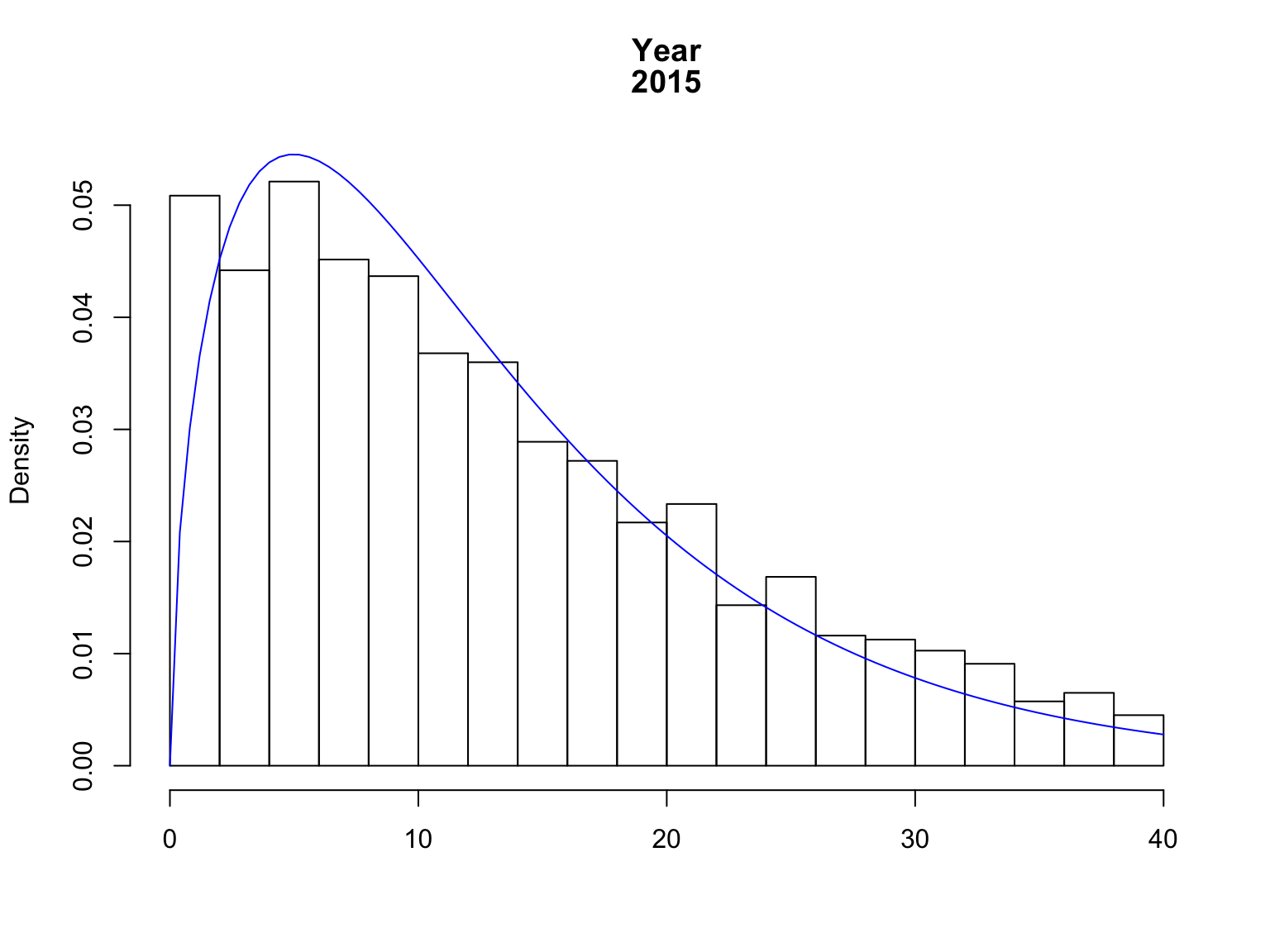}

}

\subfloat[]{\includegraphics[width=0.3\textwidth]{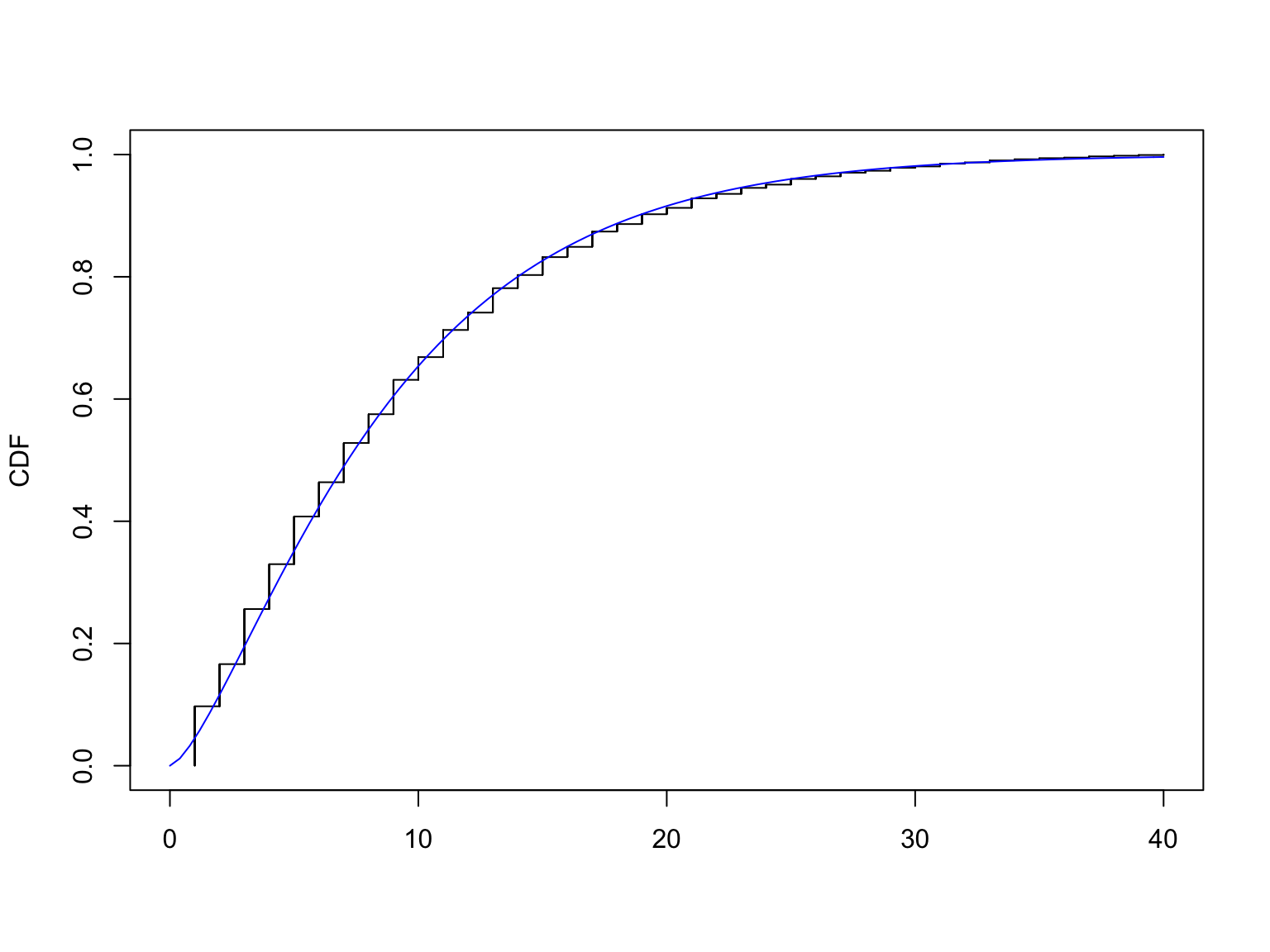}

}\hfill{}\subfloat[]{\includegraphics[width=0.3\textwidth]{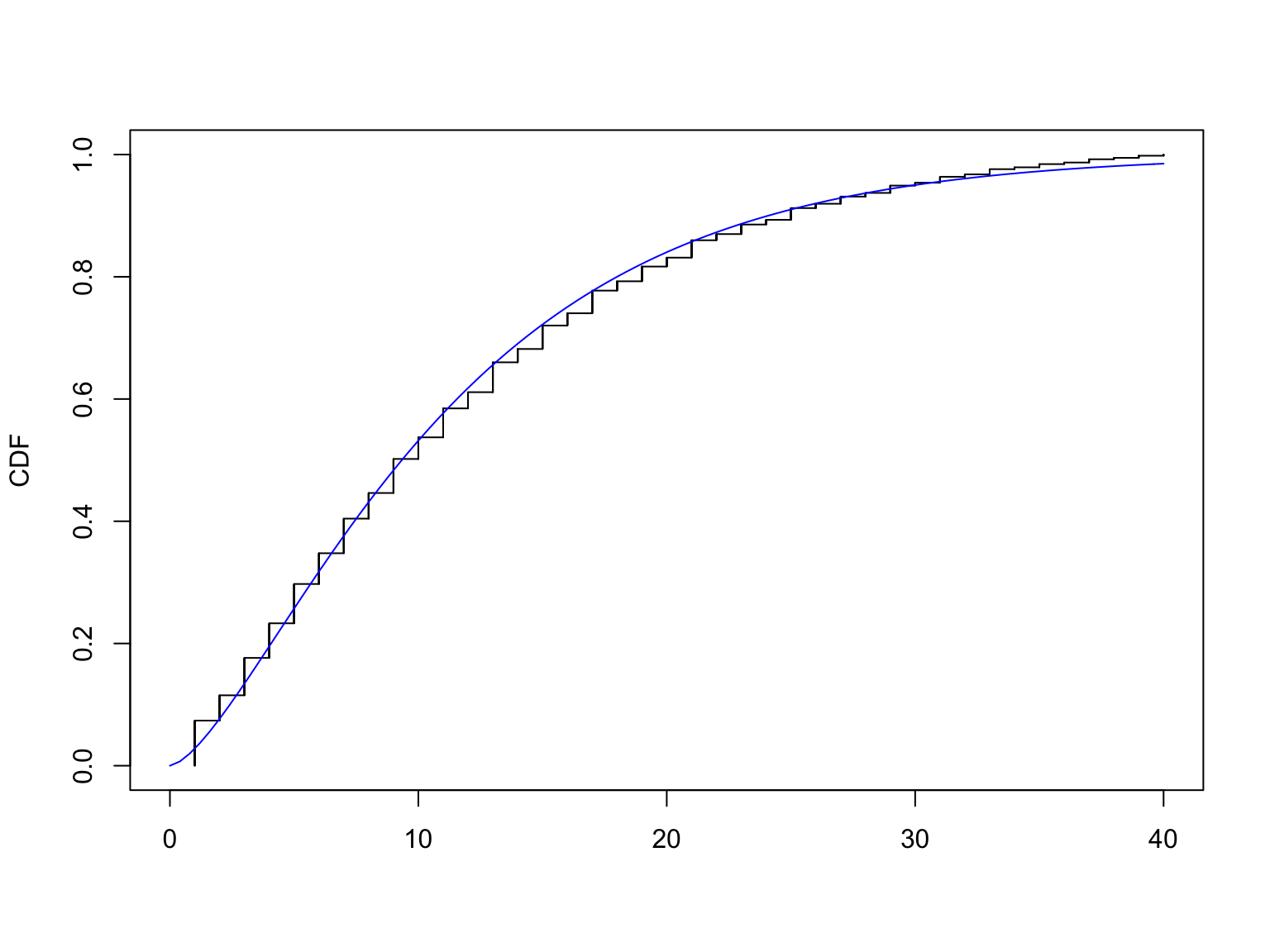}

}\hfill{}\subfloat[]{\includegraphics[width=0.3\textwidth]{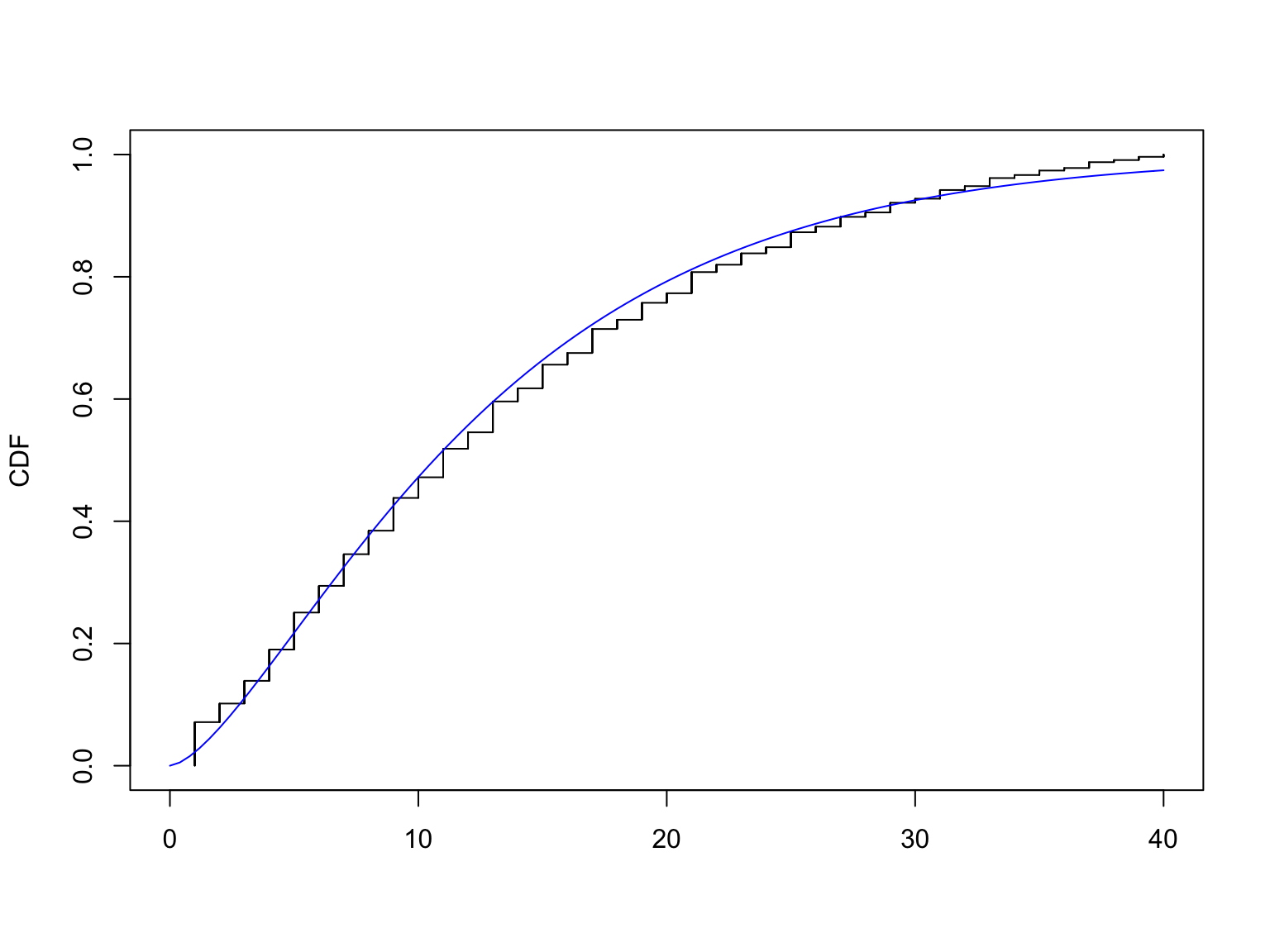}

}

\protect\caption{\label{fig:Estimates-of-Income}Estimates of Income Distribution by
Gamma$(\alpha_{t},\beta_{t})$ }

\end{figure}

\begin{figure}
\subfloat[]{\includegraphics[width=0.45\textwidth]{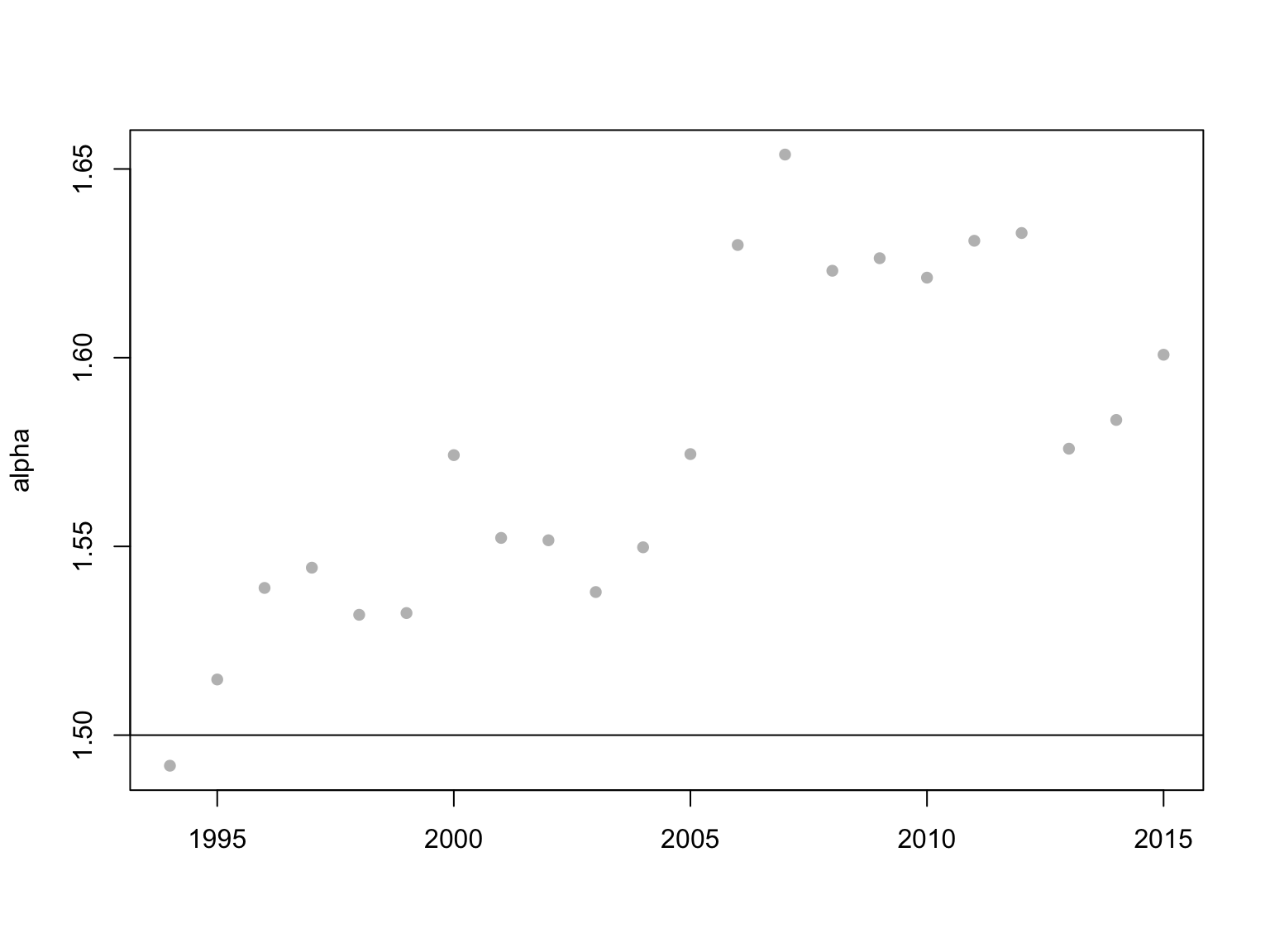}

}\hfill{}\subfloat[]{\includegraphics[width=0.45\textwidth]{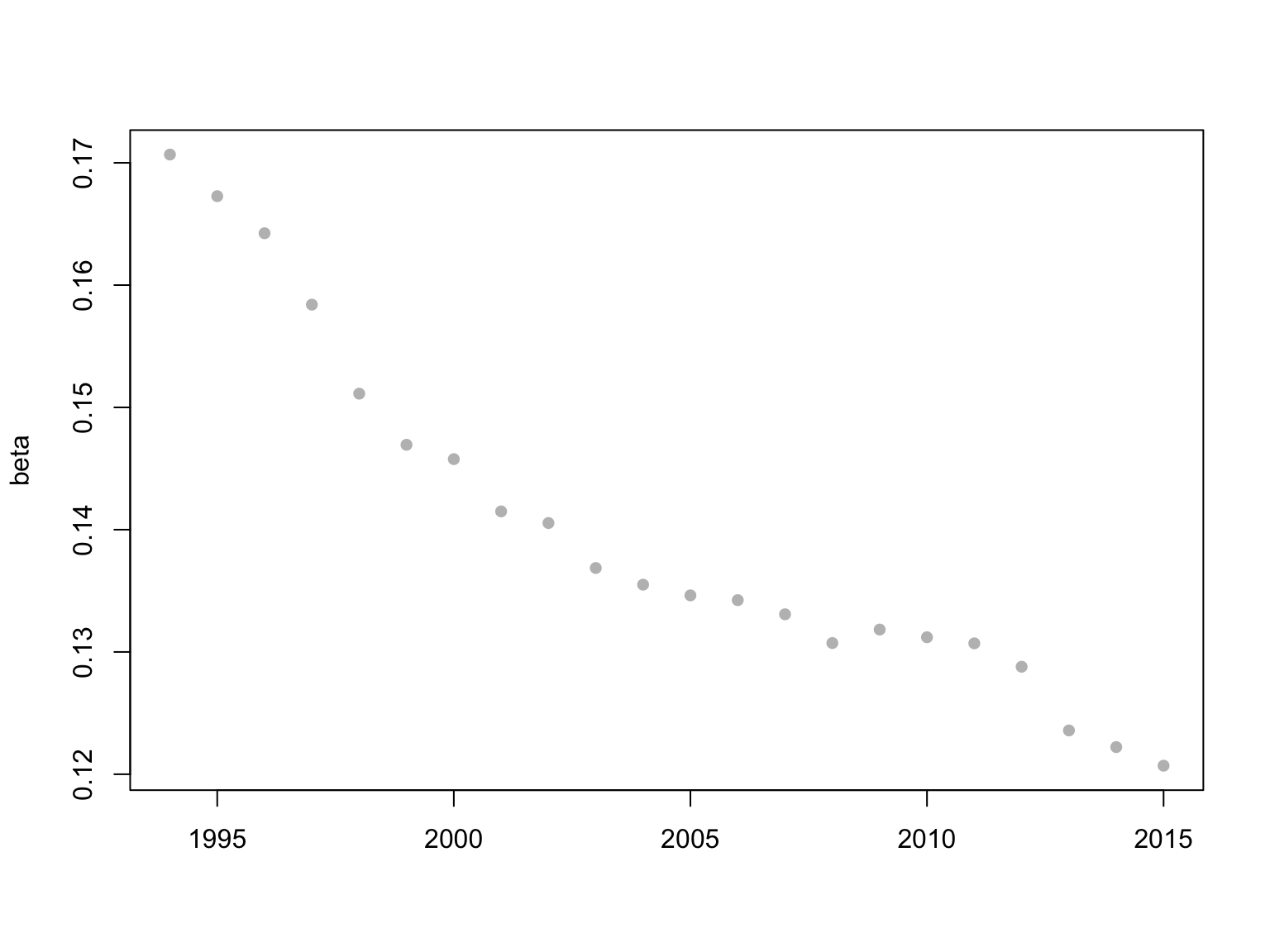}

}

\protect\caption{\label{fig:alpha-beta-94-05}$(\alpha_{t},\beta_{t})$ in $1994-2015$ }
\end{figure}

\begin{figure}
\subfloat[]{\includegraphics[width=0.45\textwidth]{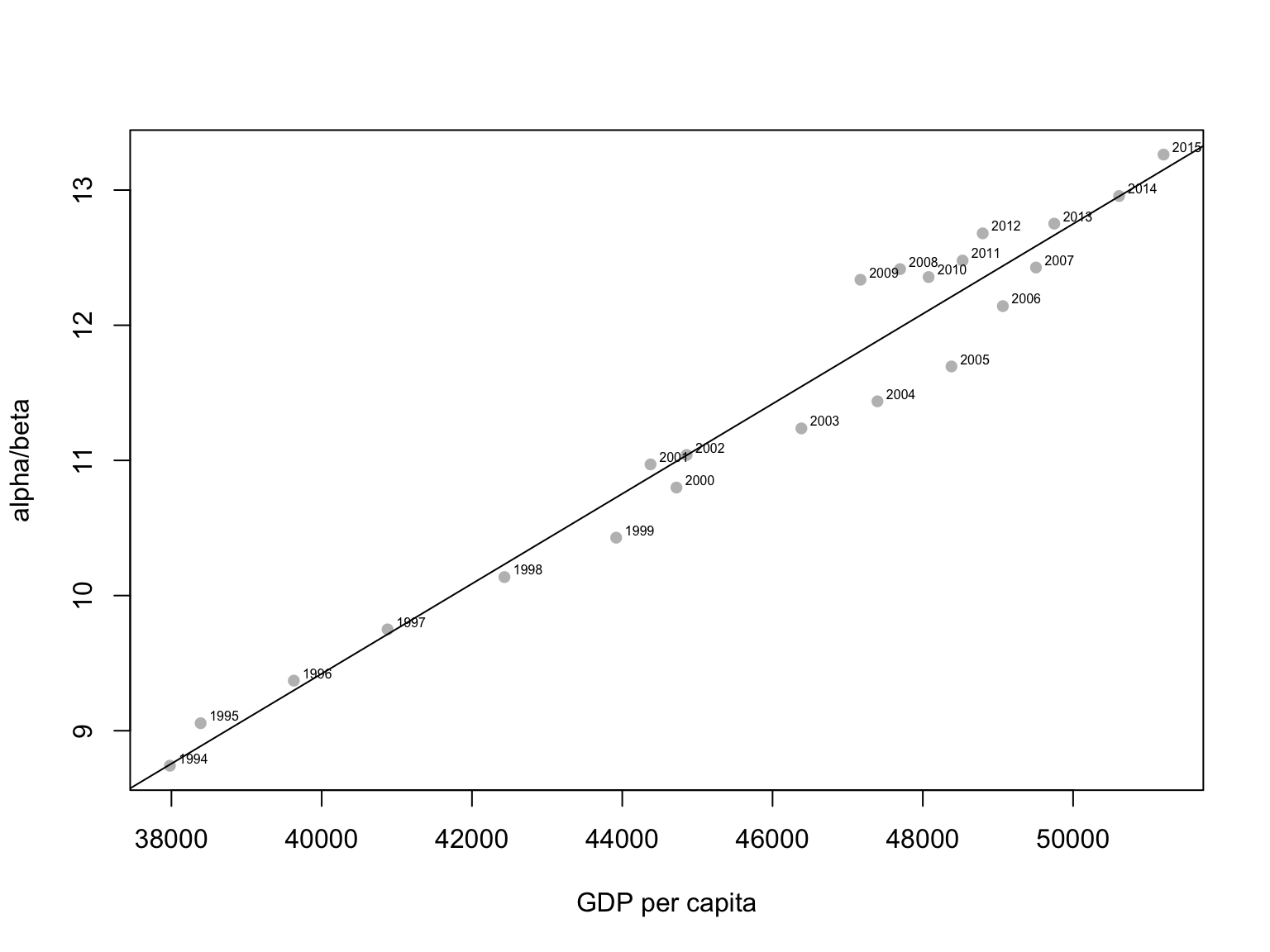}

}\hfill{}\subfloat[]{\includegraphics[width=0.45\textwidth]{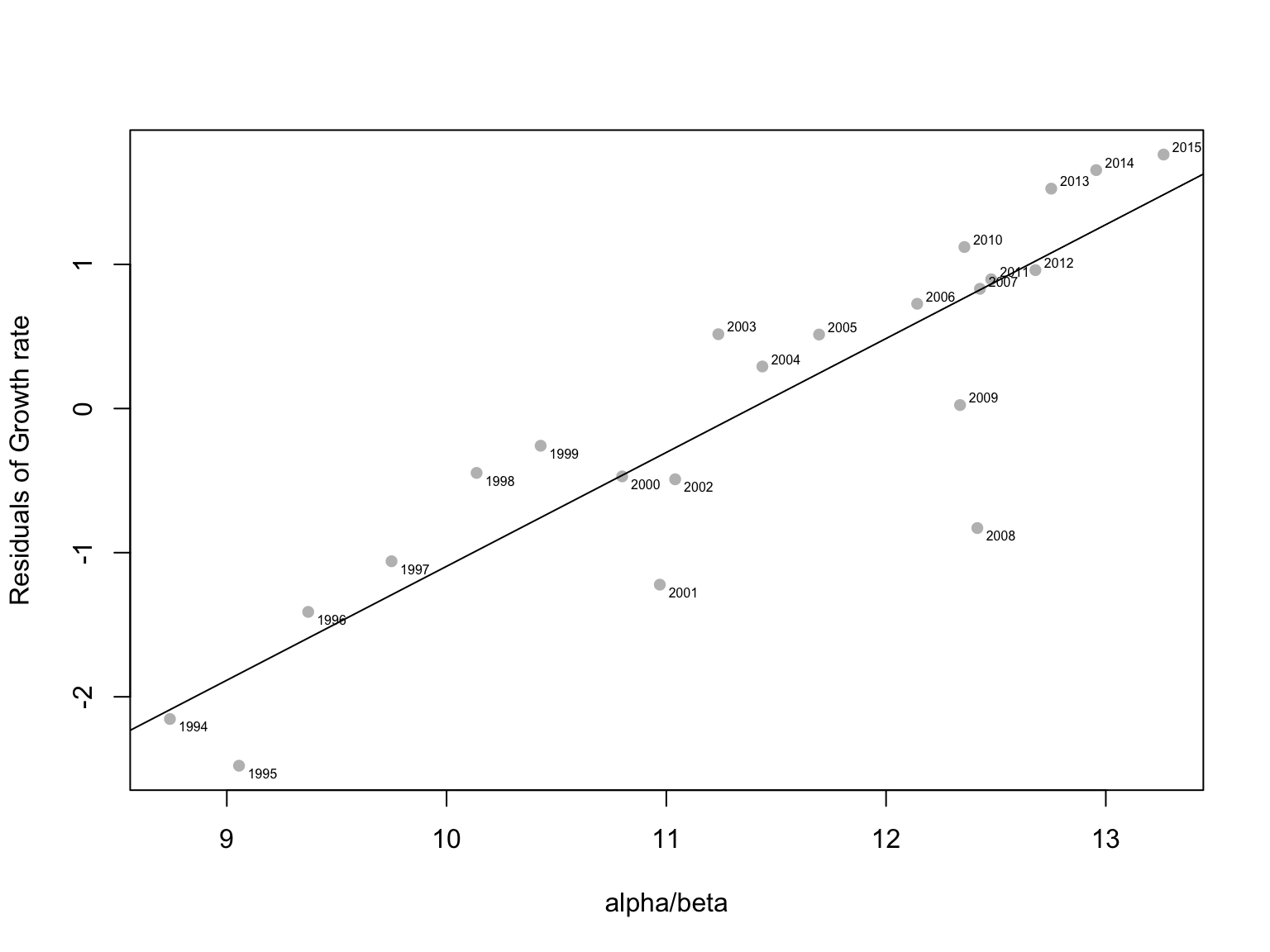}

}

\protect\caption{\label{fig:Dynamics-of-alpha-beta}Dynamics of $\alpha_{t}/\beta_{t}$
and GDP }
\end{figure}

The full estimation is divided into two steps. First, we need to estimate
$(\alpha_{t},\beta_{t})$ using income data from $1994$ to $2015$.
Figure \ref{fig:Estimates-of-Income} shows three selective years
for illustration. The plots show the Gamma specification fits well
with both density and cumulative distribution function. The plots
for other years provide similar results that can be found in the online
appendix. Dynamical patterns of $(\alpha_{t},\beta_{t})$ can be found
in Figure \ref{fig:alpha-beta-94-05}. The trend of $\alpha_{t}$
is increasing but it has fluctuations around late $1998$, $2003$,
$2008$ and $2012$ which correspond to early $2000$s recession and
financial crisis of $2008$. It is surprising that the drop at $2012$
has a bigger scale than that of $2008$. The trend of $\beta_{t}$
is monotone decreasing. It drops about $30\%$ (from $0.17$ in $1994$
to $0.12$ in $2015$). According to this calculation, the scale of
$\beta$ decreasing contributes more to the growth than the increase
of $\alpha$. This argument can be verified in Figure \ref{fig:Dynamics-of-alpha-beta}(a)
where it shows an almost linear trend between GDP and $\alpha_{t}/\beta_{t}$,
the mean of Gamma distribution. U.S. GDP keeps growing since $1994$.
However if we compare this result with the trends of $\alpha_{t}$
and $\beta_{t}$, we can see the growth may be mostly from the decline
of $\beta$ value. From the arguments in Section \ref{sub:Remarks-StruAgg},
we can see these phenomena may imply an enlarged inequality happening
in U.S.. 

\begin{figure}
\subfloat[]{\includegraphics[width=0.45\textwidth]{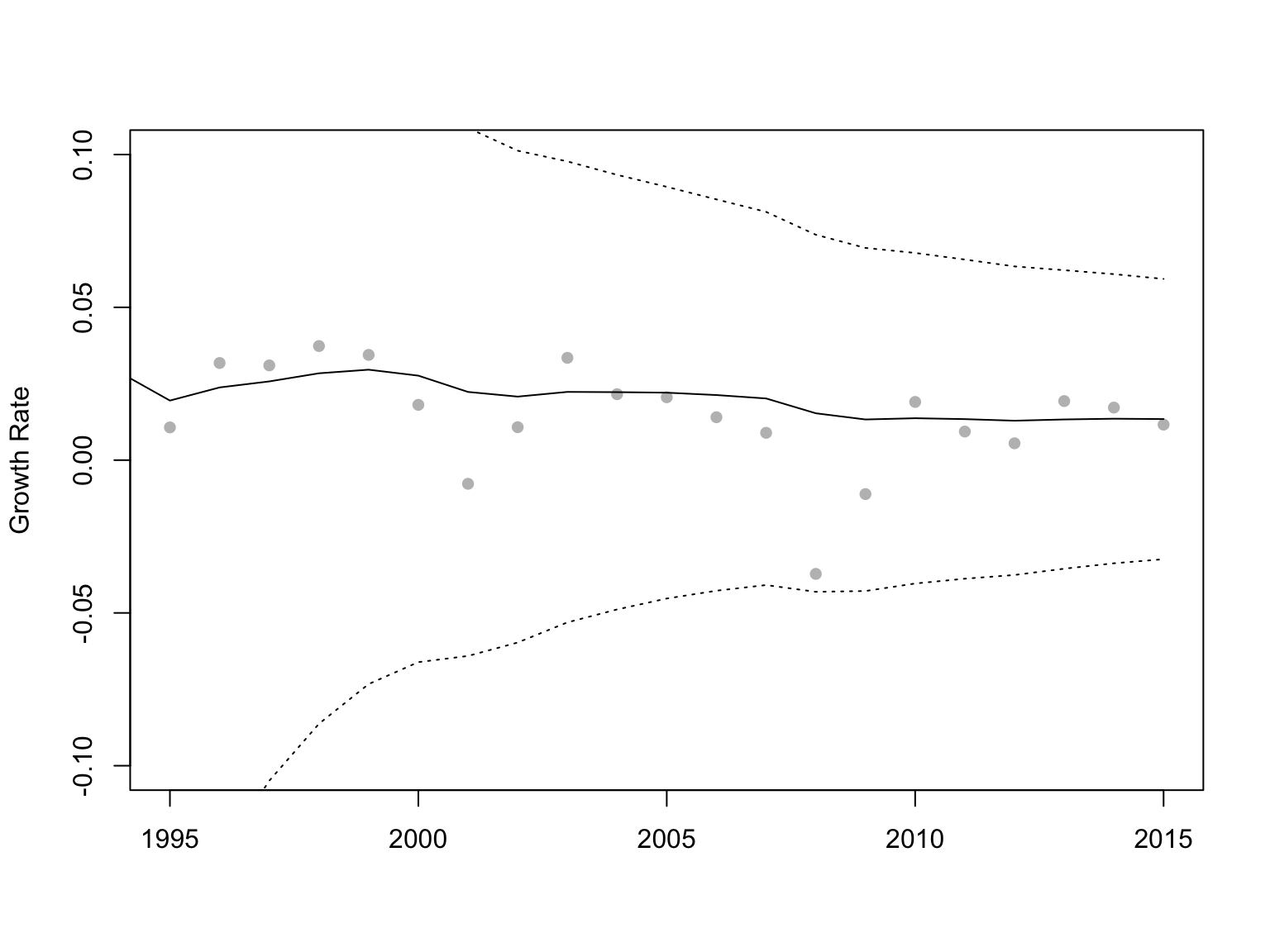}

}\hfill{}\subfloat[]{\includegraphics[width=0.45\textwidth]{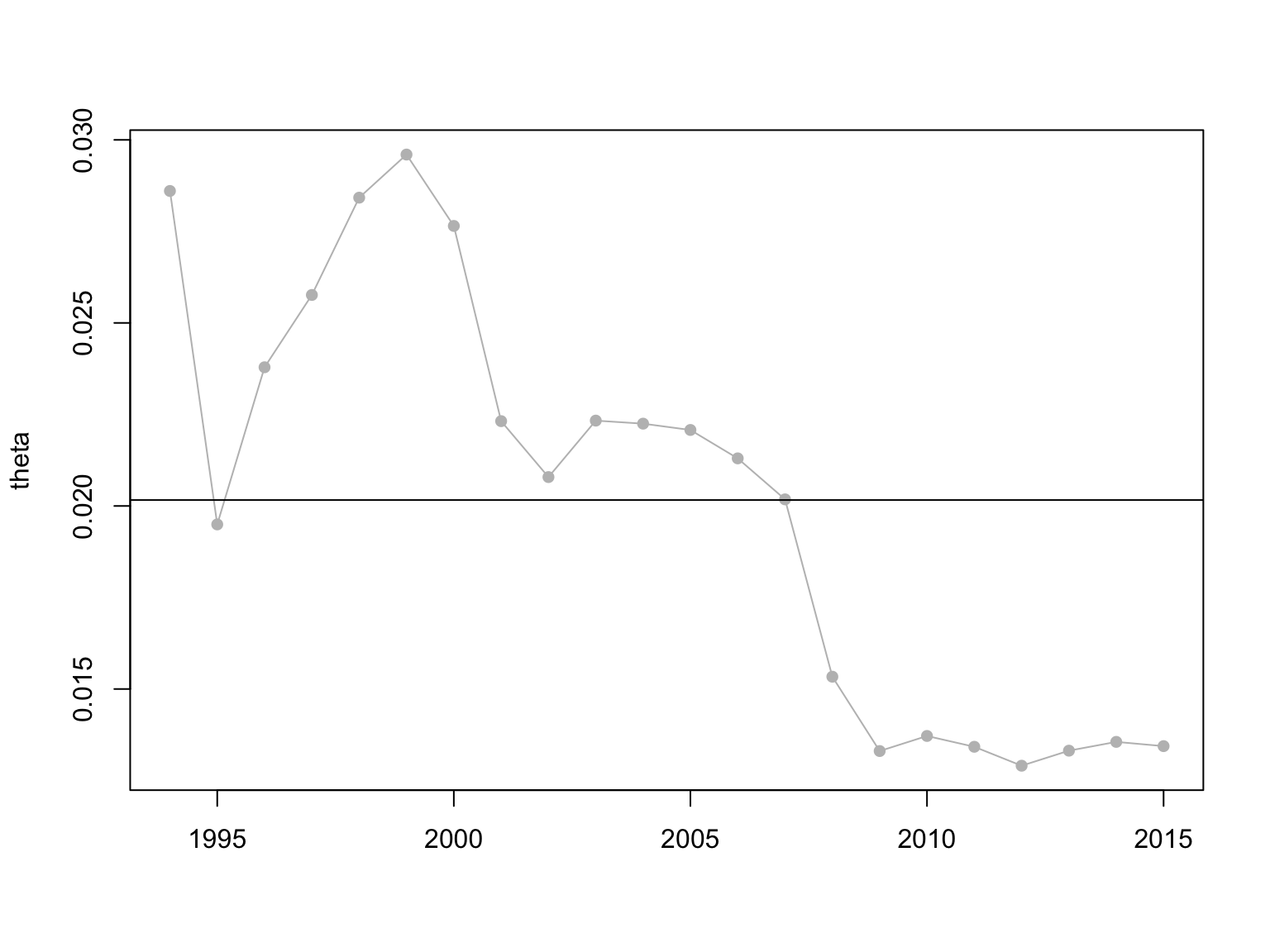}

}

\subfloat[]{\includegraphics[width=0.45\textwidth]{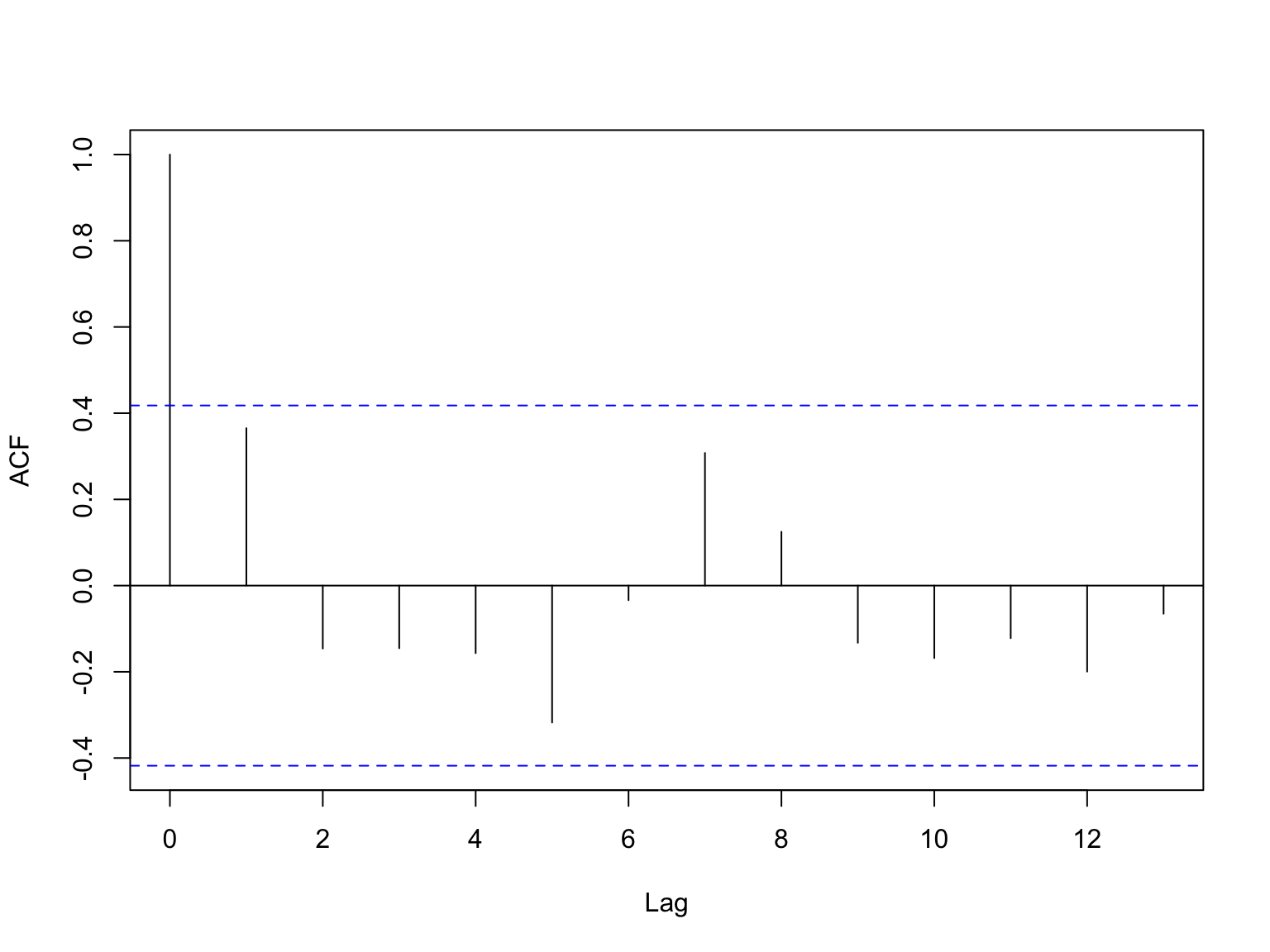}

}\hfill{}\subfloat[]{\includegraphics[width=0.45\textwidth]{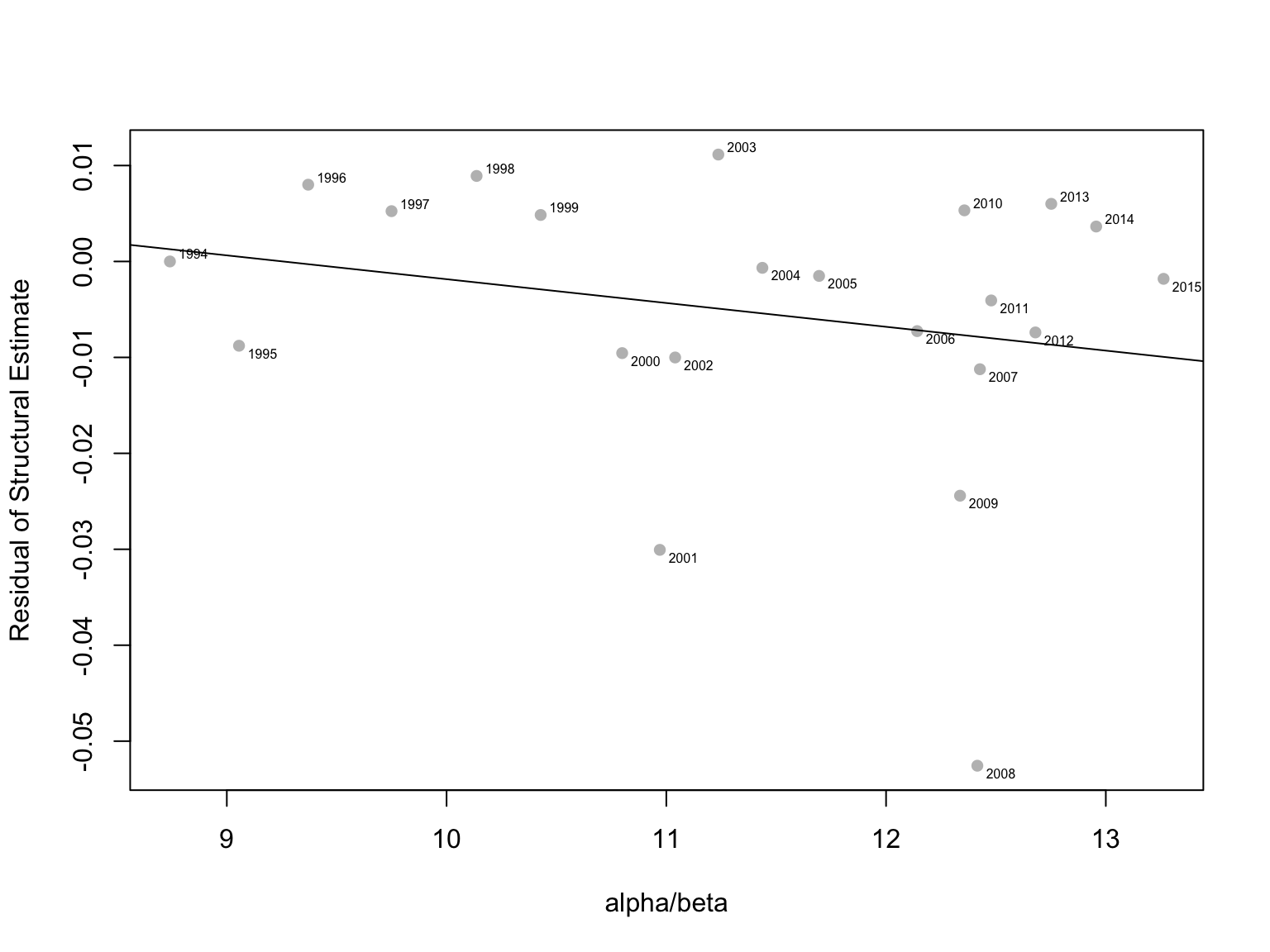}

}

\protect\caption{\label{fig:Filters}Structural Estimate}
\end{figure}

With estimated $(\alpha_{t},\beta_{t})$, the second step is to estimate
equations (\ref{eq:sf-1}) and (\ref{eq:sf-2}). The results are given
in Table \ref{Tab:SF}. The first result is from a regression of $m_{Y}(t)$
on $\alpha_{t}/\beta_{t}$. It shows a significant relation that coincides
with Figure \ref{fig:Dynamics-of-alpha-beta}(a). Equation (\ref{eq:sf-1})
also has a significant coefficient. However, $\alpha_{t}/\beta_{t}$
is shown to have a significant correlation with the residual of (\ref{eq:sf-1})
that implies an endogenous issue of (\ref{eq:sf-1}). Figure \ref{fig:Dynamics-of-alpha-beta}(b)
also illustrates this problem. With this concern, we give a structural
estimate of equation (\ref{eq:sf-2}). First we filter out the effect
of $\theta(t)$ then by using the information of the variance $\alpha_{t}/\beta_{t}^{2}$
of Gamma distribution, we remove the heteroskedasticity effect. The
residual of this estimate as shown in the in Table \ref{Tab:SF} is
uncorrelated with $\alpha_{t}/\beta_{t}$. Figure \ref{fig:Filters}(a),
(c), (d) show that forward filter estimate has captured the heteroskedasticity
and that its residual is stationary and uncorrelated with $\alpha_{t}/\beta_{t}$. 

The estimated values of time varying parameter $\theta_{t}$ can be
found in Figure \ref{fig:Filters}(b). By extracting high order information
interacting with the income distribution, dynamical pattern of $\theta$
gives us a very different picture about recent growth in the U.S..
As we have seen in Figure \ref{fig:alpha-beta-94-05}, the growth
of stationary aggregates $\alpha_{t}/\beta_{t}$ is mainly intrigued
by reducing $\beta_{t}$. The decrease of $\beta$ contributes to
the GDP growth thorough the non-linear aggregation. Once we extract
these stationary and non-stationary effects, the filtered growth rate
of U.S. GDP is significantly decreasing after $2003$ and it remains
at a relatively low level. Although Figure \ref{fig:Filters}(b) gives
a different dynamical patterns from the existing figures of U.S. GDP
growth rate, the mean of this structural parameter, $0.02$, is close
to the averaged official figures.%
\footnote{By the second equation of (\ref{eq:StructuralEq}), we can see that
$\theta_{t}$ is compatible with the role of the coefficient of $m_{Y}(t)$
in the reduced form equation (\ref{eq:rf-1}). The estimated value
of this reduced form coefficient is $0.201$ which is about ten times
of the mean of $\theta_{t}$. However, one should keep in mind that
volatility of mean fields contributes a negative effect in the reduced
form.%
} This closeness results from the economic expansion in late 90s.

Overall, these empirical results give interesting alternative interpretations
to the recent development in the U.S.. We need to emphasize that these
empirical estimates have make some simplifications from the theoretical
models. These simplifications so far have not generated unrobustness.
But we expect more sophisticated approaches to improve these estimates.

\section{\label{sec:Bird's-eye-View}Bird's-eye View}

\subsection{Deterministic and Stochastic Economies}

Section \ref{sec:A-non-abstract-model} and \ref{sec:Abstract-Models}
discuss two different types of economies. It is obvious that the stochastic
economy is more realistic and closer to the economy where we live.
On the other hand, deterministic economy describes an illusion that
is closer to the utopian images where there is no unequal status and
people produce ample goods. If we want our society to evolve to a
utopian state, can this deterministic economy be the destination of
this evolution sequence?

To reach this destination, the deterministic economy needs to reach
the same level of production as the stochastic one. That is to say,
the deterministic individual growth function $f(\cdot)$, need to
produce as much as the stochastic counterpart $f^{(a)}(\cdot)$ does.
Unfortunately, the functional class of $f$ is smaller than that of
$f^{(a)}$ which means $f$ cannot incorporate with some growth patterns
that are feasible under $f^{(a)}$. In fact, $f$ can be a special
case of $f^{(a)}$ where $a$ is a deterministic sequence indexing
every period of time. Even for deterministic function, $f$ can not
be embedded in many interesting cases. For example, some deterministic
logistic maps or cobweb models can be ruled out from the class of
$f$, because they can induce infinite many bifurcations that give
the randomness to the growth. 

Another important point is that in the evolution process variation
happens in an uncertainty way, at least, to our best understanding
so far. Innovations and their external effects are rather crucial
to our growth and evolution but they are unpredictable. Thus instead
of viewing deterministic economy as a utopia limit of our evolution,
it is better to think it as the garden of Eden, the initial state
of a stochastic economy from which the uncertainty appeared.

\subsection{Nonlinearity }

Section \ref{sec:Aggregation} points out the difference between stationary
and non-stationary aggregates. As non-stationarity comes from nonlinear
aggregation functions, one concern is the role of non-linearity. It
is well known that many nonlinear patterns can be found in population
growth, urbanization, social interactions, etc. All these nonlinear
processes in some senses seem to increase our productivities. Meanwhile,
these processes also cause many problems.

The same dilemma can be found in the model in Section \ref{sec:Aggregation}.
It shows that uncertainty level is propagated by nonlinearity. Volatility
or even higher order information may contribute to the final aggregation.
Laws of individual or stationary aggregates are distorted during the
nonlinear aggregation. One can infer that the higher the nonlinearity,
more severe the distortion is. It is not so clear the total effect
of this type of distortion. We know the innovation is highly unpredictable,
so we tend to believe that the likelihood of innovation is higher
in a nonlinear agglomeration. But the occurrences of war and crisis
also reflect some nonlinear features, can we believe that non-stationary
aggregations also conspire these events? If non-stationary aggregation
accompanies with both angels and evils, it is our duties to prevent
the worst case scenarios.

\subsection{$\alpha$-Growth and $\beta$-Growth }

Section \ref{sec:Structural-Aggregation} uses $(\alpha,\beta)$ to
connect individual probabilistic laws and aggregates'. It points out
the different roles of $\alpha$ and $\beta$ in determining inequality.
The growth driven by $\alpha$ tends to create a more equal economy
than the one driven by $\beta$. A policy maker should be in favor
of proposals that may increasing $\alpha$ rather than decreasing
$\beta$.

However, $\alpha$ and $\beta$ are deeper parameters that can be
embedded in people's characters (risk loving or risk aversion), social
norms (conservative or liberal), or even culture and religion (creative
and independent or stable and united). The direct policy on income
or wealth may hardly affect these intrinsic characters making up of
$(\alpha,\beta)$. For example decreasing $\beta$ can be an unconscious
action as shown in the equality paradox in Section \ref{sec:Structural-Aggregation}.
The core of this paradox is that $\beta$ is uncorrelated with $\alpha$
and decreasing $\beta$ happens naturally in the evolution of the
economy. So far it is hard to identify the connection between $\alpha$
and $\beta$. Thus the way of drafting a policy to support $\alpha$-growth
is in the mist. Deeper policy that can affect these parameters may
relate to some fundamentals of society such as education and religious
systems. 

An economy of considering only $\beta$-growth could be in danger.
When $\beta$ becomes small enough, the economy will become highly
unequal. In this case, even if the economy scale can reach to a very
high level (by assuming $\alpha$ not decreasing), the structure is
fragile because the majority concentrate on one side and they are
far from the elites whose fraction of the population decreases with
$\beta$. This type of social structure often leads to a revolution,
a war, or a creative destruction. Should we need to resist a $\beta$-growth?
If yes, then how?

\section{Conclusive Remarks}

A new perspective, new laws, and new models are proposed. Meanwhile,
a new paradox, new challengings and new expectations are initiated.
Lights are shining onto an unfamiliar and perhaps treacherous direction
that however can not and should not be concealed. Effective work in
this direction necessarily calls for a shift from classical devices
to novel ones. New theories deal simultaneously with heterogeneity,
endogeneity, and dynamics. New interpretation considers to treat the
nature, the forces and the human factors that determine the characteristics
and trends of our evolution as a whole. New machineries integrate
these complexities in reality and are ready for a venture into fields
beyond recognized patterns. A fascinating risk awaits in front. But
one intuition becomes clear. Growth and inequality, two symbiosis
contraries in the names of goodwill and ``evil'', will float with
us towards the brave uncertain new world, that has such people expecting
it.

\bibliographystyle{amsplain}

\newpage{}

\begin{landscape} \begin{table}[!htbp]  \centering   \resizebox{0.6\textwidth}{!}{\begin{minipage}{\textwidth} \center \caption{Results of Reduced Form Estimates}      \begin{tabular}{@{\extracolsep{5pt}}lcccc} \label{Tab:RF}  \\[-1.8ex]\hline  \hline \\[-1.8ex]   & \multicolumn{4}{c}{\textit{Dependent variable:}} \\  \cline{2-5}  \\[-1.8ex] & \multicolumn{2}{c}{$m_{Y}(t+1)-m_{Y}(t)$} & \multicolumn{2}{c}{$\sigma_{Y}^{2}(t+1)-\sigma_{Y}^{2}(t)$} \\  \\[-1.8ex] & Growth \eqref{eq:rf-1} & Growth rate \eqref{eq:rf-1} & Growth  \eqref{eq:rf-2} & Growth rate \eqref{eq:rf-2}\\  \hline \\[-1.8ex]   Coef-1 $-$0.010 & $-$0.016$^{**}$ &  &  \\     & (0.007) & (0.007) &  &  \\     & & & & \\    Coef-2 & $-$0.001$^{***}$ & $-$40.777$^{***}$ &  &  \\     & (0.0002) & (8.778) &  &  \\     & & & & \\    Coef-3 &  &  & $-$0.00002 & $-$0.00000 \\    &  &  & (0.00005) & (0.00000) \\   & & & & \\    Coef-4 &  &  & 0.772$^{***}$ & 0.780$^{***}$ \\     &  &  & (0.107) & (0.107) \\     & & & & \\    Constant & 665.123$^{**}$ & 0.176$^{**}$ & $-$7,306.789 & 0.00004 \\     & (306.751) & (0.070) & (99,917.010) & (0.0004) \\     & & & & \\  \hline \\[-1.8ex]   Observations & 87 & 87 & 86 & 86 \\   R$^{2}$ & 0.281 & 0.254 & 0.386 & 0.390 \\   Adjusted R$^{2}$ & 0.264 & 0.236 & 0.371 & 0.375 \\   Residual Std. Error & 240.671 (df = 84) & 0.005 (df = 84) & 145,910.300 (df = 83) & 0.0001 (df = 83) \\   F Statistic & 16.404$^{***}$ (df = 2; 84) & 14.308$^{***}$ (df = 2; 84) & 26.036$^{***}$ (df = 2; 83) & 26.535$^{***}$ (df = 2; 83) \\   \hline   \hline \\[-1.8ex]   \textit{Note:}  & \multicolumn{4}{r}{$^{*}$p$<$0.1; $^{**}$p$<$0.05; $^{***}$p$<$0.01} \\   \end{tabular}  \end{minipage}}  \end{table} \end{landscape}

 \begin{table}[!htbp] \centering   \resizebox{0.6\textwidth}{!}{\begin{minipage}{\textwidth} \center \caption{Results of Structural Form Estimates}  \label{Tab:SF}    \begin{tabular}{@{\extracolsep{5pt}}lcccc}  \\[-1.8ex]\hline  \hline \\[-1.8ex]   & \multicolumn{4}{c}{\textit{Dependent variable:}} \\  \cline{2-5}  \\[-1.8ex] & \multicolumn{2}{c}{$\alpha_t/\beta_t$} & Residual of & Residual of \\  \\[-1.8ex] & pre-analysis & \eqref{eq:sf-1} & \eqref{eq:sf-1} & \eqref{eq:sf-2}\\  \hline \\[-1.8ex]   $m_Y(t)$ & 0.0003$^{***}$ &  &  &  \\    & (0.00002) &  &  &  \\    & & & & \\   Coef-1 &  & $-$35.695$^{**}$ &  &  \\    &  & (15.500) &  &  \\    & & & & \\   $\alpha_t/\beta_t$ &  &  & 0.790$^{***}$ & $-$0.002 \\    &  &  & (0.091) & (0.002) \\    & & & & \\   Constant & $-$3.891$^{***}$ & 11.916$^{***}$ & $-$8.999$^{***}$ & 0.023 \\    & (0.712) & (0.348) & (1.043) & (0.027) \\    & & & & \\  \hline \\[-1.8ex]  Observations & 22 & 22 & 22 & 22 \\  R$^{2}$ & 0.959 & 0.210 & 0.790 & 0.051 \\  Adjusted R$^{2}$ & 0.957 & 0.170 & 0.780 & 0.003 \\  Residual Std. Error (df = 20) & 0.279 & 1.220 & 0.559 & 0.015 \\  F Statistic (df = 1; 20) & 464.048$^{***}$ & 5.303$^{**}$ & 75.425$^{***}$ & 1.071 \\  \hline  \hline \\[-1.8ex]  \textit{Note:}  & \multicolumn{4}{r}{$^{*}$p$<$0.1; $^{**}$p$<$0.05; $^{***}$p$<$0.01} \\  \end{tabular}       \end{minipage}} \end{table} 

\newpage{}

\appendix

\section{Proofs}

\subsection{Proof of Lemma \ref{lem:Markov}}
\begin{proof}
We know that $X_{t}(\cdot)=f^{(a_{t})}f^{(a_{t-1})}\cdots f^{(a_{1})}(\cdot)$
and that $(a_{1},\dots,a_{t})$ are i.i.d.. $X_{t}$ and $a_{t}$
generate the same $\sigma$-algebra $\mathcal{B}_{t}$. Then for any
$a'_{t-2}\in\mathcal{B}_{t-2}$, $a_{t-2}^{'}$ and $a_{t}$ are independent,
so 
\[
\Pr\left\{ X_{t}\in A\left|a'_{t-2},\: X_{t-1}\right.\right\} =\Pr\left\{ f^{(a_{t})}(X_{t-1})\in A\left|a'_{t-2},\: X_{t-1}\right.\right\} =\Pr\left\{ f^{(a_{t})}(X_{t-1})\in A\left|X_{t-1}\right.\right\} 
\]
which is the Markov property $\Pr\{X_{t}\in A|\mathcal{B}_{t-1},\mathcal{B}_{t-2}\}=\Pr\{X_{t}\in A|\mathcal{B}_{t-1}\}$.
The result follows.
\end{proof}

\subsection{Proof of Lemma \ref{lem:Convergence}}
\begin{proof}
Assumption \ref{assu:growth_2} implies that we can select a subsequence
$\{a_{t}^{'}\}\subset\{a_{t}\}$ so that $\tilde{X}_{t}(\cdot)=f^{(a_{t}^{'})}(\cdot)$
and $\tilde{X}_{t}(\cdot)$ is monotone. We have $\tilde{X}_{t}(\underline{\Omega})\leq\tilde{X}_{t+1}(\underline{\Omega})$
almost surely. Then given $y$, we have $\Pr(\tilde{X}_{t}(\underline{\Omega})\leq y)\geq\Pr(\tilde{X}_{t+1}(\underline{\Omega})\leq y)$.
Since given $\underline{\Omega}<y<\overline{\Omega}$, $\Pr(\tilde{X}_{t}(\underline{\Omega})\leq y)$
is monotone in $t$ and $\Pr(\tilde{X}_{t}(\underline{\Omega})\leq y)$
is bounded between $0$ and $1$, we know that for each $y$, $\Pr(\tilde{X}_{t}(\underline{\Omega})\leq y)$
must converge as $t\rightarrow\infty$.%
\footnote{Monotone sequence on a bounded domain always converges.%
} 

Similar argument holds for $\Pr(\tilde{X}_{t}(\overline{\Omega})\leq y)$.
As $y<\overline{\Omega}$, it is obvious that $\Pr(\overline{\Omega}\leq y)=0$.
Thus we have $\Pr(\tilde{X}_{t}(\overline{\Omega})\leq y)\leq\Pr(\tilde{X}_{t+1}(\overline{\Omega})\leq y)$.
Monotonicity and boundedness tell that $\Pr(\tilde{X}_{t}(\overline{\Omega})\leq y)$
must converge as $t\rightarrow\infty$. 

Because the index set $\mathcal{A}$ only contains i.i.d. random variables,
any other sequence $\{a_{t}\}$ has the same distribution as $\{a_{t}^{'}\}$,
then any $X_{t}(\cdot)=f^{(a_{t})}(\cdot)$ will induce the same convergent
results for $\Pr(X_{t}(\underline{\Omega})\leq y)$ and $\Pr(X_{t}(\overline{\Omega})\leq y)$.
\end{proof}

\subsection{Proof of Proposition \ref{prop:stationary}}
\begin{proof}
Lemma \ref{lem:Markov} shows that $X_{t}(\omega)$ is a Markov process.
For a Markov process, let $P(x,\mathcal{X})=\Pr(f^{(a)}(x)\in\mathcal{X})$
be a transition probability where $a$ has the same distribution as
$a_{1}$. By definition, $X_{t}(\cdot)=\overset{t-\mbox{times}}{\overbrace{f^{(a)}f^{(a)}\cdots f^{(a)}}}(\cdot)$.
We have an expression of a transition probability for $t$-times 
\[
P^{(t)}(X_{0}=\underline{\Omega},\mathcal{X})=\int_{\Omega}P(x,\mathcal{X})dP^{(t-1)}(X_{0}=\underline{\Omega},x).
\]
As $P(x,\mathcal{X})$ is a bounded, continuous and real-valued functions
on $x$, by Helly-Bray theorem ($P^{(t)}\rightarrow\mathbb{P}$) and
Lemma \ref{lem:Convergence}, we have 
\[
\mathbb{P}\left(X^{*}(\underline{\Omega})\in\mathcal{X}\right)=\int_{\Omega}P(x,\mathcal{X})d\mathbb{P}\left(X^{*}(\underline{\Omega})=x\right)
\]
as $t\rightarrow\infty$. This means that $\mathbb{P}\left(X^{*}(\underline{\Omega})\in\mathcal{X}\right)$
is a stationary distribution. Similar argument holds for $\mathbb{P}\left(X^{*}(\overline{\Omega})\in\mathcal{X}\right)$.
\end{proof}

\subsection{Proof of Theorem \ref{thm:UniqueStationary}}
\begin{proof}
Given a subsequence $\{a_{t}^{'}\}\subset\{a_{t}\}$, $\tilde{X}_{t}(\cdot)=f(\tilde{X}_{t-1}(\cdot))$
and $f$ is monotone, Proposition \ref{prop:stationary} shows that
$\Pr(\tilde{X}_{t}(\underline{\Omega})\in\mathcal{X})$ and $\Pr(\tilde{X}_{t}(\overline{\Omega})\in\mathcal{X})$
converge. As $\Pr(\tilde{X}_{t}(\underline{\Omega})\in\mathcal{X})\geq\Pr(\tilde{X}_{t}(\underline{\Omega})\in\mathcal{X})$,
we have
\[
\mathbb{P}\left(X^{*}(\underline{\Omega})\in\mathcal{X}\right)\geq\mathbb{P}\left(X^{*}(\overline{\Omega})\in\mathcal{X}\right).
\]
We now prove $\mathbb{P}\left(X^{*}(\underline{\Omega})\in\mathcal{X}\right)\leq\mathbb{P}\left(X^{*}(\overline{\Omega})\in\mathcal{X}\right)$.

With Assumption \ref{assu:Equality-prob}, let $\tau$ denote the
time when $X_{\tau}(\underline{\Omega})$ and $X_{\tau}(\overline{\Omega})$
will enter two separated sets. Let $\tau=\inf\left\{ t:\: X_{t}(\underline{\Omega})\in\mathcal{X},\; X_{t}(\overline{\Omega})\in\Omega\backslash\mathcal{X}\right\} $.
This hitting time $\tau$ is random and $\Pr(\tau<\infty)=1$. Then
use the transition probability $P(x,\mathcal{X})=\Pr(f^{(a)}(x)\in\mathcal{X})$
and stationary property, we have
\begin{align*}
\Pr\left(X_{t}(\underline{\Omega})\in\mathcal{X},\: t>\tau\right)\overset{(a)}{=} & \sum_{j=1}^{t-1}\left\{ \int P(x,\mathcal{X})dP(X_{t-j}(\underline{\Omega}),x|\tau=j)\right\} \Pr(\tau=j)\\
\overset{(b)}{\leq} & \sum_{j=1}^{t-1}\left\{ \int P(x,\mathcal{X})dP(X_{t-j}(\underline{\Omega}),x|\tau=j)\right\} ,
\end{align*}
where $\overset{(a)}{=}$ comes from the iterative expression $P(X_{2}(\underline{\Omega}),\mathcal{X})=\int_{\Omega}P(x,\mathcal{X})dP(X_{1}(\underline{\Omega}),x)$
and the conditioning of $\tau$, $\overset{(b)}{\leq}$ comes from
$\Pr(\tau=j)<1$ for some $j$s. 

Similar approach is applied to $X_{t}(\overline{\Omega})$:
\begin{align*}
\Pr\left(X_{t}(\overline{\Omega})\in\mathcal{X},\: t>\tau\right)= & \sum_{j=1}^{t-1}\left\{ \int P(x,\mathcal{X})dP(X_{t-j}(\overline{\Omega}),x|\tau=j)\right\} \Pr(\tau=j)\\
\overset{(c)}{\geq} & \sum_{j=1}^{t-1}\left\{ \int P(x,\mathcal{X})dP(X_{t-j}(\overline{\Omega}),x|\tau=j)\right\} ,
\end{align*}
where $\overset{(c)}{\geq}$ comes from the fact that when $t>\tau$,
$X_{t}(\overline{\Omega})\in\Omega\backslash\mathcal{X}$ so $\Pr(\tau=j)=0$.
Because the subset $\mathcal{X}$ is arbitrarily chosen, we can choose
the set ``closer'' to $\overline{\Omega}$ so that $P(X_{t}(\overline{\Omega}),\mathcal{X})\geq P(X_{t}(\underline{\Omega}),\mathcal{X})$
when $t\rightarrow\infty$. Combine this result with the previous
two inequalities and take $t\rightarrow\infty$ so that $\Pr\{t>\tau\}=1$,
we have
\begin{align*}
\Pr\left(X^{*}(\underline{\Omega})\in\mathcal{X}\right)\leq\sum_{j=1}^{\infty}\left\{ \int P(x,\mathcal{X})dP(X_{t-j}(\underline{\Omega}),x)\right\}  & \leq\sum_{j=1}^{\infty}\left\{ \int P(x,\mathcal{X})dP(X_{t-j}(\overline{\Omega}),x)\right\} \\
\leq & \Pr\left(X^{*}(\overline{\Omega})\in\mathcal{X}\right)
\end{align*}
so $\mathbb{P}\left(X^{*}(\underline{\Omega})\in\mathcal{X}\right)\leq\mathbb{P}\left(X^{*}(\overline{\Omega})\in\mathcal{X}\right)$,
the result follows.
\end{proof}

\subsection{Proof of Theorem \ref{thm:inifiniteDivisible}}
\begin{proof}
By Theorem \ref{thm:UniqueStationary}, we know $X^{*}(\omega^{1})\overset{d}{=}X^{*}(\omega^{2})\cdots\overset{d}{=}X^{*}(\omega^{N})$.
$\overset{d}{=}$ means equivalence in distribution. A simple arrangement
\[
X^{*}(\omega^{1})+X^{*}(\omega^{2})-X^{*}(\omega^{2})\overset{d}{=}X^{*}(\omega^{1})
\]
does not violate the equality in distribution. The distributional
expression of above equality is a convolution of $X^{*}(\omega^{1})+X^{*}(\omega^{2})$
and $-X^{*}(\omega^{2})$ such that
\[
\int_{\Omega}\mathbb{P}_{X^{*}(\omega^{1})+X^{*}(\omega^{2})}\left(X^{*}(\omega^{1})+x\right)\mathbb{P}_{X^{*}(\omega^{2})}(-x)dx=\mathbb{P}(X^{*}(\omega^{1}))
\]
$\mathbb{P}_{X^{*}(\omega^{1})+X^{*}(\omega^{2})}$ and $\mathbb{P}_{X^{*}(\omega^{2})}$
stand for the stationary distribution functions of $X^{*}(\omega^{1})+X^{*}(\omega^{2})$
and $X^{*}(\omega^{2})$. Without loss of generality, we just consider
$\Omega\in\mathbb{R}_{+}^{d}\cup\{0\}$. Note that $\mathbb{P}_{X^{*}(\omega^{2})}(-x)=0$
for $\Omega\in\mathbb{R}_{+}^{d}$ and $\mathbb{P}_{X^{*}(\omega^{2})}(0)\geq0$
by the property of production functions. Then the convolution degenerates
to
\[
\mathbb{P}_{X^{*}(\omega^{1})+X^{*}(\omega^{2})}\left(X^{*}(\omega^{1})\right)=\mathbb{P}(X^{*}(\omega^{1})).
\]
So $X^{*}(\omega^{1})+X^{*}(\omega^{2})\sim\mathbb{P}$. This argument
can be extended to a sum of all $X^{*}(\omega^{i})$. The result follows.
\end{proof}

\subsection{Proof of Theorem \ref{thm:ReducedFormAgg} }
\begin{proof}
Given the mean field equation (\ref{eq:mean_field}), one expands
$a_{1}(y)$ around $m_{Y}(t)$ via Taylor series
\begin{align*}
\int a_{1}(y)\mathbb{Q}(y,t)dy=a_{1}(m_{Y}(t))+ & \overset{(*)}{\overbrace{\left[\int(y-m_{Y}(t))\mathbb{Q}(y,t)dy\right]\left.\frac{\partial a_{1}(m_{Y}(t))}{\partial y}\right|_{y=m_{Y}(t)}}}\\
+\frac{1}{2}\left[\int(y-m_{Y}(t))^{2}\mathbb{Q}(y,t)dy\right] & \left.\frac{\partial^{2}a_{1}(m_{Y}(t))}{\partial y^{2}}\right|_{y=m_{Y}(t)}+\dots
\end{align*}
where $a_{1}(m_{Y}(t))=\int(y'-m_{Y}(t))\mathbb{W}(y'|y)dy'$. By
the definition $\int y\mathbb{Q}(y,t)dy=m_{Y}(t)$, we know $(*)$
is zero. Thus the expression becomes
\[
\int a_{1}(y)\mathbb{Q}(y,t)dy=a_{1}(m_{Y}(t))+\frac{1}{2}\overset{=\sigma_{Y}^{2}(t)}{\overbrace{\left[\int(y-m_{Y}(t))^{2}\mathbb{Q}(y,t)dy\right]}}a_{1}^{(2)}(m_{Y}(t))+\dots
\]
where $a_{1}^{(2)}(m_{Y}(t))$ denotes $\partial^{2}a_{1}(m_{Y}(t))/\partial y^{2}|_{y=m_{Y}(t)}$
and $\sigma_{Y}^{2}(t)$ is variance of $m_{Y}(t)$. Similarly, one
can deduce the second moment of the mean field
\begin{align*}
\int y^{2}\frac{\partial\mathbb{Q}(y,t)}{\partial t}dy & =\int\int y^{2}\left[\mathbb{W}(y|y')\mathbb{Q}(y',t)-\mathbb{W}(y'|y)\mathbb{Q}(y,t)\right]dydy'\\
= & \int\int\left(y'^{2}-y^{2}\right)\mathbb{W}(y'|y)\mathbb{Q}(y,t)dydy'\\
= & \int\int\left[\left(y'-y\right)^{2}+2y(y'-y)\right]\mathbb{W}(y'|y)\mathbb{Q}(y,t)dydy'\\
= & \int a_{2}(y)\mathbb{Q}(y,t)dy+2\int y\, a_{1}(y)\mathbb{Q}(y,t)dy.
\end{align*}
As $\sigma_{Y}^{2}(t)=\int(y-m_{Y}(t))^{2}\mathbb{Q}(y,t)dy=\int y^{2}\mathbb{Q}(y,t)dy-m_{Y}^{2}(t)$,
there is 
\begin{align*}
\frac{d\sigma_{Y}^{2}(t)}{dt} & =\int y^{2}\frac{\partial\mathbb{Q}(y,t)}{\partial t}dy-2m_{Y}(t)\frac{dm_{Y}(t)}{dt}\\
 & =\int a_{2}(y)\mathbb{Q}(y,t)dy+2\int y\, a_{1}(y)\mathbb{Q}(y,t)dy-2m_{Y}(t)\int a_{1}(y)\mathbb{Q}(y,t)dy\\
 & =\int a_{2}(y)\mathbb{Q}(y,t)dy+2\int(y-m_{Y}(t))\, a_{1}(y)\mathbb{Q}(y,t)dy\\
 & \overset{(a)}{=}\int a_{2}(y)\mathbb{Q}(y,t)dy+2\int\Biggl\{(y-m_{Y}(t))\, a_{1}(m_{Y}(t))\\
 & +\left.\left.\frac{\partial a_{1}(m_{Y}(t))}{\partial y}\right|_{y=m_{Y}(t)}\left(y-m_{Y}(t)\right)^{2}+\dots\right\} \mathbb{Q}(y,t)dy\\
 & \overset{(b)}{=}a_{2}(m_{Y}(t))+2a_{1}^{(1)}(m_{Y}(t))\sigma_{Y}^{2}(t)+\dots
\end{align*}
where $\overset{(a)}{=}$ comes from a Taylor expansion of $a_{1}(y)$
around $m_{Y}(t)$ and $\overset{(b)}{=}$ comes from another Taylor
expansion of $a_{2}(y)$ around $m_{Y}(t)$. Expansion terms that
outvie $(y-m_{Y}(t))^{2}$ have not displayed in the expression. 
\end{proof}

\subsection{Proof of Theorem \ref{thm:StructuralAgg}}
\begin{proof}
A heuristic proof of (I) and (II). 

(I) By Assumption \ref{assu:RateAndSize}, we consider two states
transition at small time interval, a monotone growth from $x'$ to
$x=1$ or to stay at position $x=0$:
\begin{align*}
\mathbb{P}_{\Delta t}(x=1|x') & =\alpha\Delta t+o((\Delta t)^{2}),\\
\mathbb{P}_{\Delta t}(x=0|x') & =1-\alpha\Delta t+o((\Delta t)^{2}).
\end{align*}
Let $\mathbb{P}(x_{1},t)=\Pr\{X_{t}(\omega)=x_{1}\}$. Then 
\begin{align*}
\mathbb{P}(x_{1},t+\Delta t)\overset{(a)}{=} & \mathbb{P}(x_{1},t)\Pr\{X_{\Delta t}(\omega)=0|x_{1}\}+\mathbb{P}(x_{1},t)\Pr\{X_{\Delta t}(\omega)=1|x_{1}\}\\
= & \mathbb{P}(x_{1},t)(1-\alpha\Delta t)+\mathbb{P}(x_{1},t)\times\alpha\Delta t+o((\Delta t)^{2})\\
\overset{(b)}{=} & \mathbb{P}(x_{1},t)(1-\alpha\Delta t)+\int\mathbb{W}(x_{1}|x')\mathbb{P}(x',t)dx'\times\alpha\Delta t+o((\Delta t)^{2})
\end{align*}
where $\overset{(a)}{=}$ comes from transition from $\mathbb{P}(x_{1},t)$
to $\mathbb{P}(x_{1},t+\Delta t)$ and $\overset{(b)}{=}$ comes from
\[
\mathbb{P}(x_{1},t)\overset{(c)}{=}\lim_{\Delta t\rightarrow0}\int\left(\mathbb{P}_{\Delta t}(x_{1}|x')\mathbb{P}(x',t)\right)dx'\overset{(d)}{=}\int\mathbb{W}(x_{1}|x')\mathbb{P}(x',t)dx'.
\]
$\overset{(c)}{=}$ comes from Chapman-Kolmogorov equation (see (\ref{eq:master_3}))
and $\overset{(d)}{=}$ uses the instantaneous rate $\mathbb{W}(x_{1}|x')$
in (\ref{eq:master_1}) to approximate the transition probability
when $\Delta t\rightarrow0$. 

Arrange the expression of $\Pr\{X_{t+\Delta t}(\omega)=1\}$ and take
the limit of $\Delta t$,
\begin{align*}
\lim_{\Delta t\rightarrow0} & \frac{\mathbb{P}(x_{1},t+\Delta t)-\mathbb{P}(x_{1},t)}{\Delta t}\\
 & =\alpha\left(\int\mathbb{W}(x_{1}|x')\mathbb{P}(x',t)dx'-\mathbb{P}(x_{1},t)\right)+\lim_{\Delta t\rightarrow0}o(\Delta t)\\
\frac{d}{dt}\mathbb{P}(x_{1},t) & =\alpha\left(\int\mathbb{W}(x_{1}|x')\mathbb{P}(x',t)dx'-\mathbb{P}(x_{1},t)\right).
\end{align*}
Assumption \ref{assu:RateAndSize} imposes the rate of size change,
$c$, to be linear in $t$. By change of variable $x=x_{1}+ct$, one
has
\[
\frac{\partial}{\partial t}\mathbb{P}(x,t)+c\frac{\partial}{\partial x}\mathbb{P}(x,t)=\alpha\left(\int\mathbb{W}(x|x')\mathbb{P}(x',t)dx'-\mathbb{P}(x,t)\right).
\]
The result follows.%
\footnote{An alternative way of attaining this expression is to use Poisson
semigroup.%
}

(II) As the right hand side of (\ref{eq:StationaryGamma}) is a convolution
integral, it can be compactly represented as
\[
\left[\int_{0}^{x}udu\right]\frac{\partial\mathbb{P}(x)}{\partial x}=\alpha\left[\mathbb{V}*\mathbb{P}\right](x).
\]
On the right hand side of the equation, Laplace transform of $[\mathbb{V}*\mathbb{P}](x)$
is $\mathbb{V}_{\mathcal{L}}(s)\mathbb{P}_{\mathcal{L}}(s)$ where
$\mathbb{V}_{\mathcal{L}}(s)$ and $\mathbb{P}_{\mathcal{L}}(s)$
are the Laplace transform of $\mathbb{V}(x)$ and $\mathbb{P}(x)$
respectively. Note that $\mathbb{V}_{\mathcal{L}}(s)=\beta/(s+\beta)$
by the property of Laplace transform of exponential distribution. 

On the left hand side of this equation, Laplace transform is 
\[
\int\left\{ \left[\int udu\right]\frac{\partial\mathbb{P}(x)}{\partial x}e^{-sx-su}\right\} dx=\int_{0}^{\infty}\left(\int_{0}^{\infty}e^{-su}du\right)\left(x\frac{\partial\mathbb{P}(x)}{\partial x}e^{-sx}\right)dx.
\]
Since $-\int_{0}^{\infty}e^{-su}du=1/s$, after simplification, the
above equation becomes 
\[
-\int\frac{x}{s}\frac{\partial\mathbb{P}(x)}{\partial x}e^{-sx}dx=\int\frac{1}{s}\frac{\partial\mathbb{P}(x)}{\partial x}\frac{\partial e^{-sx}}{\partial s}dx=\frac{\partial\left\{ \frac{1}{s}\int\left(\frac{\partial\mathbb{P}(x)}{\partial x}e^{-sx}dx\right)\right\} }{\partial s}\overset{(a)}{=}-\frac{\partial\mathbb{P}_{\mathcal{L}}(s)}{\partial s},
\]
where $\overset{(a)}{=}$ comes from differentiation by parts: 
\[
-\frac{1}{s}\int_{0}^{\infty}\left(\frac{\partial\mathbb{P}(x)}{\partial x}\right)e^{-sx}dx=\int_{0}^{\infty}\mathbb{P}(x)e^{-sx}dx-\left[\mathbb{P}(x)e^{-sx}\right]_{0}^{\infty}=\mathbb{P}_{\mathcal{L}}(s),
\]
So the equation becomes 
\[
\frac{\partial\mathbb{P}_{\mathcal{L}}(s)}{\partial s}=-\frac{\alpha\beta\mathbb{P}_{\mathcal{L}}(s)}{s+\beta}.
\]
Dividing $\mathbb{P}_{\mathcal{L}}(s)$ on both sides gives 
\[
\frac{1}{\mathbb{P}_{\mathcal{L}}(s)}\frac{\partial\mathbb{P}_{\mathcal{L}}(s)}{\partial s}=\frac{\partial\ln(\mathbb{P}_{\mathcal{L}}(s))}{\partial s}=-\alpha\frac{\partial\ln((s+\beta)/\beta)}{\partial s}.
\]
Taking the integral w.r.t. $s$ gives 
\[
\ln(\mathbb{P}_{\mathcal{L}}(s))=-\alpha\ln\left(\frac{s+\beta}{\beta}\right).
\]
The solution is $\mathbb{P}_{\mathcal{L}}(s)=\beta^{\alpha}[s+\beta]^{-\alpha}$
which is the Laplace transform of Gamma distribution. To see this,
assume $X$ follows Gamma distribution, then its Laplace transform
is
\[
\mathbb{E}\left[e^{-sX}\right]=\int e^{-sx}\frac{\beta^{\alpha}}{\Gamma(\alpha)}x^{\alpha-1}e^{-\beta x}dx=\frac{\beta^{\alpha}}{\Gamma(\alpha)(\beta+s)^{\alpha}}\underbrace{\int\left((\beta+s)x\right)^{\alpha-1}e^{-(\beta+s)x}d\left((\beta+s)x\right)}_{(*)}.
\]
Substitute $u=(\beta+s)x$ into the $(*)$ term, it becomes $\Gamma(\alpha)$.
So $\mathbb{P}_{\mathcal{L}}(s)=\beta^{\alpha}[s+\beta]^{-\alpha}$.
The result follows.
\end{proof}

\subsection{Proof of Corollary \ref{cor:PowerLaw}}
\begin{proof}
By infinite divisible distribution from Theorem \ref{thm:inifiniteDivisible},
one knows 
\[
\mathbb{P}(bx)=k(b)\mathbb{P}(x)
\]
for any integer $b=1,2,3,\dots$, as $bx\sim\mathbb{P}$ and $x\sim\mathbb{P}$,
$k(\cdot)$ is a normalized function. Set $x=1$, $\mathbb{P}(b)=k(b)\mathbb{P}(1)$,
one knows $k(b)=\mathbb{P}(b)/\mathbb{P}(1)$. So $\mathbb{P}(bx)=\mathbb{P}(b)\mathbb{P}(x)/\mathbb{P}(1)$.
By taking the differentiation w.r.t. $b$ on both side, one has $x(\dot{\mathbb{P}}(bx))=\mathbb{P}(x)\dot{\mathbb{P}}(b)/\mathbb{P}(1)$
where $\dot{\mathbb{P}}$ stands for differentiation w.r.t. its argument.
Set $b=1$,
\[
x\frac{d\mathbb{P}(x)}{dx}=\frac{\dot{\mathbb{P}}(1)}{\mathbb{P}(1)}\mathbb{P}(x)=\alpha\mathbb{P}(x)
\]
one has (\ref{eq:PowerLaw}) and $\alpha=\dot{\mathbb{P}}(1)/\mathbb{P}(1)$.
Integrate above equation
\begin{equation}
\ln\mathbb{P}(x)=\frac{\mathbb{P}(1)}{\dot{\mathbb{P}}(1)}\ln x+\mbox{Constant}.\label{eq:Powerlaw-1}
\end{equation}
One has $\mbox{Constant}=\ln\mathbb{P}(1)$ by setting $x=1$. Taking
the exponential of both sides in (\ref{eq:Powerlaw-1}) gives $\mathbb{P}(x)=\mathbb{P}(1)x^{\alpha}.$
The results follows. 
\end{proof}

\subsection{Proof of Theorem \ref{thm:StructuralForm}}
\begin{proof}
The first equation in (\ref{eq:StructuralEq}) comes from Assumption
\ref{assu:AggOneToOne} that $g(\cdot)$ is one-to-one. Thus by change
of measures, we have 
\[
m_{Y}(t)=\int y\mathbb{Q}(y,t)dy=\int\left[g(x)\frac{\mathbb{Q}(y,t)dy}{\mathbb{P}(x,t)dx}\right]\mathbb{P}(x,t)dx.
\]
Thus $\mathbb{E}_{t}\left[g\left(X_{t}(\omega)\right)\mathbb{L}(t)\right]=m_{Y}(t)$.

To show the second equation in (\ref{eq:StructuralEq}), we need to
represent the mean field equation (\ref{eq:mean_field}): 
\begin{align*}
\frac{dm_{Y}(t)}{dt}= & \int\int(y'-y)\mathbb{W}(y'|y)\mathbb{Q}(y,t)dy'dy\\
= & \int\int y'\mathbb{W}(y'|y)\mathbb{Q}(y,t)dy'dy-\int\int y\mathbb{W}(y'|y)\mathbb{Q}(y,t)dy'dy\\
\overset{(a)}{=} & \int\int y'\mathbb{W}(y'|y)\mathbb{Q}(y,t)dy'dy-\int\int y(A-\mathbb{W}(x'|x))\mathbb{Q}(y,t)dy'dy+o(1)\\
\overset{(b)}{=} & \int\int\left(y'\mathbb{W}(y'|y)-yA\right)\mathbb{Q}(y,t)dy'dy+\left\{ \int y\mathbb{Q}(y,t)dy\right\} \left\{ \int\mathbb{W}(x'|x)dg(x')\right\} +o(1)\\
= & \theta_{t}\, m_{Y}(t)+e_{t}+o(1)
\end{align*}
where $\overset{(a)}{=}$ uses the transition rate $\mathbb{W}(x'|x)$
in Theorem \ref{thm:StructuralAgg}, in $\overset{(b)}{=}$ we replace
the aggregation $y'$ with $g\left(x'\right)$. Let $e_{t}=\int\int\left(y'\mathbb{W}(y'|y)-yA\right)\mathbb{Q}(y,t)dy'dy$,
$A=\sup\left|\mathbb{W}(y'|y)+\mathbb{W}(x'|x)\right|$ is the sup-norm
of two transition kernels, $\theta_{t}=\int\mathbb{W}(x'|X_{t})dg(x')$
uses the kernel conditional on the state value $x=X_{t}(\omega)$,
$m_{Y}(t)=\int y\mathbb{Q}(y,t)dy$. The result follows.

The third expression in (\ref{eq:StructuralEq}) is a direct result
of Theorem \ref{thm:inifiniteDivisible} and \ref{thm:StructuralAgg}.
\end{proof}

\section{Miscellaneous}

\subsection{\label{sub:Solow-Model}Solow Model}

Provided that the quasi-linear utility of individual $i$ is 
\[
U(K_{t}^{i},L_{t}^{i})=f_{P}(K_{t}^{i},L_{t}^{i})-B_{S,t}^{i}K_{t}^{i}-A_{S,t}^{i}L_{t}^{i},
\]
the rational household maximizes the utility with the following FOC
\[
B_{S,t}^{i}=f_{P}^{'}(X_{t}^{i}),\qquad A_{S,t}^{i}=f_{P}(X_{t}^{i})-f_{P}^{'}(X_{t}^{i})X_{t}^{i},
\]
where $f_{P}^{'}$ is the derivative of $f_{P}$ w.r.t. $X_{t}^{i}=K_{t}^{i}/L_{t}^{i}$.
In this expression, $B_{S,t}^{i}$ and $A_{S,t}^{i}$ are the return
to capital and the cost rate respectively for household $i$. This
FOC gives an affine representations for the production function as
well as the growth function
\begin{align*}
f_{P}(X_{t}^{i})= & A_{S,t}^{i}+B_{S,t}^{i}X_{t}^{i},\\
f(X_{t}^{i})= & \beta_{1}A_{S,t}^{i}+\left(\beta_{1}B_{S,t}^{i}+1-\beta_{2}\right)X_{t}^{i}.
\end{align*}

\subsection{\label{sub:Master-Equation}Derivation of Equation (\ref{eq:MasterEq}) }
\begin{proof}
Given a small time interval $\Delta t$, one expands the transition
probability $\mathbb{Q}_{\Delta t}(y_{3}|y_{2})$ in a Taylor series
over zero
\begin{equation}
\mathbb{Q}_{\Delta t}(y_{3}|y_{2})=\delta(y_{2}-y_{3})+\Delta t\,\mathbb{W}(y_{3}|y_{2})+O((\Delta t)^{2}).\label{eq:master_1}
\end{equation}
Dirac delta function $\delta(y_{2}-y_{3})$ expresses that the probability
to stay at the same state $y_{3}=y_{2}$ equals one for any test function
$\phi$ whereas the probability to change state $y_{3}\neq y_{2}$
equals zero such that $\int_{-\infty}^{\infty}\phi(y)\delta(y)dy=\phi(0)$.
\[
\mathbb{W}(y_{3}|y_{2})=\lim_{h\rightarrow0}\frac{\partial\mathbb{Q}_{h}(y_{3}|y_{2})}{\partial h}
\]
is the time derivative of the transition probability at $h=0$. 

Integrating right hand side of equation (\ref{eq:master_1}) w.r.t.
$y_{3}$ may be larger than one. So we need to normalize $\delta(y_{2}-y_{3})$
by subtracting a factor $\int\mathbb{W}(y_{3}|y_{2})dy_{3}$ so that
\begin{equation}
\mathbb{Q}_{\Delta t}(y_{3}|y_{2})=\left(1-\int\mathbb{W}(y_{3}|y_{2})dy_{3}\right)\delta(y_{2}-y_{3})+\Delta t\,\mathbb{W}(y_{3}|y_{2})+O((\Delta t)^{2}).\label{eq:master_2}
\end{equation}
The Chapman-Kolmogorov equation is
\begin{equation}
\mathbb{Q}_{h+\Delta t}(y_{3}|y_{1})=\int\left[\mathbb{Q}_{\Delta t}(y_{3}|y_{2})\mathbb{Q}_{h}(y_{2}|y_{1})\right]dy_{2}.\label{eq:master_3}
\end{equation}
Substitute (\ref{eq:master_2}) into (\ref{eq:master_3}),
\begin{align*}
\mathbb{Q}_{h+\Delta t}(y_{3}|y_{1})= & \int\mathbb{Q}_{h}(y_{2}|y_{1})\delta(y_{2}-y_{3})dy_{2}-\Delta t\int\int\,\mathbb{W}(y_{3}|y_{2})dy_{3}\delta(y_{2}-y_{3})\mathbb{Q}_{h}(y_{2}|y_{1})dy_{2}\\
 & +\Delta t\,\int\mathbb{W}(y_{3}|y_{2})\mathbb{Q}_{h}(y_{2}|y_{1})dy_{2}+O((\Delta t)^{2})\\
= & \mathbb{Q}_{h}(y_{3}|y_{1})-\Delta t\int\mathbb{W}(y_{2}|y_{3})\mathbb{Q}_{h}(y_{3}|y_{1})dy_{2}+\Delta t\,\int\mathbb{W}(y_{3}|y_{2})\mathbb{Q}_{h}(y_{2}|y_{1})dy_{2}.
\end{align*}
Re-arrange the expression and take $\Delta t\rightarrow0$
\[
\lim_{\Delta t\rightarrow0}\frac{\mathbb{Q}_{h+\Delta t}(y_{3}|y_{1})-\mathbb{Q}_{h}(y_{3}|y_{1})}{\Delta t}=\int\left[\mathbb{W}(y_{3}|y_{2})\mathbb{Q}_{h}(y_{2}|y_{1})-\mathbb{W}(y_{2}|y_{3})\mathbb{Q}_{h}(y_{3}|y_{1})\right]dy_{2}.
\]
The result follows.
\end{proof}

\subsection{\label{sub:Derivation-mean-field}Derivation of Equation (\ref{eq:mean_field})}
\begin{proof}
Take the time derivatives on both sides of $m_{Y}(t)=\int y\mathbb{Q}(y,t)dy$
\begin{align*}
\frac{dm_{Y}(t)}{dt} & =\int y\frac{\partial\mathbb{Q}(y,t)}{\partial t}dy\\
= & \int\int y\left[\mathbb{W}(y|y')\mathbb{Q}(y',t)-\mathbb{W}(y'|y)\mathbb{Q}(y,t)\right]dydy'\\
\overset{(a)}{=} & \int\int(y'-y)\mathbb{W}(y'|y)\mathbb{Q}(y,t)dydy'=\int a_{1}(y)\mathbb{Q}(y,t)dy,
\end{align*}
where $\overset{(a)}{=}$ is from
\[
\int\int y'\mathbb{W}(y'|y)\mathbb{Q}(y,t)dydy'=\int\int y\mathbb{W}(y|y')\mathbb{Q}(y',t)dydy'
\]
by Fubini theorem.
\end{proof}

\subsection{\label{sub:Discussion-of-Memoryless}Discussion of Memoryless Transition
Rate}

The idea of memoryless transitions can be traced back to \cite{Finetti1929a}
who considers 
\[
\Pr\{a>t+\Delta t|a>t\}=\frac{\Pr\{a>t+\Delta t\}}{\Pr\{a>t\}}=\Pr\{a>\Delta t\}.
\]
Taking logarithm and denoting $k(t)=\log\Pr\{a>t\}$, one has $k(t+\Delta t)=k(\Delta t)+k(t)$
a linear decreasing function. Thus $k(t)$ should be a linear function
of $t$ such as $k(t)=-\beta t$ which means $\Pr\{a\leq t\}=1-e^{-\beta t}$,
a exponential distribution whose density kernel is $\beta e^{-\beta t}$.

\subsection{\label{sub:Zeta-Distribution}Zeta Distribution}

The solution of (\ref{eq:PowerLaw}) is $\mathbb{P}(x)=\mathbb{P}(1)x^{-\gamma}$
such that $x=1,2,3,\dots$ and $\mathbb{P}(1)$ is a constant. By
the definition of Riemann Zeta function $\zeta(\gamma)=\sum_{x=1}^{\infty}x^{-\gamma}$,
one knows 
\[
\sum_{x=1}^{\infty}\mathbb{P}(x)=1=\mathbb{P}(1)\sum_{x=1}^{\infty}x^{-\gamma}=\mathbb{P}(1)\zeta(\gamma)=\mathbb{P}(1)\zeta(-\alpha).
\]
Thus $\mathbb{P}(x)=(x^{\gamma}\zeta(\gamma))^{-1}=x^{\alpha}/\zeta(-\alpha)$
which is the zeta distribution, in particular, $\mathbb{P}(1)=1/\zeta(\gamma)$.

\subsection{\label{sub:Scheme-of-Estimates}Scheme of Estimates}

\paragraph*{Reduced estimate:}

\begin{algorithmic} \State Step 1. Regress $m_Y(t+1)-m_Y(t)$ on $m_Y(t)$ for all $t=1,\dots, T$.
\State Step 2. Use residual $\hat{\varepsilon}_{t}$ to construct $\sigma^{2}_{t}$.
\State Step 3. Regress $m_Y(t+1)-m_Y(t)$ on $m_Y(t)$ and $\sigma^{2}_{t}$.
\While {$\hat{\varepsilon}_{t}$ is not a white noise}     \State Update $\sigma^{2}_{t}$ and go to Step 3. \EndWhile     \State Step 4. Regress $\sigma^{2}_{t}- \sigma^{2}_{t-1}$ on $m_Y^{2}(t-1)$ and $\sigma^{2}_{t-1}$ for all $t=1,\dots, T$. \State Step 5.   Heteroskedasticity check of $\hat{\upsilon}_{t}$. \end{algorithmic}In our analysis, the ``while loop'' after Step 3 is not executed
in all cases because white noise patterns are quite significant.

\paragraph*{Structural estimate:}

\begin{algorithmic} \State Step 1. Use income distribution dataset to estimate $\alpha_{t}$ and $\beta_{t}$ for all $t=1,\dots, T$.
\State Step 2. Regress $\alpha_{t}/\beta_{t}$ on $m_{Y}(t)$ for all $t=1,\dots, T$.
\If {$\mbox{Coef}_{1}$ is significant} enter Filtering:
\State Set $\theta(t)=m_Y(t)$ at $t=1$. Recursively compute the mean and variance of $$\theta(t) + \frac{\varepsilon_{2,t}}{\sum_{i=1}^{N}X_{t}(\omega^{i})}$$ for $t=2,\dots,T$.
\Else \State Income data is not qualified for structural estimate. End.
\EndIf \State Step 3. Store filter mean and variance estimates.  \State Step 4. Regress $\alpha_{t}/\beta_{t}$ on filter residuals for all $t=1,\dots, T$ and check endogeneity. \end{algorithmic}

\subsection{\label{sub:Forward-Filter}Forward Filter}

This is an analogy to the forward-algorithm in Kalman filter. Recall
equation (\ref{eq:sf-2}) 
\[
\ln m_{Y}(t+1)-\ln m_{Y}(t)=\overset{(*)}{\overbrace{\theta(t)+\varepsilon_{1,t}}}+\overset{(**)}{\overbrace{\frac{\varepsilon_{2,t}}{\sum_{i=1}^{N}X_{t}(\omega^{i})}}}.
\]
Let $L_{Y}(t)$ denote $\ln m_{Y}(t+1)-\ln m_{Y}(t)$, let $\tilde{\theta}$
denote $(*)$ term, let $\tilde{e}$ denote $(**)$ term so that 
\[
L_{Y}(t)=\tilde{\theta}+\tilde{e}
\]
where $\tilde{\theta}$ is normally distributed with time varying
mean and variance. By the property of Gamma distribution, we know
that $\mathbb{E}[\sum_{i=1}^{N}X_{t}(\omega^{i})]=\alpha_{t}/\beta_{t}$
and $\mbox{Var}[\sum_{i=1}^{N}X_{t}(\omega^{i})]=\alpha_{t}/\beta_{t}^{2}$.
So $\tilde{e}$ is also conditional normal distributed $\tilde{e}|X\sim\mathcal{N}(0,\sigma^{2}(X))$
where $\sigma^{2}(X_{t})=\beta_{t}^{2}/\alpha_{t}$. At $t=0$, let
the initial $\tilde{\theta}\sim\mathcal{N}(\theta_{0},1)$. 

Given observations $\mathcal{M}_{Y}(t):=(L_{Y}(1),\dots L_{Y}(t))$
and conditioning on $X$, we consider
\[
(\tilde{\theta}|\mathcal{M}_{Y}(t),X)\sim\mathcal{N}(\theta_{t},C_{t}|\mathcal{M}_{Y}(t),X).
\]
The distribution of $\tilde{\theta}$ can be updated by Baye's rule
such that
\begin{align*}
\mathcal{N}(\theta_{t},C_{t}|\mathcal{M}_{Y}(t),X) & \propto\prod_{i=1}^{t}\mathcal{N}(L_{Y}(i),\sigma^{2}(X)|\tilde{\theta},X)\mathcal{N}(\tilde{\theta},1)\\
 & =\prod_{i=1}^{t}\frac{1}{\sqrt{2\pi}\sigma(X)}\exp\left\{ -\frac{1}{2\sigma^{2}(X)}(L_{Y}(i)-\tilde{\theta})^{2}\right\} \times\\
 & \frac{1}{\sqrt{2\pi}}\exp\left\{ -\frac{1}{2}(\tilde{\theta}-\theta_{0})^{2}\right\} \\
 & \propto\exp\left\{ -\frac{1}{2\sigma^{2}(X)}\left(\sum_{i=1}^{t}L_{Y}(i)^{2}-2\tilde{\theta}\sum_{i=1}^{t}L_{Y}(i)+t\tilde{\theta}^{2}\right)-\right.\\
 & \left.\frac{1}{2}(\tilde{\theta}^{2}-2\tilde{\theta}\theta_{0}+\theta_{0}^{2})\right\} \\
 & \overset{(a)}{\propto}\exp\left\{ -\frac{1}{2\sigma^{2}(X)}\left(\sum_{i=1}^{t}L_{Y}(i)^{2}-2\tilde{\theta}t\bar{L}_{Y}+t\tilde{\theta}^{2}-\right.\right.\\
 & \left.\left.\sigma^{2}(X)\left(\tilde{\theta}^{2}-2\tilde{\theta}\theta_{0}+\theta_{0}^{2}\right)\right)\right\} \\
 & \propto\exp\left\{ -\frac{1}{2\sigma^{2}(X)}\left((t+\sigma^{2}(X))\tilde{\theta}^{2}-2(t\bar{L}_{Y}+\sigma^{2}(X)\theta_{0})\tilde{\theta}\right)\right\} \\
 & \propto\exp\left\{ -\frac{(t+\sigma^{2}(X))}{2\sigma^{2}(X)}\left(\tilde{\theta}^{2}-2\frac{t\bar{L}_{Y}+\sigma^{2}(X)\theta_{0}}{t+\sigma^{2}(X)}\tilde{\theta}\right)\right\} \\
 & \overset{(b)}{\propto}\exp\left\{ -\frac{(t+\sigma^{2}(X))}{2\sigma^{2}(X)}\left(\tilde{\theta}^{2}-\frac{t\bar{L}_{Y}+\sigma^{2}(X)\theta_{0}}{t+\sigma^{2}(X)}\right)^{2}\right\} 
\end{align*}
where $\bar{L}_{Y}=\sum_{i=1}^{t}L_{Y}(i)/t$. In step $(a)$ because
$\sum_{i=1}^{t}L_{Y}(i)^{2}$ and $\theta_{0}^{2}$ do not relate
to $\tilde{\theta}$, they will not be used to construct updated distribution.
Similar argument is applicable to $(b)$ where the quadratic form
is constructed. Note that $(\tilde{\theta}|\mathcal{M}_{Y}(t),X)\sim\mathcal{N}(\theta_{t},C_{t}^{2})$.
From the last expression, we have 
\[
\theta_{t}=\frac{t\bar{L}_{Y}+\sigma^{2}(X)\theta_{0}}{t+\sigma^{2}(X)},C_{t}=\frac{\sigma^{2}(X)}{t+\sigma^{2}(X)}
\]
Some simplifications give
\[
\theta_{t}=\frac{1}{1+\sigma^{2}(X)/t}\bar{L}_{Y}+\frac{\sigma^{2}(X)/t}{1+\sigma^{2}(X)/t}\theta_{0}
\]
the updated mean is a weighted average between the sample mean $\bar{L}_{Y}$
and the initial state $\theta_{0}$, with weights depending on and
$\ensuremath{\sigma^{2}}(X)$. We can have a recursive expression
\begin{equation}
\theta_{t}=\theta_{t-1}+\frac{C_{t-1}}{C_{t-1}+\sigma^{2}(X)}(L_{Y}(t)-\theta_{t-1}),C_{t}=\left(\frac{1}{\sigma^{2}(X)}+\frac{1}{C_{t-1}}\right)^{-1}\label{eq:KalmanFilter}
\end{equation}
$\theta_{t}$ is obtained by correcting the previous estimate $\theta_{t-1}$.
Once we have the recursive expression, we can update the conditional
variance by setting $\sigma^{2}(X_{t})=\beta_{t}^{2}/\alpha_{t}$
in (\ref{eq:KalmanFilter}).
\end{document}